\documentclass[preprint, 10.5pt]{elsarticle}
\usepackage[english]{babel}
\usepackage{array}
\usepackage{ragged2e}
\usepackage{tabularx}
\usepackage[table]{xcolor}
\newcolumntype{L}[1]{>{\hsize=#1\hsize\raggedright\arraybackslash}X}%
\newcolumntype{R}[1]{>{\hsize=#1\hsize\raggedleft\arraybackslash}X}%
\newcolumntype{C}[1]{>{\hsize=#1\hsize\centering\arraybackslash}X}%
\newcolumntype{Y}{>{\raggedright\arraybackslash}X}
\usepackage{lineno,hyperref}
\usepackage{amssymb}
\usepackage{amsmath}
\usepackage{mathtools, cuted}
\usepackage{lipsum}
\usepackage{gensymb}
\usepackage{color}
\usepackage{graphicx}
\usepackage{wrapfig}
\usepackage{subcaption}
\usepackage{float}
\usepackage{natbib}
\usepackage{geometry}
\journal{Elsevier}

\biboptions{sort&compress}

\begin{document}

\begin{frontmatter}

\title{An Auto-Generated Geometry-Based Discrete Finite Element Model for Damage Evolution in Composite Laminates with Arbitrary Stacking Sequence}

\author[cumae]{Jiakun Liu\corref{cor1}}
\ead{jl2938@cornell.edu}
\author[cumae]{Stuart Leigh Phoenix}
\ead{slp6@cornell.edu}

\cortext[cor1]{Principal corresponding author}

\address[cumae]{Sibley School of Mechanical and Aerospace Engineering, Cornell University, Ithaca, NY 14853, USA}

\begin{abstract}
Stiffness degradation and progressive failure of composite laminates are complex processes involving evolution and multi-mode interactions among fiber fractures, intra-ply matrix cracks and inter-ply delaminations. 
This paper presents a novel finite element model capable of explicitly treating such discrete failures in laminates of random layup. 
Matching of nodes is guaranteed at potential crack bifurcations to ensure correct displacement jumps near crack tips and explicit load transfer among cracks.
The model is entirely geometry-based (no mesh prerequisite) with distinct segments assembled together using surface-based tie constraints, and thus requires no element partitioning or enrichment.
Several numerical examples are included to demonstrate the model's ability to generate results that are in qualitative and quantitative agreement with experimental observations on both damage evolution and tensile strength of specimens.
The present model is believed unique in realizing simultaneous and accurate coupling of all three types of failures in laminates having arbitrary ply angles and layup.
\end{abstract}


\end{frontmatter}


\section{Introduction} \label{sec1:intro}

Damage growth in composite laminates involves: (i) fiber or yarn failures in a ply (seen as transverse cracks across bundles of fibers, yarns or tows under tension, or as kink bands in compression), (ii) intra-ply matrix cracks between fibers and yarns, and (iii) inter-ply delaminations. Besides varying ply angles and stacking sequence, the statistical properties of the fiber and matrix materials, associated with these failure modes and their mutual interaction, together dominate damage initiation, evolution, and ultimate failure of such laminates. 

However, due to the complexity of geometric discontinuities and the multi-mode interacting mechanisms between cracks involved, realistic analytical modeling and numerical analysis of progressive damage in composites has been slow to progress. Though various modeling methods have attempted to account for at least one type of such failure modes, high-fidelity modeling of the initiation and interactive evolution of all three disparate and discrete damage modes remains a challenging problem. 

\begin{figure}[H]
	\centering
	\includegraphics[width=0.8\textwidth]{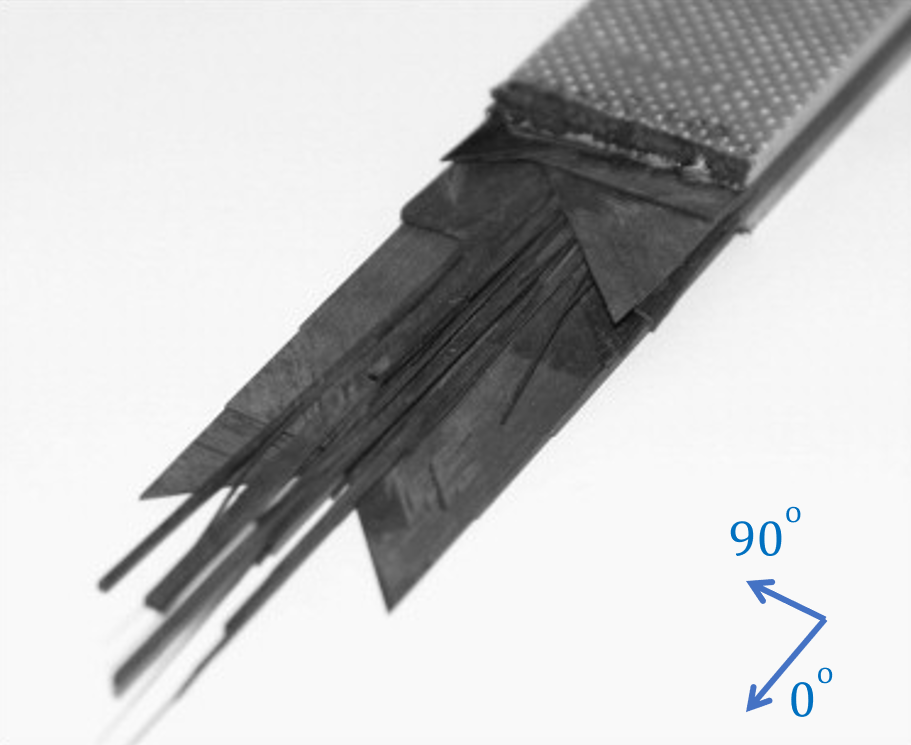}
	\caption{Tensile failure of an un-notched $[45_4/90_4/\text{--}45_4/0_4]_s$ rectangular laminate specimen, from Wisnom et al. \cite{Wisnom2008} }
	\label{fig:wisnom-exp}
\end{figure}

As an example, Fig. \ref{fig:wisnom-exp} shows the resulting failure pattern in an unnotched $[45_4/90_4/\text{--}45_4/0_4]_s$ laminate specimen in an experimental study by Wisnom et al. \cite{Wisnom2008} As illustrated, significant fiber or fiber bundle fractures accompanied by longitudinal matrix cracks between them occurred in $0^{\circ}$ plies, resulting in the failed fiber bundles/yarns being pulled out, while the failure in $\pm 45^{\circ}, 90^{\circ}$ plies was dominated by dispersed matrix cracks between fibers and yarns. Moreover, the intra-ply fiber and matrix cracks were also strongly coupled to inter-ply delaminations throughout the laminate, which resulted in localized failure 'blocks' and separated plies. Interactions between these discrete types of failures usually also affect the local stress distribution and damage progression in laminates. The damage pattern and ultimate strength are also affected by features such as: statistically distributed material strength, initial structural imperfections, specimen shape and geometry, inclusion of notches or holes.

One review of finite element modeling of progressive damage growth in composite laminates is that of Tay et al. \cite{Tay2008}, who also discussed the relevant issues in structural applications. Wisnom \cite{Wisnom2010} reviewed the application of cohesive interface elements for modeling polymer composites dominated by discrete matrix cracks, who also described various damage phenomena of multi-mode, interacting discrete cracks, and emphasized the necessity of properly modeling these features. Another review of recent work in laminate progressive damage analysis was given by Liu and Zheng \cite{Liu2010}.  

In the following subsections, we continue this review and discussion of existing methodologies for modeling damage evolution in laminates. To better show the importance of accommodating multiple discrete damage modes, and for organizational purposes as well, the following subsections are classified according to the types of discrete failure features included in the numerical models.

\subsection{Models of laminate failure based entirely on homogenized standard elements} \label{no-discrete}

Many failure studies used standard finite elements in conjunction with Ply Discount Method (PDM) or Continuum Damage Mechanics (CDM) approach, and thus, did not focus on the discrete nature of damage modes described above. 
Laminates in these models are constructed using standard elements, where plies are either smeared (in the thickness direction) into one single mesh layer, or a mesh of multiple layers bonded together by sharing common nodes. 
Then, CDM-based failure criteria, such as maximum strain, maximum stress, or Tsai-Wu, are often used to evaluate the consequences of damage \cite{Kaddour2004}. In cases where a single mesh layer represents multiple plies, effective stiffness properties are computed by Classical Laminate Theory (CLT), then by inverse computation based on the compliance matrix of each ply, stresses in each layer can be obtained through extrapolation for the purpose of applying failure criteria on each ply \cite{Gay2003}. 

Even though the failure criteria can be assessed in local coordinates to account for fiber and matrix damage states, failures are still described in terms of homogenized element blocks, and therefore, the capability to explicitly model discrete and discontinuous damage features is absent. When a standard element has been degraded or considered failed, its entire volumetric space in the model is involved. Thus, failing a standard element over its entire volumetric space will result in unrealistically large energy dissipation and a large spatial void, which can lead to the common numerical problem of premature 'failure localization', i.e., the damage path forms as a single row of interconnected failed standard elements, when in reality the failure process is more distributed.

In physical space, a failed standard element also contradicts the idea that in-situ, discrete cracks in laminates are essentially surface discontinuities having extremely small process zones. As shown in micrographs of damaged laminates \cite{Takeda1994, Swolfs2015}, the process zones of cracks (thickness or gap normal to the fracture surface) are generally smaller than $0.5 \text{mm}$ for delaminations between plies and even shorter for transverse matrix cracks, commonly around $5$ to $100$ microns ($\mu \text{m}$).

In CDM-based framework, several methods have been presented in literature to overcome such limitations. For example, Bazant and Oh \cite{Bazant1983} proposed a Crack Band (CB) theory initially adopted for fracture of concrete. In the model, although cracks are still treated as 'smeared crack' being represented by continuum elements, formulation of material degradation laws are modified by introducing a characteristic length, results are therefore objective with regard to element size. 
More specifically for composite laminates, Pineda and Waas \cite{Pineda2013} presented a thermodynamically-based work potential theory, where pre-peak damage laws (typically non-linear softening curves) of lamina materials are formulated by extending the Schapery Theory (ST) \cite{Schapery1990}, while post-peak failure laws (degrading curves) are based on CB theory. Formulations for element characteristic length and constitutive equations are derived based on work potential theory to ensure preservation of the total energy dissipated from mesh to mesh.

In brief, homogenized elements with CDM-based failure modes are computationally efficient, satisfactory numerical results could be obtained when combined with effective modifications. Due to the nature of smeared continuum elements, fracture surfaces and discontinuities with small process zones are absent. Therefore, CDM-based continuum elements may show limitations in certain conditions, for example, when high-fidelity representation of fracture is essential, when the discreteness of failures is critical to structural damage evolution, when multi-mode interactions between various types of cracks are key mechanisms affecting the progressive damage of laminates. 

\subsection{Models assuming only one type of discrete failure} \label{single-discrete}

Many numerical methods have been developed in the past decades to account for discontinuous features of fracture surfaces and crack evolution in general continuum and solids. These techniques have been intensively applied for modeling damage and failure evolution of composite laminates. 
Commonly used methodologies are:

(i) Debonding and modification on elements nodal connectivity, such as Virtual Crack Closure Technique (VCCT) \cite{Rybicki1977, Krueger2004}.
(ii) Element enrichment techniques based on the Extended Finite Element Method (XFEM) and its extensions \cite{Moes1999, Moes2002}. It should be noted that fracture surfaces are still absent in XFEM-based models, since the discontinuity are described and approximated by enrichment of shape or nodal functions of elements by numerical means, i.e., no real discontinuous in finite element model. Nevertheless, XFEM is computational efficient and has shown its effectiveness in many situations. 
(iii) Insertion of elements having zero or small thickness at potential crack locations. The inserted elements can either be generated as distinct ones (in a line or thin layer) tied to the adjacent parts, or embedded inside the nodes of priorly-generated structured meshes. The feature of small or zero thickness of the inserted elements provides good representation of small process zones and discontinuous fracture surfaces. Due to its small or zero volume, they may no longer be treated as continuum elements, thus the constitutive laws of inserted elements are usually defined in terms of Cohesive Zone Modeling (CZM) technique (such as interfacial traction separation law), therefore they are commonly called as Cohesive Elements (CEs) or interface elements. 
(iv) Element partitioning, in which an element is separated or 'cut' into two or more sub-domains to form discontinuity and fracture surfaces in numerical model.
The latter two methods have been widely used and modified in many recent studies, as will be discussed in the following sections.

\subsubsection{Models allowing only discrete matrix cracks} \label{only-matrix-crack}

Arrays of intra-ply matrix cracks oriented along the fiber direction (often called transverse matrix cracks, since they are most prominent in $ 90^{\circ}$ plies) usually initiate during the early stages of mechanical loading, whether a monotone tension test, or a fatigue test with increasing load cycles. Further growth of such matrix cracks results in stress redistribution and concentration between plies, which in turn triggers additional modes of failure such as inter-ply delaminations and sequential fiber and yarn failures that finally branch out and split the laminate. Thus, intra-ply matrix cracks are often considered to be a precursor, or even a trigger, for progressive failure in laminates. 

Many studies\cite{Lundmark2006, Singh2009, Vinogradov2010, Huang2014, Carraro2015, Asadi2015, Maimi2011, Barbero2011} have focused on: (i) analytical treatments of perturbations in stress fields and stiffness degradation caused by transverse matrix cracks, (ii) constraint effects and interactions between systems of matrix cracks in multi-directional laminates, and (iii) effective stiffness reduction and its correlation to increased crack density from increasing loads or load cycles in fatigue. 

In these studies, such transverse matrix cracks are generally assumed to be pre-existing short cracks that have subsequently grown across the entire widths of their respective off-axis plies. The intra-ply crack density is typically viewed in terms of fully developed, parallel matrix cracks with specific spacing. In these analyses, the crack densities are treated as a homogenized process, rather than through an explicit description of the evolution of the crack system, originating from discrete and dispersed short local matrix cracks. Several recent studies have applied a multi-scale modeling strategy to predict effective crack density evolution in plies, but mostly in symmetric laminates with special and convenient ply angles \cite{Jagannathan2016, Carraro2017, Montesano2018}.

Interestingly, inter-ply delaminations are basically ignored in these numerical studies, as their focus is narrowed to already-complicated matrix crack systems. In reality, propagation of matrix cracks leads to the onset of delaminations at matrix crack fronts, affecting and changing local stiffness and stress re-distribution profile. Then as damage accumulates, delaminations extend and connect matrix cracks from abutting plies, thus forming local networks of discontinuities including yarn or tow splitting. The consequence can be constraining or shielding effects, which usually emerge in later stages of laminate damage growth \cite{Carraro2017}. These various kinds of complicated interacting mechanisms remain active until complete failure of the laminate. 

Experimental observations also have confirmed the important role of inter-ply damage modes (delaminations) involved in crack initiation and growth, as well as the apparent strong interaction between inter-ply damages and matrix cracks. Contributions made by delaminations can cause failure to occur much sooner than expected, especially if the specimen has initial notches or under a fatigue loading condition. For example, in studies of notched laminates \cite{Hallett2009, Li2009}, crucial localized delaminations emerged early around open hole edges as a result of the associated stress concentration, soon after which delaminations propagated quickly and became the main mechanism causing ultimate laminate failure. In another study of unnotched laminates by Tong et al.\cite{Tong1997}, local delaminations were found to be distinct features of fatigue damages, wherein delaminations occurred at nearly all ply interfaces when certain specimens were tested in loading cycles. However, for the identical laminate specimens, such delamination phenomena were largely absent in quasi-static loading tests.

Although numerical models with only intra-ply discrete matrix cracks can yield reasonable predictions of stress fields and effective laminate stiffness degradation in certain cases, the inability to model interacting intra and inter-ply discrete damages limits their use to conditions where such laminates have low matrix crack densities and no initial stress concentrations from notches or open holes.

\subsubsection{Models allowing only discrete inter-ply delaminations} \label{only-delamination}

Many numerical models have been developed that focused only on delamination initiation and propagation. For modeling of pure delaminations, CE and VCCT are the two mostly used methods. An overview of early work on such models was given by Tay \cite{Tay2003}, a brief assessment and discussion of recent strategies for modeling delaminations can be found in \cite{Heidari-Rarani2019}.
If CDM-based intra-ply damages are also included in the model, accurate and direct coupling between intra-ply discrete damages and inter-ply delaminations will be absent, since the homogenization nature of ply continuum elements leads to the loss of key information required to realistically describe damage interactions \cite{Talreja2006, Ling2009}. As also pointed out in several studies, CDM-based intra-ply damage modes were found to be inaccurate in predicting crack paths \cite{VanDerMeer2009}, and unable to predict the correct fiber stress relaxation behavior resulting from longitudinal spitting \cite{Iarve2005}. In general, models with only discrete inter-ply delaminations may be essentially limited to cases where the transverse matrix cracks and fiber failures are negligible, or when the interaction between them is not a concern, such as in a two-point Mode-I inter-laminar fracture test.

\subsection{Models with multiple interacting discrete cracks} \label{multi-discrete}

In order to treat multiple cracks and their bifurcations or intersections, improvements on the traditional finite element methods are necessary to address the limitations of treating such cracks as though they are non-interacting.

Many studies use XFEM or CDM-based methods for intra-ply damages (usually as a continuation of their earlier work on single ply models), then attempt to add modified numerical treatment into multi-ply models to couple intra-ply damages with inter-ply delaminations. Despite the absence of real physical intra-ply fracture surfaces and discontinuities, it is a step forward from the models mentioned in Section \ref{no-discrete} and \ref{single-discrete}, since at least some sort of interacting mechanisms are applied. 

For example, Joseph et al. \cite{Joseph2018} applied Crack Band model in continuum elements for intra-ply failure modes, and small-thickness interface elements for inter-ply delaminations, respectively. In order to simulate the coupling between different types of damage modes, an 'intra-ply failure triggered delamination initiation' function is added to the model. Simply speaking, inter-ply interface elements along the intra-ply matrix crack paths are degraded, as if their failures are triggered by the neighboring matrix cracks from abutting plies. A global crack tracking method is implied, more specifically, a user-defined routine to artificially block damage initiation within a defined distance around failed elements, thus largely avoiding mesh-sensitive 'failure localization' commonly seen in CDM-based elements as described above (Section \ref{no-discrete}).  
Sun, Hu and co-authors \cite{Sun2013, Hu2016} applied XFEM-based structured elements to model intra-ply transverse matrix cracks, and small-thickness cohesive elements between plies to model delaminations. During the numerical analysis, ply elements are classified according to their distances from a nearby crack tip (if any), then modified enrichment formulations are implied accordingly. Additionally, inter-ply cohesive elements around a matrix crack tip are also enriched to account for the interaction between two intersecting cracks. Simply speaking, both ply continuum and inter-ply cohesive elements at crack intersecting locations are enriched by special formulations to approximate the coupling between different failure modes. Higuchi et al. \cite{Higuchi2017} used a similar XFEM/CZM coupled approach with added modification on the interfacial constitutive laws, wherein all crack sites are pre-inserted to limit the number of internal variables. 

In the above studies, special modifications are made into the traditional XFEM and/or CDM-based elements to account for coupling between different failure modes. However, in these finite element methods, the ply elements are still not partitioned into sub-domains, i.e., in-ply fracture surfaces and discontinuities are absent, and displacement fields around cracks are either represented by failed continuum elements (smeared crack bands), or approximated by enrichment functions. Consequently, features such as interaction between fracture surfaces, crack openings and closure, are not captured. This drawback also means the inability to couple intersecting discontinuities in a explicit manner.
Admittedly, XFEM and CDM-based methods are computationally efficient, and as demonstrated by most of the aforementioned studies, can generate excellent and robust results in various situations. But the point here is that real and high-fidelity representation of fracture surfaces and discontinuities could be essential and critical for many applications, especially under situations where multi-mode interaction between and coupling of crack networks are fundamental phenomenas governing the progressive failure of composite laminates.

Actually, many studies adopted element partitioning methods or interface elements to have the capacity of modeling real fracture surfaces.
For element partitioning in FEM, two similar methods have been developed, which are based on an early work by Hansbo and Hansbo \cite{Hansbo2004}:
In the first, Belytschko and co-authors proposed a Phantom Node Method (PNM) \cite{Song2006}, wherein 'phantom' nodes are introduced to overlap an element's true nodes to add additional degrees of freedom (DOF), thereby allowing the crack to be modeled by the superposition of two sub-domains divided by a crack.
Then second is by Ling et al. who proposed an Augmented Finite Element Method (A-FEM) \cite{Ling2009}, wherein an element can be directly partitioned into two sub-elements to account for a single discontinuity.
Prabhakar and Waas \cite{Prabhakar2013} also presented a Continuum-decohesive Finite Element (CDFE) for intra-ply cracks in composite laminates, wherein an initially-continuum element can be cut across the integration point or centroid to form two sub-domains. In such kinds of methods, once real fracture surfaces are formed, new stiffness matrix and force vectors of sub-elements are computed accordingly, user-defined interacting behavior between fracture surfaces, such as cohesive traction-separation laws, may be added.   

\begin{figure}[H]
	\centering
	\begin{subfigure}[h]{0.45\textwidth}
		\centering
		\includegraphics[width=\textwidth]{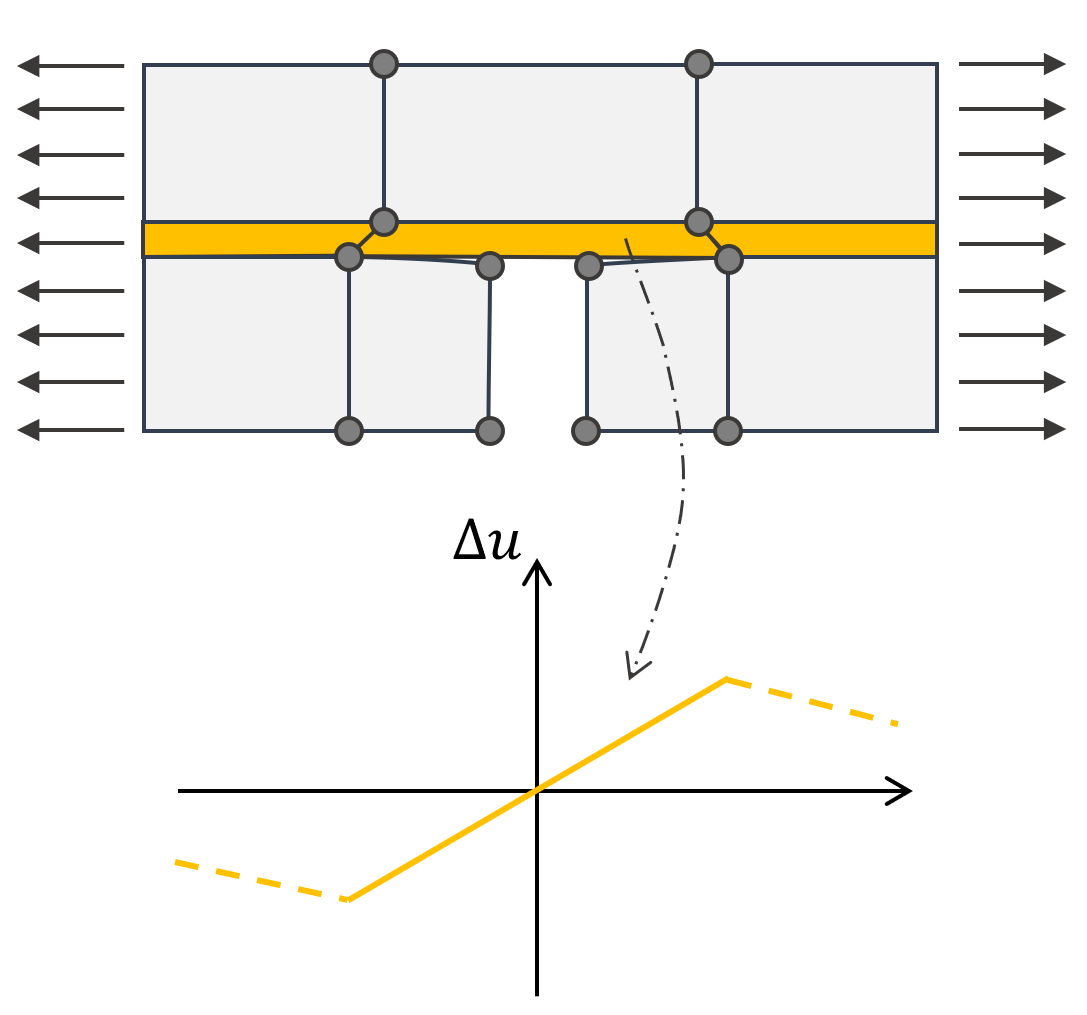}
		\caption{Non-matching mesh at a possible crack bifurcation.}
		\label{fig:non-matching-ce}
	\end{subfigure}
	\hfill
	\begin{subfigure}[h]{0.45\textwidth}
		\centering
		\includegraphics[width=\textwidth]{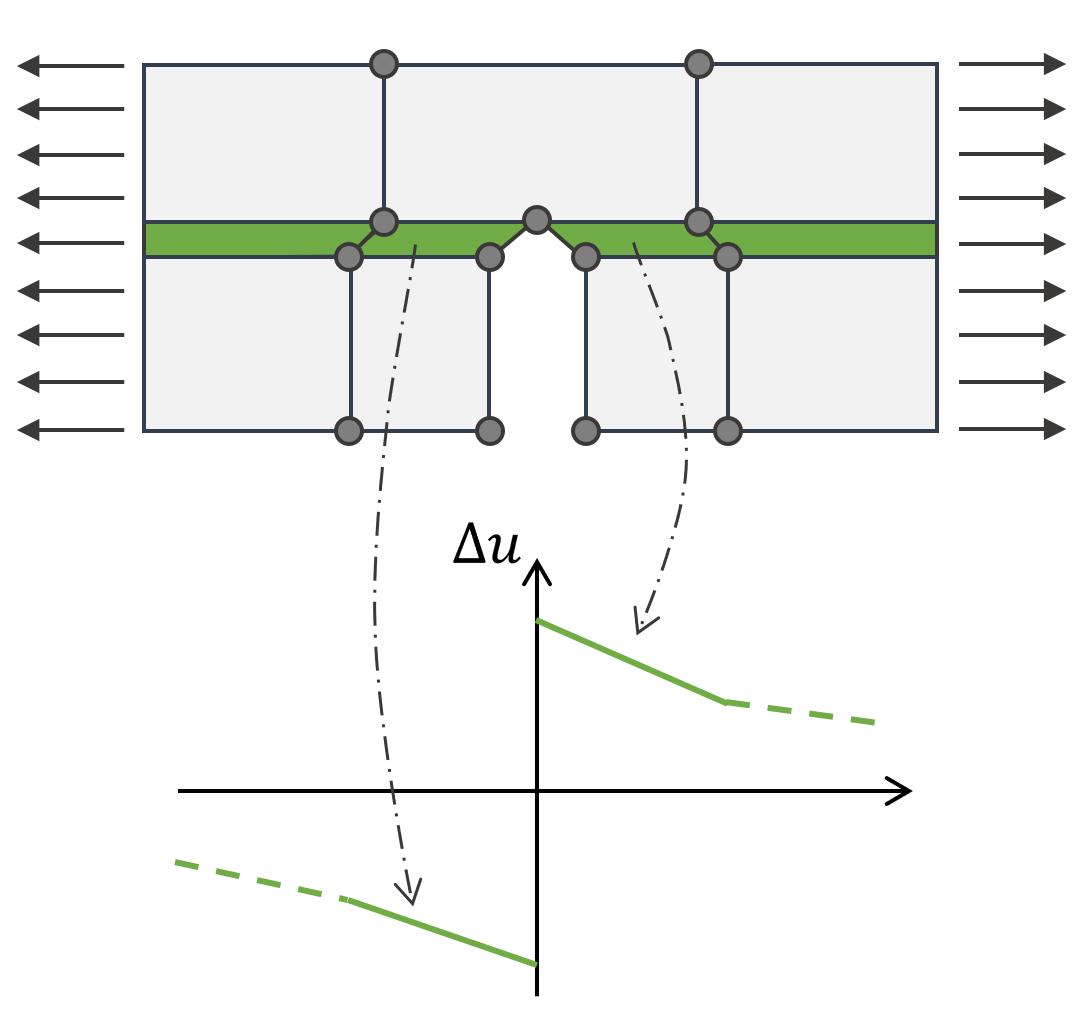}
		\caption{Matching mesh at a possible crack bifurcation with bonded elements.}
		\label{fig:matching-ce-bonded}
	\end{subfigure}
	\caption{Comparison of a non-matching and matching mesh at a possible crack bifurcation to form a T-shaped crack, after Fang et al. \cite{Fang2011}
		 (a) Non-matching mesh fails to properly capture the displacement jump along the crack tip, thus resulting in erroneous strains and load transfer at a potential crack bifurcation; 
		(b) A matching mesh using a separable element with bonded elements sharing nodes (FNM \cite{Chen2014} or A-FEM \cite{Fang2011}), thus capturing the displacement jump at a crack tip and high local shear stresses driving the potential formation of a T-shaped crack.}
	\label{fig:comparison-of-matching-and-non-matching-ce}
\end{figure}

However, it remains a tricky problem to couple two or more discrete cracks properly (when multiple fracture surfaces are intersecting and forming crack networks). 
Fang et al. \cite{Fang2011} showed that a standard cohesive element next to partitioned sub-elements leads to non-matching mesh at crack tip, and thus, is unable to capture the displacement jump across the crack surface. As shown in Fig. \ref{fig:non-matching-ce}, this shortcoming, due to non-matching mesh, results in erroneous shear strains and load transfer at crack tip, therefore fails to model potential crack bifurcations (such as a T-shaped crack). To address this issue, they developed an Augmented Cohesive Zone (ACZ) element \cite{Fang2011} using A-FEM, wherein the coordinate information of a crack tip is stored and passed onto its neighboring element. This way, the partitioning line is determined by the tip of an approaching crack from an adjacent element, thus creating a matched mesh allowing for the possibility of a crack bifurcation, as shown in Fig. \ref{fig:matching-ce-bonded}. 
Chen et al. \cite{Chen2014} later proposed a Floating Node Method (FNM) similar to PNM, but instead uses floating points on an element's edges to track discontinuities, thereby partitioning elements according to the current and predicted crack pattern. The extra flexibility introduced by floating nodes leads to several advantages over PNM, for example, the elimination of errors in the mapping of straight cracks in PNM/XFEM. 
Detailed description of the aforementioned methods can be found in various publications, in which several two-dimensional (2D) numerical examples are also provided. In 2D, only matrix cracks and delaminations need to be considered, and in theses cases the two types of failures are orthogonal to each other.

In 3D applications, situation is even more complicated, since the varying angles between plies and all three types of damage modes must be properly treated. Just like in 2D, it is critical to have matching meshes at crack bifurcations to capture displacement jumps and crack interactions accurately. Several recent studies have presented numerical models where at least two types of failures with explicit and discontinuous fracture surfaces are supported for composite laminates in 3D space:

Higuchi et al. \cite{Higuchi2017a} used a model where several rows of cohesive elements are inserted in certain plies to account for matrix cracks, and independent thin layers of cohesive elements are put between plies to model delaminations. This combination gives real fracture surfaces of matrix cracks and delaminations, which is found to yield numerical predictions with higher accuracy than a comparable model where all matrix cracks are CDM-based. However, since delaminations are modeled by independent layers of cohesive elements constrained between plies with surface-based tie constraints, this method will inevitably lead to non-conformal mesh (with dissimilar mesh patterns) between inter- and intra-ply elements regardless of ply angles. Consequently, displacement jump and explicit coupling between cracks can not be captured due to the non-matching meshes at potential crack bifurcations, as explained above (Fig. \ref{fig:comparison-of-matching-and-non-matching-ce}). Besides, fiber failures in the model are based on smeared crack model, and thus are still CDM-based without explicit fracture surfaces.    

Lu et al. \cite{Lu2018} proposed a Separable Cohesive Element (SCE) based on floating node method (FNM) to model delaminations. Such an element is embedded between ply elements, and can be partitioned according to the existence and pattern of cracks from adjacent plies (one from above and one from below), where intra-ply cracks are formed by partitioning ply elements using FNM.
This way, the partitioned sub-domains of an inter-ply SCE will match the pattern of one or tow cracks formed in adjacent ply elements (partitioned by floating nodes), thus forming matching meshes between inter-ply and intra-ply cracks.
However, a SCE is only capable of matching at most two cracks \cite{Lu2018}, meaning that only one type of intra-ply failure is allowed for a specific ply. Conceivably, even if both types of intra-ply failures were to be considered (upon further model development), an inter-ply SCE must be partitioned at least twice (in two numerical increments) to form matching meshes between itself (delamination) and both types of intra-ply failures. More specifically, this has to involve two matrix cracks plus two fiber fractures from adjacent plies (four in total), in the form of two intersecting 'T-shaped' or 'X-shaped' crack networks. 
Therefore, the simultaneous treatment and coupling of all three types of failures may not be explicitly captured, since the desired element partitioning requires more than one numerical increment. Such limitations could lead to more gradual and milder damage interactions between cracks, which could be an issue in cases where extremely rapid damage growth is physically likely, and where virtually instantaneous local interaction between cracks is critical to the damage evolution and stress redistribution. Such increment-by-increment element partitioning mechanism may also affect the pattern and propagation of crack networks. 
In a later section, a numerical example of laminate failure will be discussed, in which fiber breaks, that took place when a specimen failed abruptly in the final stage of a static loading test, are found to have strong effect on the failure modes of the specimen.
Besides, in FNM and its extensions, floating nodes must be inserted on all edges of a separable element \cite{Chen2014, Lu2018}. At the outset, this requires historical internal variables and floating nodes throughout the simulation to pass information between elements. As a reference, before a SCE is partitioned in a numerical increment, sixteen possible geometrical scenarios must be evaluated so that the SCE can be properly partitioned to match the pattern of up to two cracks from abutting plies. Once an element has been partitioned into sub-domains, new node-connectivity constraints and more floating nodes have to be added, new stiffness and assembling functions require computation as well.

In another work, Joosten et al. \cite{Joosten2018} applied preprocessing tools to embed zero-thickness elements among structured standard elements to anticipate all three types of failure modes. The embedded elements are not physically created, but instead, are defined and numerically formed by the nodes of priorly-generated ply elements. Simply speaking, creating 'springs' between nodes of bonded continuum ply elements.
Then these 'interface springs' are categorized into three groups (defined with corresponding constitutive laws and failure criteria) based on the damage mode involved: fiber fractures, matrix cracks and delaminations. The formation of such interface elements in null space will ensure matching nodes at potential crack bifurcations, thus direct and simultaneous coupling between multiple cracks can be captured. However, since fiber fracture surfaces are orthogonal to transverse matrix cracks, while each ply must have a fiber-aligned mesh determined by its layup angle, this method will inevitably result in small or skewed elements with poor mesh qualities, which consequently results in the inability to generate elements in laminates with arbitrary layup. A possible way to avoid undesirably skewed meshes is to use surface-based contact relations between structured elements of different plies, but as pointed out by the authors themselves \cite{Joosten2018}, displacement jumps may not be represented properly between intra and inter-ply cracks because of the involved master/slave pairing mechanism in the surface-based contact formulations. 
Actually, the surface-based contact behavior in FEM is based on the interpolation of nodes on the involved surfaces, one of the two is usually defined as a 'slave surface', such that its DOFs follow the other 'master surface'. As result, for any surface-based contact pair, matching meshes at potential crack bifurcations are generally absent, displacement jumps and explicit load transfer between cracks are not properly treated because those connectivity informations and stress computations are based on interpolation of surrounding nodes. 
Thus, this method is suitable primarily for cases where the aspect ratio of fiber-aligned standard elements are well-controlled, and where the discrete cracks can be represented by the edges or diagonals of the standard elements, namely, if the angle between plies is either $0^{\circ}$ or a multiple of $45^{\circ}$. Consequently, the laminate ply angles presented in their study were limited to $0^{\circ}, \pm 45^{\circ}, 90^{\circ}$. 

Another earlier method by Bouvet et al. \cite{Bouvet2009,Bouvet2013} also used zero-thickness interface elements embedded in priorly-generated ply mesh to account for matrix cracks and delaminations. In their model, each ply is meshed into 'little strips' with one continuum element representing the width and thickness of one strip, then zero-thickness interface elements (intra-ply 'springs' of null length) are defined between neighboring strips (within the same layer) to capture potential matrix cracks; similarly, inter-ply 'springs' are defined between locally overlapping nodes from adjacent plies to account for effective inter-ply delaminations at a local representative region. However, fiber failures are still modeled by CDM-based methods, wherein the stiffness of a continuum element is rapidly degraded to simulate longitudinal fiber/strip damage. Moreover, since the creation of 'matrix crack springs' requires priorly-generated meshes consist of 'strips' that must be fiber-aligned with respect to each ply angle, meanwhile the definition of 'delamination springs' requires overlapping nodes arranged in a grid pattern from adjacent layers, this method consequently lead to similar restrictions on the choice of layup angles as mentioned above. And it turns out that the ply angles in their studies were also limited to $0^{\circ}, \pm 45^{\circ}, 90^{\circ}$.   

\subsection{Summary}\label{intro-summary}

To summarize, it has been shown in this section that the commonly used techniques have limitations in modeling discrete damages and their couplings in composite laminates with arbitrary layup.
A comprehensive review of recent studies indicates that: it is not only important to represent discrete failure modes with real fracture surfaces, but also essential to have all three types of cracks; at the same time, it is critical to have matching meshes at potential crack bifurcations to ensure explicit load transfer and couplings between cracks, which is necessary not only at locations where crack bifurcations have already occurred, but also where they can be anticipated as damages propagate.    

In brief, there are currently two types of methods approaching these tasks: The first method is to partition elements 'on-the-go', such as by using Floating Node Method (FNM) for intra-ply cracks and Separable Cohesive Element (SCE) based on FNM for inter-ply delaminations. The second method is to define interface 'spring' elements between nodes of priorly-generated structured mesh, or, to insert small-thickness cohesive elements along the paths where cracks are anticipated as damages occur and propagate. However, an inevitable drawback of the second method is that un-conformal and/or highly skewed meshes will occur between plies that have certain layup angles, and thus, the method is not suitable for modeling laminates with non-standard stacking sequences. While in the first method, a delamination SCE is only capable to form matching nodes with one type of intra-ply damage mode, i.e., either matrix cracks or fiber fractures will be absent in a ply. Besides, simultaneous formation and coupling of discrete cracks may be improperly handled due to the nature of such method, which is likely a critical issue when laminate damage grows rapidly, though extremely small increment may be adopted.

It should be noted that these two types of methods are summarized based on representative recent studies attempting to account for interactions between discrete cracks, similar or comparable modeling techniques may also be seen in several other studies. There are multiple notable authors who have done extensive amount of work, for example, Yang and co-authors. \cite{Cox2006, Ling2009,Fang2011}, Iarve et al. \cite{Iarve2005,Iarve2017,Iarve2018} and Van Der Meer et al.\cite{VanDerMeer2009,VanderMeer2010,VanDerMeer2012,VanDerMeer2013}, are among the early few groups to approach discrete damage modeling in composite laminates. Indeed, Yang is one of the contributors in Augmented Finite Element Method (A-FEM) mentioned above \cite{Ling2009,Fang2011}, which is a essential fundamental of the first type of method (element partitioning), Yang also contributed in the second type of method (element embedding) referred here as work by Joosten et al. \cite{Joosten2018}. 
Besides the discrete modeling approach with real fracture surfaces, many CDM-based and XFEM-based methods have been intensively modified and developed in-house, thus showing excellent capability and computational efficiency in various applications. For example, the work by Thorsson and Waas and co-authors in the studies of laminates involving impact damages \cite{Thorsson2018a, Thorsson2018b} and progressive failures \cite{Pineda2013, Nguyen2016}, Waas also contributed to the element partitioning methods as mentioned above (CDFE) \cite{Prabhakar2013}. 

Nevertheless, the motivation raised in this paper is to create a new method with the capability of modeling high-fidelity fracture surfaces of all three types of failure modes in composite laminates having arbitrary layup angles, meanwhile ensuring matching nodes among crack bifurcations for accurate coupling and load transfer between cracks, and to do so with a code that is computationally efficient. For this purpose, a novel modeling method, somehow differs from the current aforementioned existing methodologies and approaches, has been developed, as will be demonstrated in subsequent sections. 

The remainder of this paper is organized as follows: Section \ref{sec2:AGDM} provides a detailed description of the proposed modeling method. Section \ref{sec3:numerical_example} gives some numerical examples that successfully capture the progressive failure behavior and strengths seen in experimental samples, thus illuminating the various modeling challenges and failure phenomenas described above. Section \ref{sec4:conclusion} gives some conclusions.

\section{Auto-Generated Geometry-Based Discrete Model (AGDM)}\label{sec2:AGDM}

The methodology of AGDM is to represent the volume of a laminated composite object of interest in terms of discrete 'parts', which can be either 'conventional parts' representing resin-impregnated yarns or yarn segments (viewed primarily as 3D continuum structures), or, small-thickness 'interface parts' (almost 2D planar structures). These 'interface parts' can represent the interface (commonly voids or resin rich zones) between two adjacent yarns or between two plies, but they can also represent the location of a weak spot along a yarn (in reality a local collection of weaker fibers) where such a yarn can break in half. Essentially these 'interface parts' anticipate the occurrences and various locations of one of the three crack types (failure modes) and their possible interactions. All such parts will typically be partitioned or subdivided into smaller segments, not only to cover all potential failure sites (within a reasonable length-scale of local load-transfer between the various parts), but also taking care to ensure matching meshes at all potential crack bifurcations. Finally, all parts are translated to their appropriate locations in the space occupied by the laminate object, and then assembled together by explicitly defined surface-based tie constraints among certain surfaces to form the whole laminate in modeling space. These operations are realized by in-house developed preprocessing codes (Python scripting commands) within the commercial FEM software Abaqus \cite{DassaultSystemes2016}.  

Given the roles of the various 'parts' in the laminate failure process, and to have a consistent and easily-remembered naming convention, the continuum 'conventional parts' are called 'yarns', or when subdivided, 'yarn segments', which may be understood as bundles of thousands of fibers (e.g., carbon) bonded together by a matrix of some volume fraction (e.g., epoxy at 40 percent), which is virtually the same volume fraction as in the laminated composite as a whole. 
On the other hand, the very thin 'interface parts' are inserted to allow for the possible formation of discrete failures of different types, or namely, cracks that may occur. These include: (i) impregnated yarn failures (usually modeled as discrete cracks across yarns), (ii) intra-ply matrix cracks (developing and growing along the interface between two yarns), (iii) delaminations between plies (cracks growing at ply-to-ply interface). 

In anticipation of their possible triggering of such cracks, interface parts are called 'cracklets', to instill the notion that initially they are not cracks, and mechanically they are virtually taking zero volume, and are not affecting the mechanical properties of composites. But these embedded cracklets may trigger upon certain criteria, and become segments of cracks and material discontinuities. Specifically, there are (i) 'yarn cracklets' for yarn tensile fractures, (ii) 'matrix cracklets' anticipating interface damages between adjacent yarns in a ply and subsequent intra-ply matrix cracks propagated through linking of failed matrix cracklets, (iii) 'delamination cracklets' anticipating the triggering of inter-ply delaminations and their growth from such cracklet linking. Load transfer from failure of a cracklet may trigger failure of a neighboring cracklet of the same or a different type, thus generating laminate damage in the form of a growing networks of cracks. Matching mesh are ensure at any potential crack intersections in 3D space. 

In the following subsections, we provide details of constructing the various 'parts' and 'cracklets' of the model assembly based on their order of generation. 

$\bullet$ Initial ply discretization in Section \ref{sec2.1:ply_initiate}; 

$\bullet$ Yarn segmentation, and insertion of yarn cracklets in Section \ref{sec2.2:yarn_segmentation}; 

$\bullet$ Partition of yarn interfaces into matrix cracklets in Section \ref{sec2.3:matrix_segmentation}

$\bullet$ Partition of ply interfaces into delamination cracklets in Section \ref{sec2.4:delamination_segmentation}; 

$\bullet$ Assembly of all parts and overall constraints in Section \ref{sec2.5:assmebly}

\subsection{Initial ply discretization of yarn and yarn interface}\label{sec2.1:ply_initiate}

Key input parameters describing the geometry and layup angles of the laminate plies, taken to be a rectangular prism, are listed in Table \ref{tab:geo_input}.

\begin{table}[H]
	\centering
	\caption{Geometrical information required to generate the discrete laminate model} \label{tab:geo_input}
	\begin{tabular}{l l}
		\hline
		Input & Definition \\[1ex]
		\hline \hline
		$W $ & Width of the rectangular laminate specimen \\
		$L $ & Length of the rectangular laminate specimen \\
		$N $ & Total number of plies in the laminate specimen \\
		$h_i \,$ & Thickness of the $i$th ply  ($i=1,2,...,N$) \\
		$\theta_i \,$ & Layup angle of the $i$th ply ($\theta_i \in [-90^{\circ}, 90^{\circ} ]$)\\
		$d_i \,$ & Spacing of potential matrix cracks in the $i$th ply (gives matrix crack density) \\
		$l_i \,$ & Spacing of potential yarn fractures along the fiber direction in the $i$th ply \\
		$t_f \,$ & Nominal thickness of a yarn cracklet (potential yarn fracture) \\
		$t_m \,$ & Nominal thickness of a yarn interface (matrix cracklet) \\
		$t_d \,$ & Nominal thickness of a ply interface (delamination cracklet) \\
		\hline
	\end{tabular}
\end{table}

Upon establishing the user-input parameter values, the code begins by discretizing each ply (rectangular in the current implementation) consists of discrete yarns and yarn-to-yarn interfaces (to anticipate matrix cracks between yarns). The yarns in a ply will have the same local orientation as the layup angle of that specific ply, $\theta_i$.
The pattern of ply discretization is assumed to be symmetric about the rectangular plane of the laminate. Specifically, a 'center yarn' of each ply is set to pass through the ply such that the centroids of the center yarn in each ply are all aligned in the out-of-plane direction. Then, all remaining parts are assembled with respect to these centroids of the center yarns following the centrosymmetric rule. 

\begin{figure}[H]
	\centering
	\includegraphics[width=0.8\textwidth]{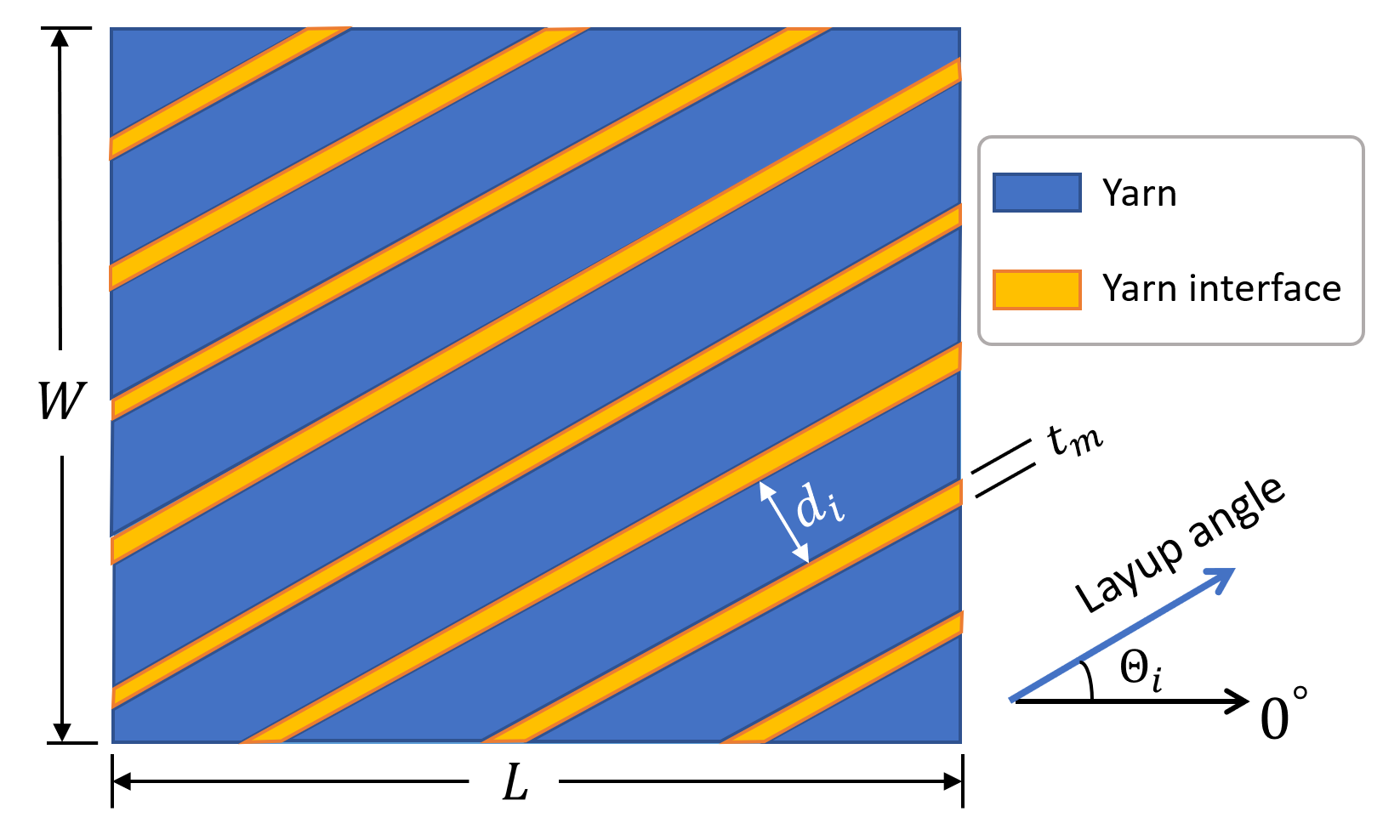}
	\caption{Initial ply-level discretization of the $i$th ply containing continuum yarns in blue and yarn interfaces in yellow (locations where matrix cracks can potentially initiate and develop).}
	\label{fig:in-ply-disrete-1}
\end{figure}

A diagram illustrating the initial ply-level discretization process is shown in Fig. \ref{fig:in-ply-disrete-1}, where the thicknesses of yarn interface, $t_m$ are enlarged for clearer view, but with the understanding that they are accounting for fracture surfaces of matrix cracks, and thus, are typically set to be much thinner than the continuum yarn parts (i.e., $t_m \ll d_i$, as mentioned in Section \ref{sec1:intro}).
For the $i$th ply, the required number and locations of discrete yarns and potential yarn interfaces are computed based on the layup angle, $\theta_i$, the matrix crack spacing of that ply, $d_i$, (approximately the fiber/epoxy yarn width) and the thickness of the yarn interfaces, $t_m$. More specifically, the coordinates of each part are computed and assigned by the code according to geometric principles, which involve finding the intersections of where the ply boundaries would 'cut' each ply part, i.e. where the outer prism boundary of the laminate cuts through each yarn and yarn interface, as shown in Fig. \ref{fig:in-ply-disrete-1} for a rectangular ply. 
Then each part for the $i$th ply will be created by extending its 2D planar representation (sketch) in the laminate stacking direction with magnitude of ply thickness $h_i$.
This way, all parts for the $i$th ply, when assembled together as in this example, will form a complete rectangular lamina ply with dimensions $W \times L \times h_i$.

\subsection{Insertion of yarn cracklets and corresponding yarn segmentation}\label{sec2.2:yarn_segmentation}

\begin{figure}[H]
	\centering
	\begin{subfigure}[htb]{0.8\textwidth}
	\includegraphics[width=\textwidth]{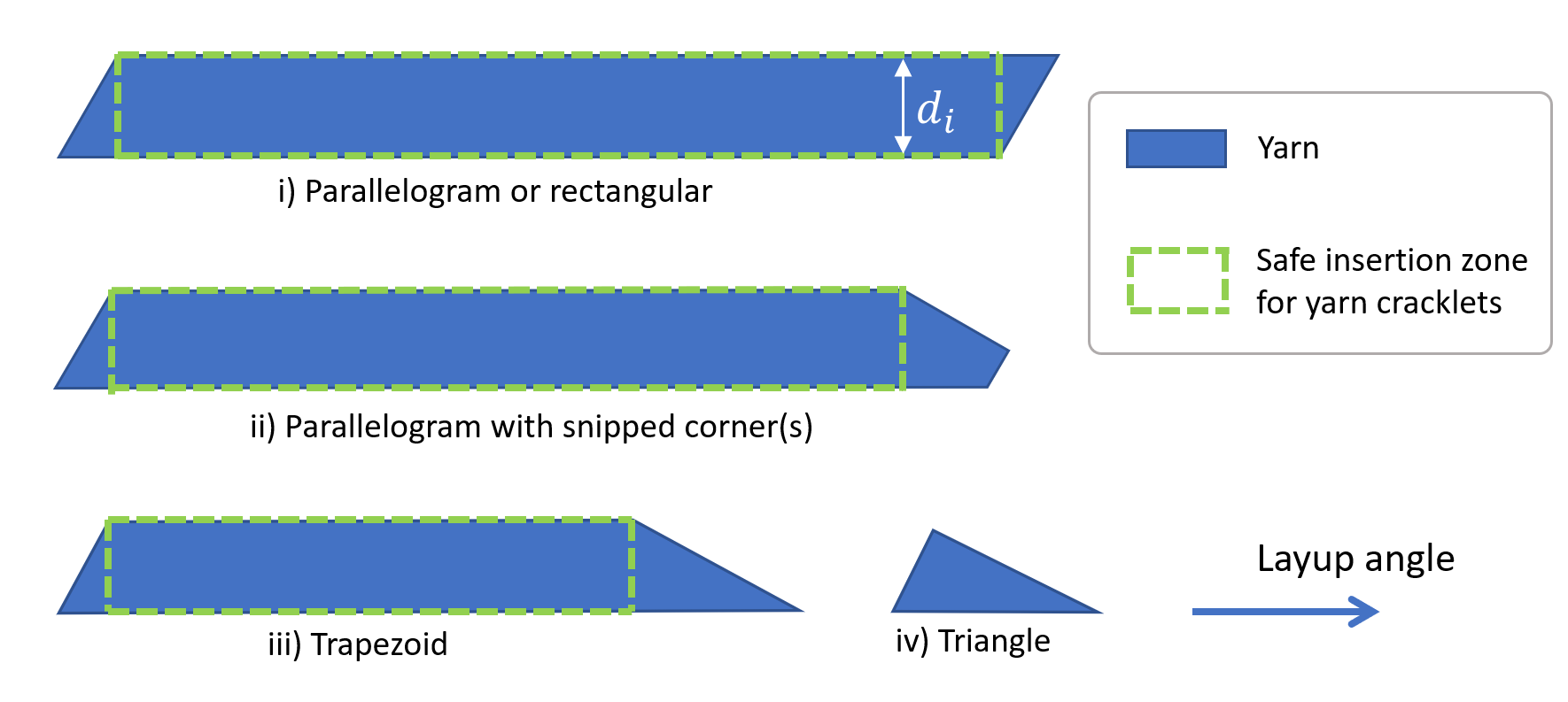}
	\caption{Four possible yarn shapes and their safe cracklet insertion zones.}
	\label{fig:yarn-shape-and-box}
	\end{subfigure}
	\begin{subfigure}[h]{0.8\textwidth}
	\centering
	\includegraphics[width=\textwidth]{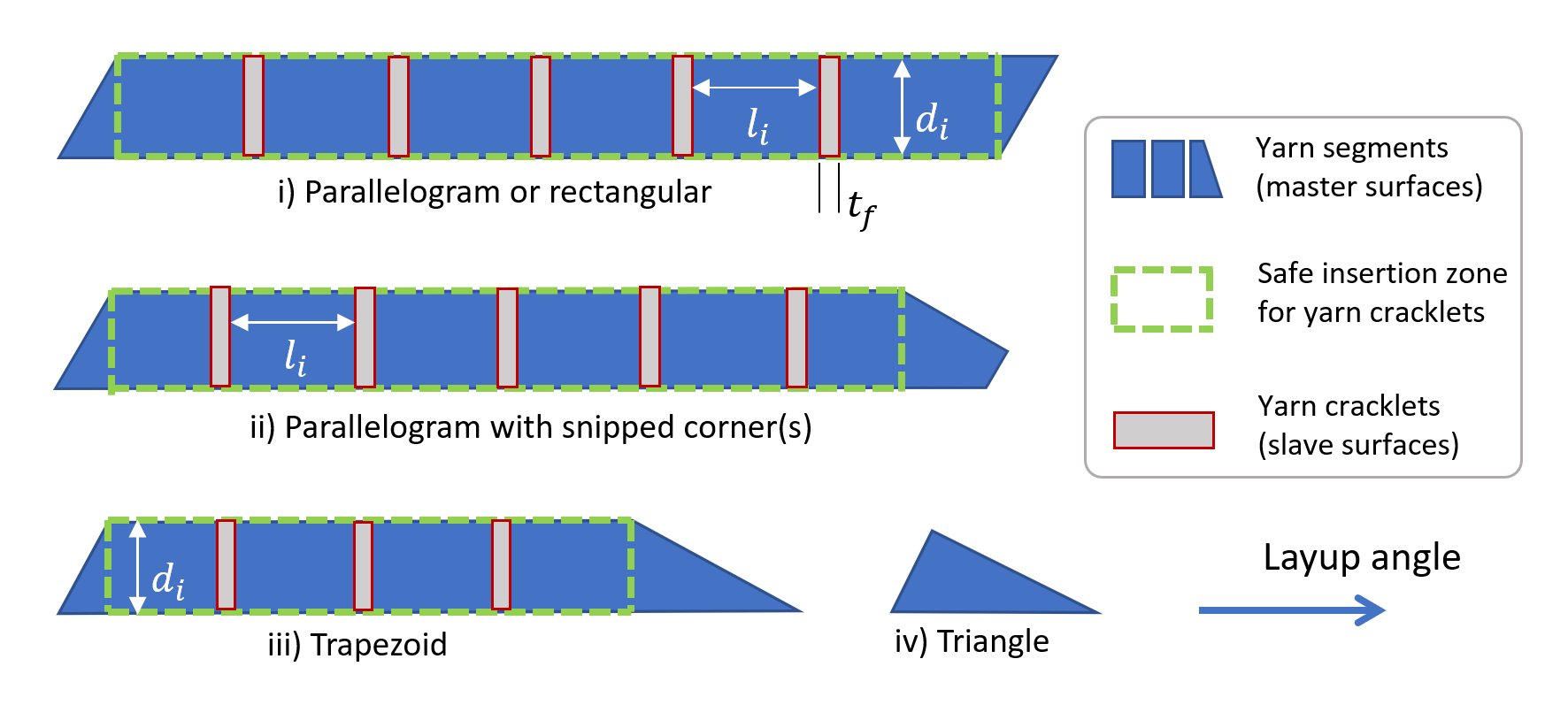}
	\caption{Yarns from Fig. \ref{fig:yarn-shape-and-box} partitioned into segments with inserted yarn cracklets}
	\label{fig:yarn-shape-and-insert}
	\end{subfigure}
\caption{Shapes of yarns and corresponding partitioning mechanism for yarn cracklet insertions (allowing potential fiber fractures).}
\end{figure}

Once the initial discretization of each ply has been performed, the next step is to insert yarn cracklets (locations of potential yarn tensile breaks) at predetermined locations along each yarn. Before doing so, it is important to point out that, the ply-level discretization typically results in several possible shapes of the yarn parts lying within the ply boundary. More specifically, a continuum yarn part can have one of four possible shapes (shown in Fig. \ref{fig:yarn-shape-and-box}):  
(i) a parallelogram (or a rectangle in $0^{\circ}$ or $\pm 90^{\circ}$ plies),
(ii) a parallelogram with snipped-off corner(s),
(iii) a trapezoid,
or (iv) a triangle. 

Note that a triangular yarn part is only expected to occur as a yarn region at a corner, thus by symmetry there will be either zero or two triangle yarns in a given ply. It is not reasonable, based on both experimental observations and mechanical means, to expect that fiber failures would occur at the very corner of plies with large layup angles. On the other side, a triangle yarn shape does not contain two parallel surfaces oriented along the fiber direction, and therefore, inserting small-thickness yarn cracklets to form potential 'fracture lines' that are perpendicular to a yarn edge (along the layup direction) will potentially lead to an extremely skewed mesh for certain layup angles. Given these facts, it is of little interest to insert yarn cracklets into triangular shaped yarn parts (at most two of them in ply corners).

From the above observations, yarns with the first three shapes are considered as candidates for insertion of yarn cracklets, and by geometrical features, these yarns all have two parallel edges representing two sides of the yarn. To establish a yarn cracklets insertion regime, the preprocessing code analyzes the coordinates of each continuum yarn part, tracks the lengths of its two parallel edges (along the ply layup direction), then defines the zone of common span of its two parallel edges (outlining a virtual rectangular 'box') as the 'safe insertion zone' for yarn cracklets (shown as a green, dashed box in Fig. \ref{fig:yarn-shape-and-box}). This way, the two surfaces of an inserted yarn cracklet, representing potential yarn fracture surfaces, will be perpendicular to the yarn's axis, thus ensuring stable and high-quality mesh generation. 

Based on this safe insertion zone concept and the initially-specified nominal (or effective) yarn cracklets spacing in the $i$th ply, $l_i$, the number and locations of required yarn cracklets for each yarn in the $i$th ply are established. Accordingly, each yarn will be partitioned into specific segments, as illustrated in Fig. \ref{fig:yarn-shape-and-insert}, where the thickness of yarn cracklets, $t_f$, is enlarged for clearer view (typically it should be set small or close to zero, i.e. $ t_f \to 0^+$). For short yarn shapes, the computed number of inserted yarn cracklet may be zero, particularly near a corner, thus no partitioning on those yarns will occur. 
Empirically, the corners and ends of test specimens are usually constrained or clamped by boundary conditions or loading units, tensile fiber fractures (triggered by a cluster of fiber failures) typically occur during the final stages in tensile test, and in general, is only observed in plies with layup angles close to the tensile load direction, and thus, are far less likely to occur near a physical corner of the laminate. For laminates with certain stacking sequences (e.g., all ply angles are larger than $45^{\circ}$ in a tensile test along $0^{\circ}$), ultimate failures are only resulted from extensive accumulation and interaction between matrix cracks and delaminations, i.e., fiber failure may not occur at all throughout test. 

Note that the number and position of yarn cracklets inserted into a yarn depends not only on the shape and length of that yarn, but also on the user-defined effective yarn fracture distance for the $i$th ply, $l_i$, which determines the maximum possible density of yarn fractures. If $l_i$ is larger than the laminate length $L$, then no yarn cracklet will be inserted into any yarn. Therefore, it is important for the user to determine a proper value of $l_i$, such that the tensile dominated failure can be captured by the model. Given that the matrix crack spacing, $d_i$, which is the yarn width in the $i$th ply, is relatively much smaller than the laminate plane dimensions ($W,L$), most yarns ought to have large aspect ratios (long and narrow yarns). Thus, as long as $l_i$ is reasonably chosen, any yarn not traversing near a laminate corner should have been partitioned and pre-inserted with yarn cracklets, and in this way there will be no 'unbreakable yarns'. Also, should any yarn span the whole length of the laminate plane but is not inserted with sufficient yarn cracklets, warning messages will be sent to the user by the code indicating a need for potential modification, i.e. the yarn cracklet spacing, $l_i$, may need to be reduced. Of course, $l_i$ can be set to a very small number (e.g., smaller than effective load transfer distance in a statistical concern) to capture any potential or small local stress redistribution, though in this case computational cost would raise. Nevertheless, insertion of yarn cracklets can still be deactivated by the user for specific usage: For example, the model without yarn fractures could be an effective toolkit for studies involving variational analysis of matrix crack systems as mentioned in Section \ref{only-matrix-crack}, wherein fiber fracture is not a concern.

Moreover, the inserted cracklets should be tied to corresponding continuum yarn segments to form a connected assembly. For a surface-based tie constraint between two surfaces in Abaqus (as with the case for most FEM tools), one surface must be defined as a master surface while the other as a slave surface, so that the DOF of nodes on the slave surface are constrained to the DOF of the master surface.
Interface elements are usually applied in conjunction with CZM to represent the interactions between larger continuum parts, and thus they usually do not directly carry external mechanical loads. Therefore, for any contact pair in AGDM, the edges of the continuum yarn segments should be set as the master surface, while the edges of inserted cracklets (for the three types of discrete cracks) should be set as slave surfaces. Any external load or boundary conditions should be defined only on master surface/nodes, namely the nodes associated with surfaces of continuum yarns.
As shown in Fig. \ref{fig:yarn-shape-and-insert}, corresponding yarn segments are defined as master surfaces, while all necessary yarn cracklet edges are set as slave surfaces. 
Such surface-based constraints are commonly used in numerical simulations with virtually zero or small-thickness cohesive elements, and where the necessity of explicit definition is very important. 

\subsection{Partition of yarn-to-yarn interfaces, based on yarn cracklet insertion locations, to define end-to-end locations of matrix cracklets}\label{sec2.3:matrix_segmentation}

As emphasized in Section \ref{multi-discrete}, it is critical to have matching meshes at potential crack bifurcations and intersections to capture accurate displacement jumps across cracks and load transfer between them. Therefore, each yarn interface should be segmented according to all potential crack bifurcations that are physically reasonable, such that whenever a yarn failure (fiber crack) occurs and generates or interacts with a matrix crack, matching nodes at such locations are guaranteed.   

\begin{figure}[H]
	\centering
	\includegraphics[width=0.8\textwidth]{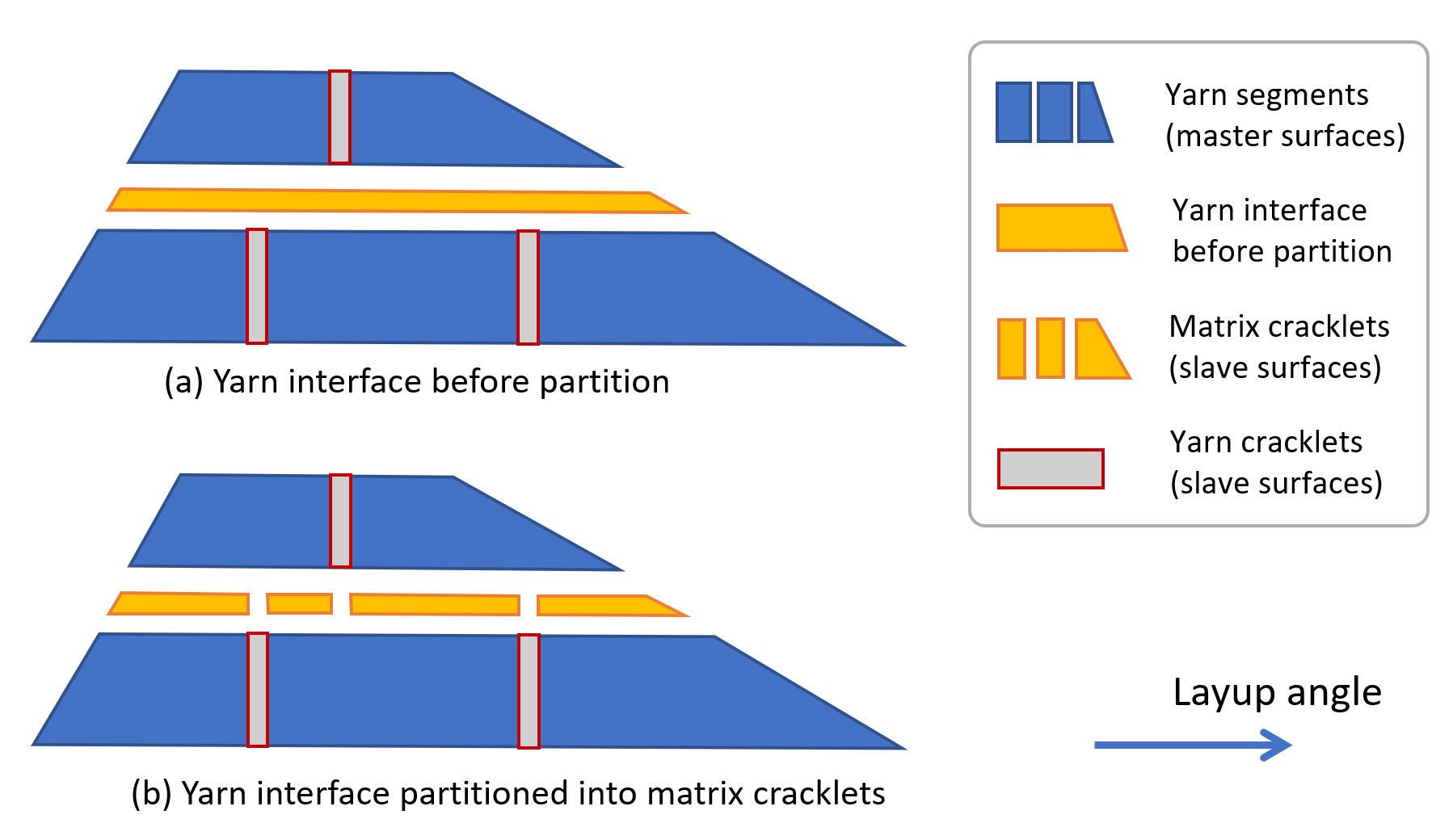}
	\caption{Partition of a yarn interface between yarns to form matching nodes at potential crack intersections (between axial yarn cracks and transverse matrix cracks).}
	\label{fig:matrix-crack-segmentation}
\end{figure}

To realize such a complex laminate cracking process, the preprocessing code stores the local and global coordinates of each yarn, as well as the position information of all inserted yarn cracklets in each yarn part. Then the segmentation of a yarn interface will be determined by the relative locations of inserted yarn cracklets in its two adjacent yarn neighbors. As shown in Fig. \ref{fig:matrix-crack-segmentation}, each interface between two yarns is partitioned (into matrix cracklets) to ensure matching nodes at all potential crack bifurcations. Again, master and slave surfaces must be defined explicitly in the model through use of a surface-based tie constraint formulations, as just explained in Section \ref{sec2.2:yarn_segmentation}.

\begin{figure}[H]
	\centering
	\includegraphics[width=0.8\textwidth]{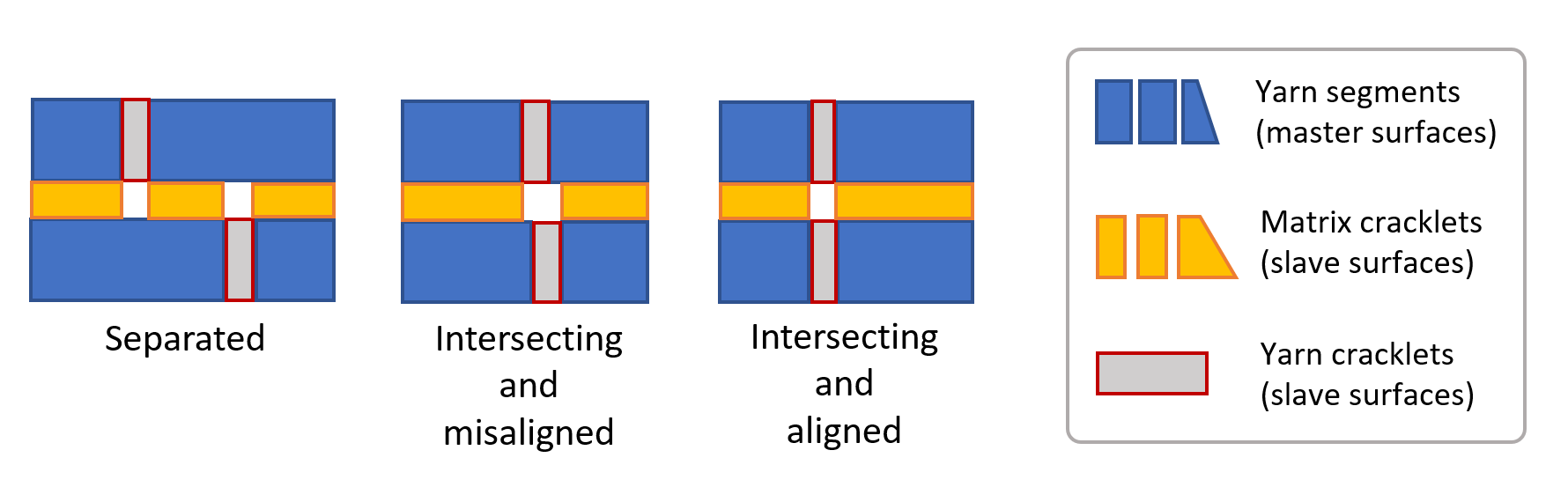}
	\caption{Possible alignment patterns between two yarn cracklets from adjacent yarns, wherein matching nodes are ensured in all cases.}
	\label{fig:fiber-fracture-alignment}
\end{figure}

Note that there are three possible alignment patterns between two yarn cracklets from adjacent yarns: separated (not intersecting), intersecting and perfectly aligned, intersecting and misaligned, as shown in Fig. \ref{fig:fiber-fracture-alignment}. The third case is common in $0^{\circ}$ and $\pm 90^{\circ}$ plies, whereas the second case is unlikely to occur since the value of $t_f$ is very small relative to the dimensions of other features. Nevertheless, all cases are checked by the preprocessing code to ensure explicit, and appropriate segmentation. Because of the geometric principals involved in partitioning yarn and interface parts, once a part has been partitioned into two separate sub-domains, four new nodes must be generated to form a partitioning line (or eight new nodes in a 3D element/part). However, only half of them may be required to form matching nodes on one side of a crack bifurcation, the remaining half of the nodes on the other side only help to generate sub-domains. Due to the tie formulation used, the nodes away from the crack tip, being slave nodes, will remain close to each other and follow the DOF of their corresponding master surface or edge, and thus, behave as one node. This way the material integrity in numerical model will not be affected given the typical small values of the cracklet thicknesses, $t_m$ and $t_f$.

\begin{figure}[H]
	\centering
	\includegraphics[width=0.8\textwidth]{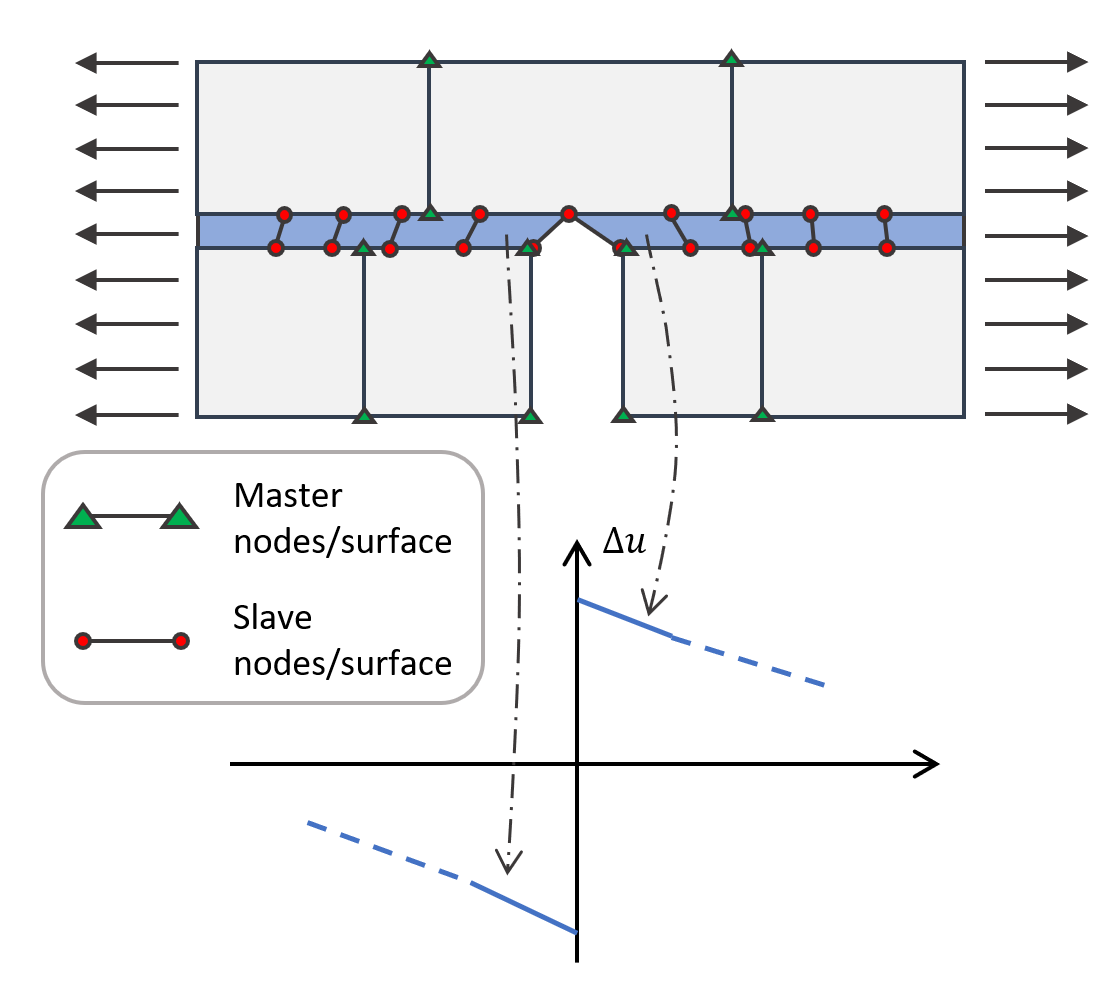}
	\caption{Accurate displacement jump is captured by matching nodes at the crack tip, while different parts are tied together by surface-based constraints (as occurs in AGDM). Nodes away from the crack bifurcation are not necessarily bonded, to compare with Fig. \ref{fig:comparison-of-matching-and-non-matching-ce} }
	\label{fig:matching-ce-tied}
\end{figure}

One important feature requires emphasis: As long as the nodes at potential crack bifurcations are matched properly, correct displacement jumps across crack tips and accurate load transfer between cracks can be explicitly modeled, meshed elements are not required to ensure such feature. Also, the nodes next to the crack tip are not necessarily bonded with other elements (i.e., do not need to share or bonded by common nodes). Instead, the nodes only a short distance away from the crack tip can be constrained to the adjacent edges using surface-based tie formulations. As demonstrated in Fig. \ref{fig:matching-ce-tied}, an accurate displacement jump can be captured explicitly by the matching nodes at the potential discontinuous fracture surfaces and crack bifurcations, while the other nodes are not necessarily bonded, as was the case shown in Fig. \ref{fig:matching-ce-bonded}.

The feature just described embodies a key fundamental difference between the proposed AGDM and the two current approaches discussed in Section \ref{multi-discrete} and summarized in Section \ref{intro-summary}. Both methods are based either on element partitioning or embedding, and thus they both require that bonded and structured mesh/elements be priorly-generated from the very beginning. In contrast, the presented method adopted in AGDM will create individual, precisely segmented domains/parts to ensure the matching of nodes at all potential crack bifurcations before mesh generation, which provides more flexibility on mesh techniques. 

For example, if the matrix crack density is set to 2 potential cracks per mm ($2 \text{/mm}$), then the ply element size must be $0.5 \text{mm}$ in both methods, whereas in AGDM, only the width of unmeshed yarn parts $d_i$ must be $0.5 \text{mm}$, while the mesh technique and yarn element size are completely controllable. Moreover, since the discretization is at the ply level, different intra-ply matrix crack densities are supported. These features and the flexibility afforded, make the model effective for various challenging applications involving stress redistribution in multi-ply laminates with matrix crack systems having various crack densities, as mentioned in Section \ref{only-matrix-crack}. Furthermore, the model supports arbitrary layup angles, $\theta_i$, and crack spacing, $d_i$, and the mesh size can be set small enough (in the whole model or at specified regions) to capture the stress profile and gradients of interest in fracture mechanics studies.

\begin{figure}[h]
	\centering
	\begin{subfigure}[h]{0.45\textwidth}
		\centering
		\includegraphics[width=\textwidth]{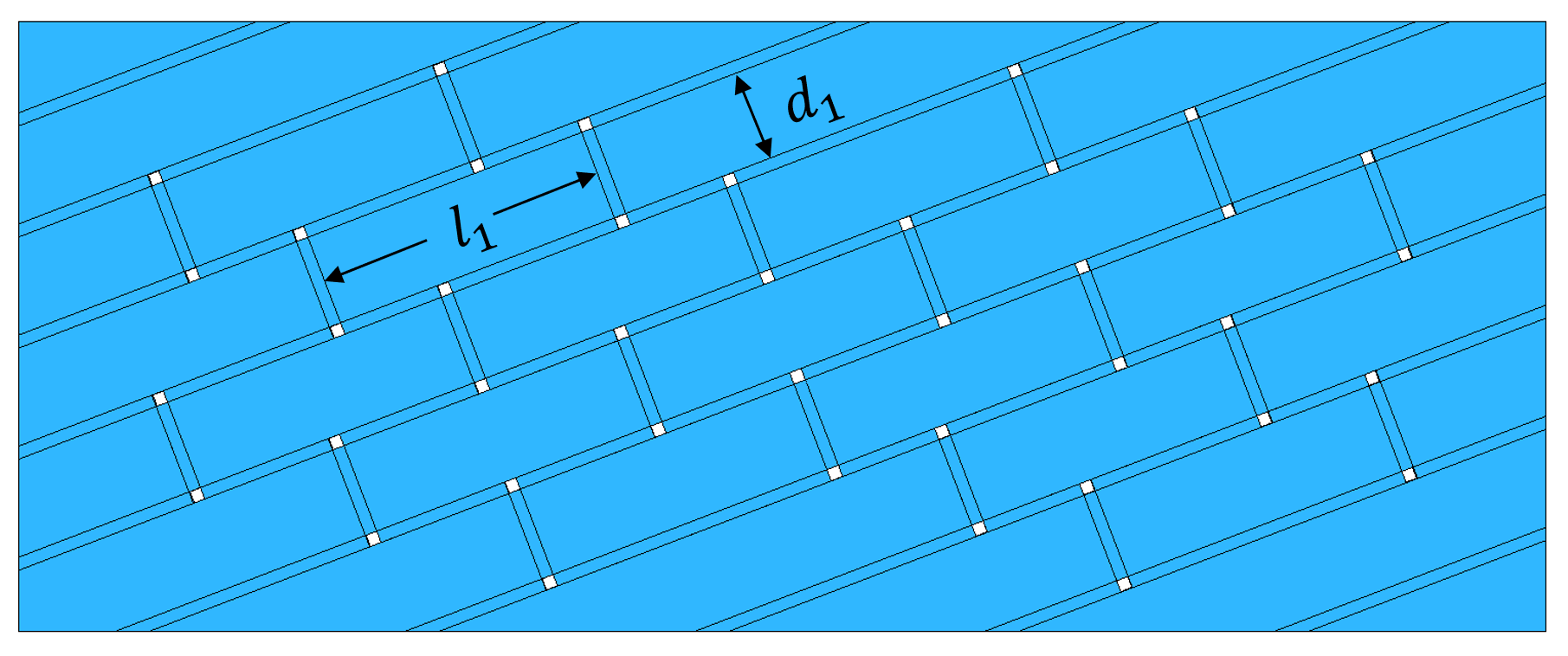}
		\caption{Ply assembly with enlarged interface thickness $t_f=t_m=0.2 \text{mm}$.}
		\label{fig:in-ply-part-example-1}
	\end{subfigure}
	\begin{subfigure}[h]{0.45\textwidth}
		\centering
		\includegraphics[width=\textwidth]{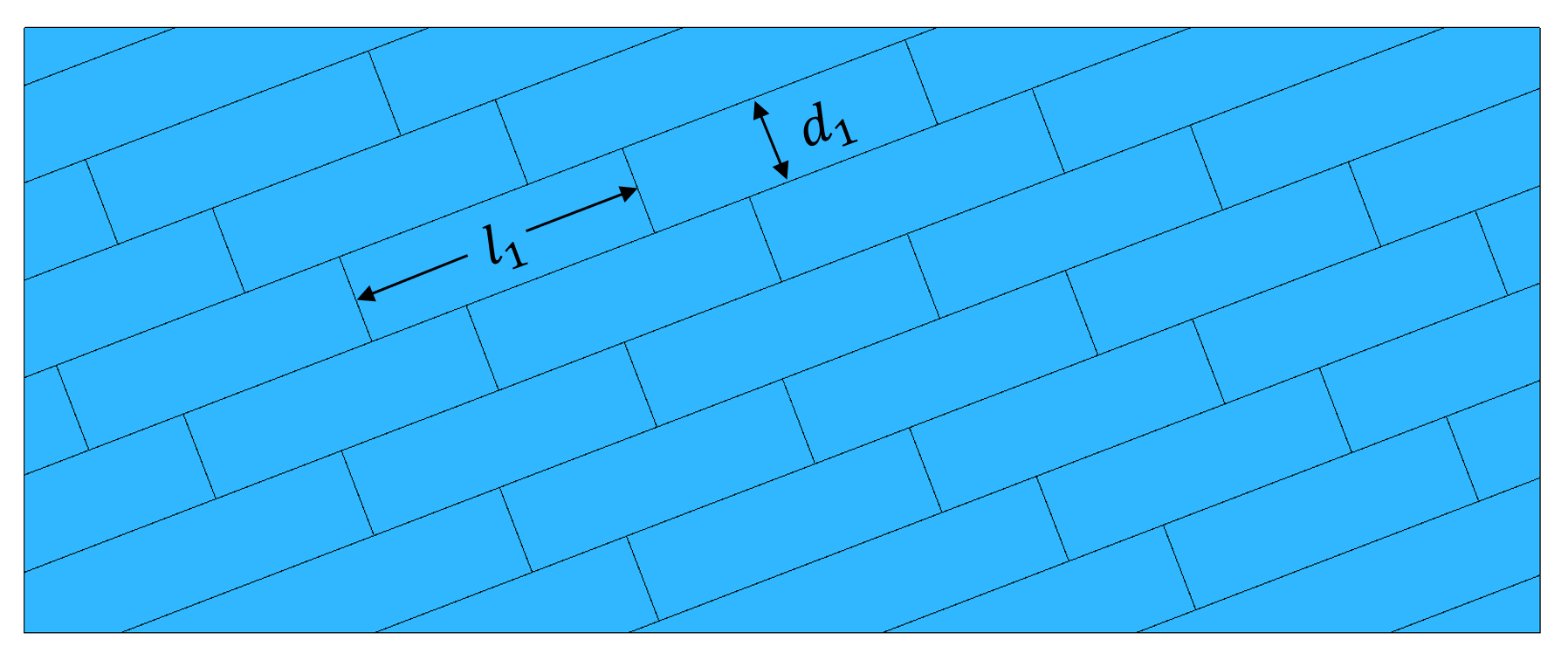}
		\caption{Ply assembly with typical small interface thickness ($t_f = t_m = 10^{-3} \text{mm} = 1 \mu\text{m}$).}
		\label{fig:in-ply-part-example-2}
	\end{subfigure}
	\begin{subfigure}[h]{0.45\textwidth}
	\centering
	\includegraphics[width=\textwidth]{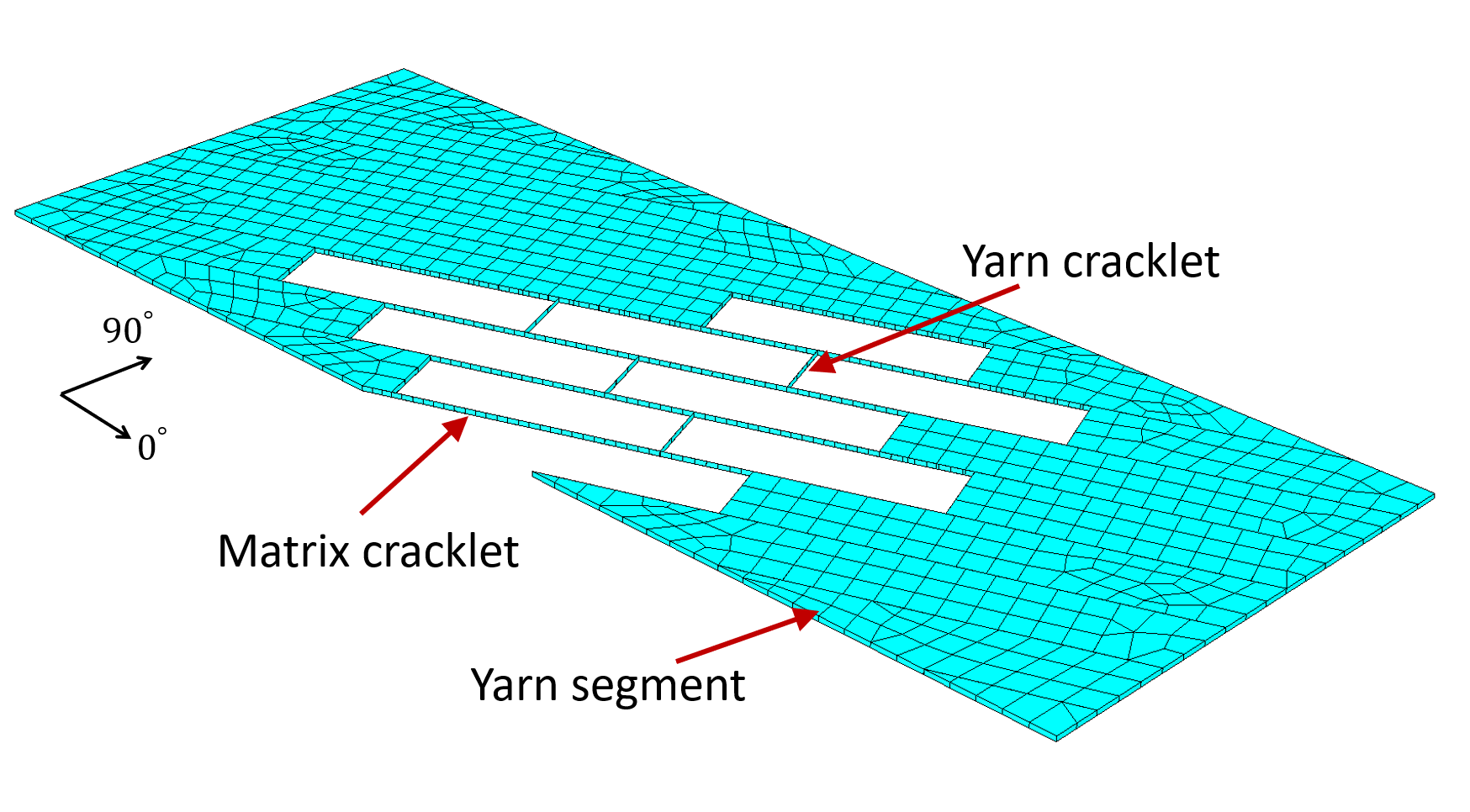}
	\caption{Mesh example of the ply with small-thickness interface ($t_f = t_m = 1 \mu\text{m}$) and with several yarn segments visually suppressed.}
	\label{fig:in-ply-mesh-example-1}
	\end{subfigure}
	\begin{subfigure}[h]{0.45\textwidth}
	\centering
	\includegraphics[width=\textwidth]{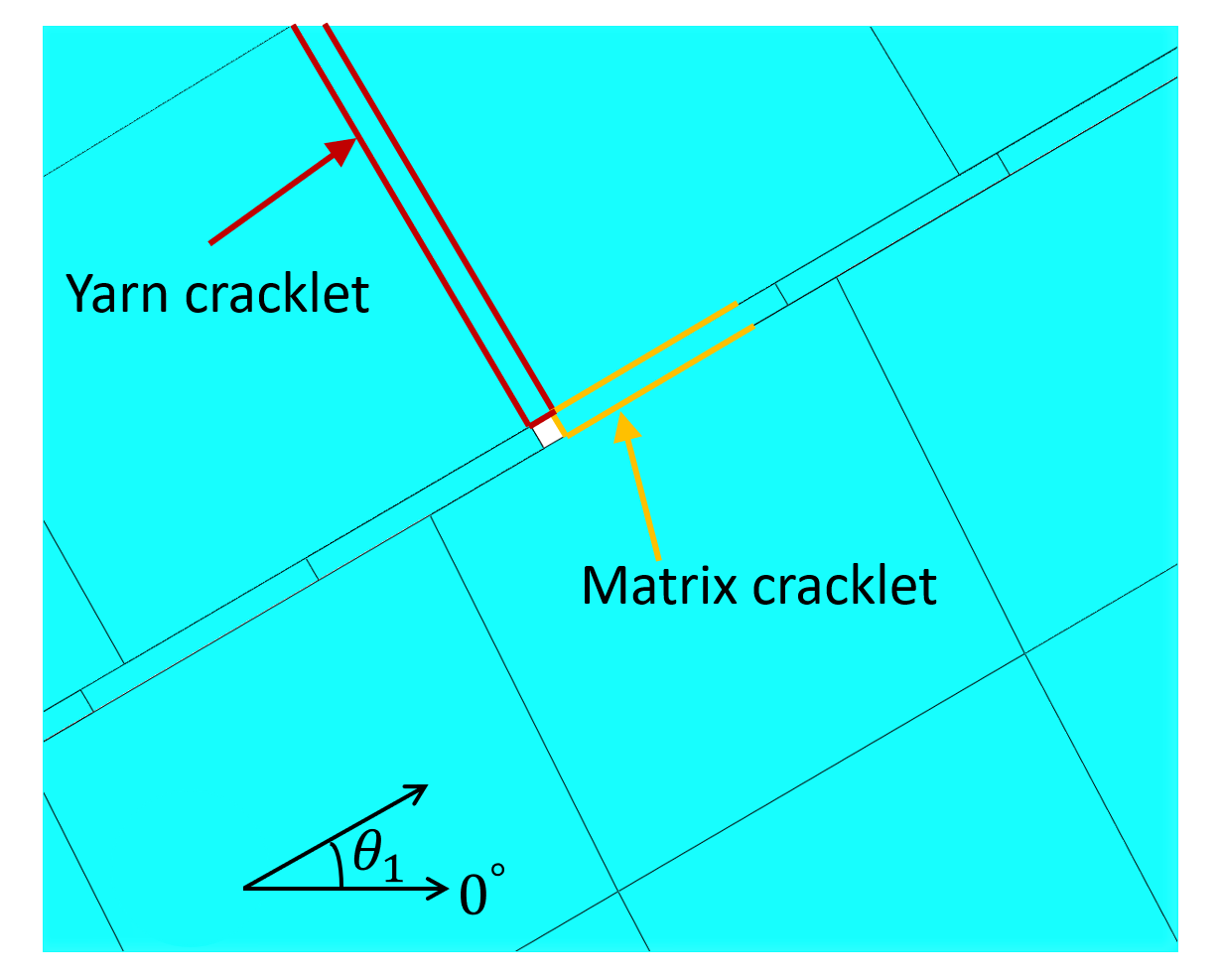}
	\caption{Magnified view of mesh at a potential intra-ply crack bifurcation.}
	\label{fig:in-ply-magnified-view-of-crack-tip}
	\end{subfigure}
	\caption{Assembly and mesh of a $ 25 \times 10 \times 0.1 \text{mm}$ ply with $\theta_1 = 21 ^{\circ}$ layup. Matrix crack spacing (yarn width) is $ d_1 = 1.5 \text{mm}$, and effective yarn fracture spacing is $ l_1 = 5.0 \text{mm}$.}
	\label{fig:in-ply-assembly-example}
\end{figure}

As an example, screen-shots of the assembled parts (yarn segments, yarn cracklets and matrix cracklets) and finite element mesh for a single ply, as generated by the preprocessing code, are shown in Fig. \ref{fig:in-ply-assembly-example}, and where the ply angle is arbitrarily chosen as $\theta_1 = 21^{\circ}$.
Fig. \ref{fig:in-ply-part-example-1} shows the parts assembled together to form the ply, where all parts have been translated on-site with necessary surface-based tie constraints assigned, and where interface thicknesses are enlarged ($t_f = t_m = 0.2 \text{mm}$) for clearer view.
Fig. \ref{fig:in-ply-part-example-2} shows the assembly of the same ply setup, but with more realistic, small-thickness interfaces ($t_f = t_m = 1.0 \mu \text{m}$), the extremely thin interface parts appear as thin single lines in this top-down view.
Fig. \ref{fig:in-ply-mesh-example-1} is a mesh view of the ply, where several yarn segments are visually suppressed to show the meshing of matrix cracklets and yarn cracklets, it is shown that the mesh size of the ply does not depend on the matrix crack spacing.
Fig. \ref{fig:in-ply-magnified-view-of-crack-tip} is a magnified view of the mesh at a potential crack bifurcation, and where matching nodes are realized as expected.

\subsection{Partition of an interface between two plies into delamination cracklets that accommodate previously-positioned, yarn and matrix cracklets }\label{sec2.4:delamination_segmentation}

Once yarn cracklets and matrix cracklets have been generated for all $N$ plies, an assembly of delamination cracklets is generated reflecting all possible inter-ply cracks that may occur. For a laminate with $N$ plies, there are $N-1$ ply interfaces to be partitioned into delamination cracklets, however, there will be differences in partitioning such interfaces, depending on the discretization and pattern of previous cracklet generation in the two adjoining plies.  

Similar to the yarn interfaces, each ply interface also must be partitioned into delamination cracklets that ensure matching nodes at all potential crack bifurcations, which is necessary to allow for formation of interconnected networks of cracks meandering through multiple plies. Delaminations can potentially interact with yarn fractures and matrix cracks from the two adjoining plies, therefore, the relative locations of all previously-inserted cracklets in adjacent plies must be taken into account in determining the locations and local geometries of all delamination cracklets. 

To review, the previously inserted cracklets in order of their insertion were yarn cracklets, followed by matrix cracklets (with segments accommodating the locations of the previous yarn cracklets), all of which are already stored in the preprocessing code. 
Then in partitioning each inter-ply interface part into delamination cracklets, the coordinates of all partitioning 'lines', which define the geometry of each delamination cracklet, must have end points and corners that are compatible with those previously-inserted intra-ply cracklets, in terms of allowing for all possible crack bifurcations to form a complex crack network. This requires linking of adjoining cracklets of any of the three cracklet types (irrespective of the order of their failure during loading). These locations and partitioning patterns are established based on geometric principles, and once determined, all surfaces of each segment will be designated as either master surfaces or slave surfaces, in order to assign suitable tie constraints. 

The procedural steps are follows:
In the first step, a complete rectangle of dimensions $W \times L $, representing an ply interface, is partitioned to form matching nodes at all potential bifurcation sites involving the delamination itself and matrix cracklets from two adjacent plies. In the second step, further partitioning is carried out to ensure that there are also matching nodes at all potential intersections with inserted yarn cracklets from the two adjacent plies. This way, matching nodes will be guaranteed at any potential crack bifurcations, regardless of the number and crack type (damage mode) involved.

\begin{figure}[H]
	\centering
	\begin{subfigure}[h]{0.45\textwidth}
		\centering
		\includegraphics[width=\textwidth]{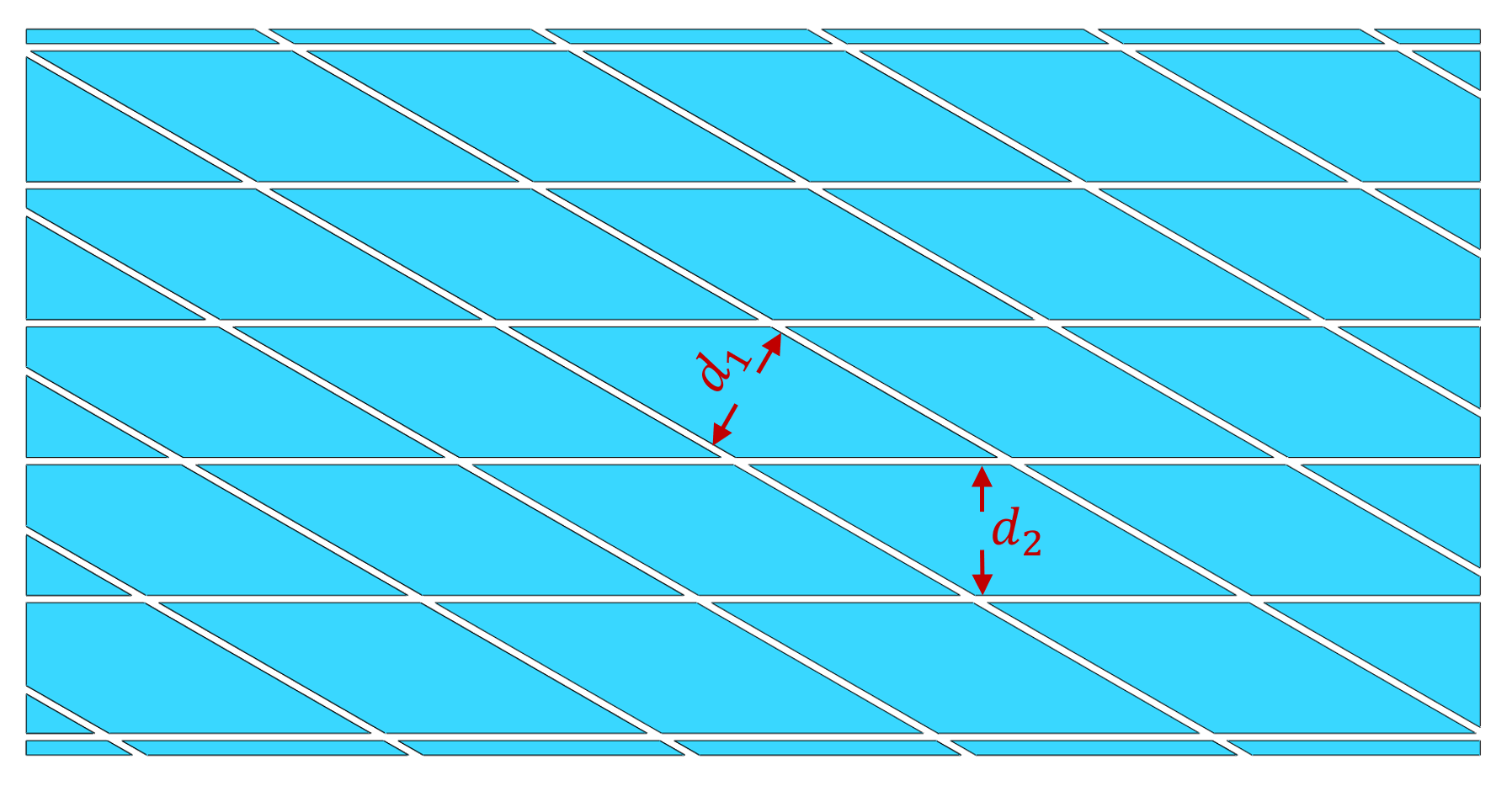}
		\caption{First-step partition of a ply interface part.}
		\label{fig:delamination-part-1}
	\end{subfigure}
	\begin{subfigure}[h]{0.45\textwidth}
		\centering
		\includegraphics[width=\textwidth]{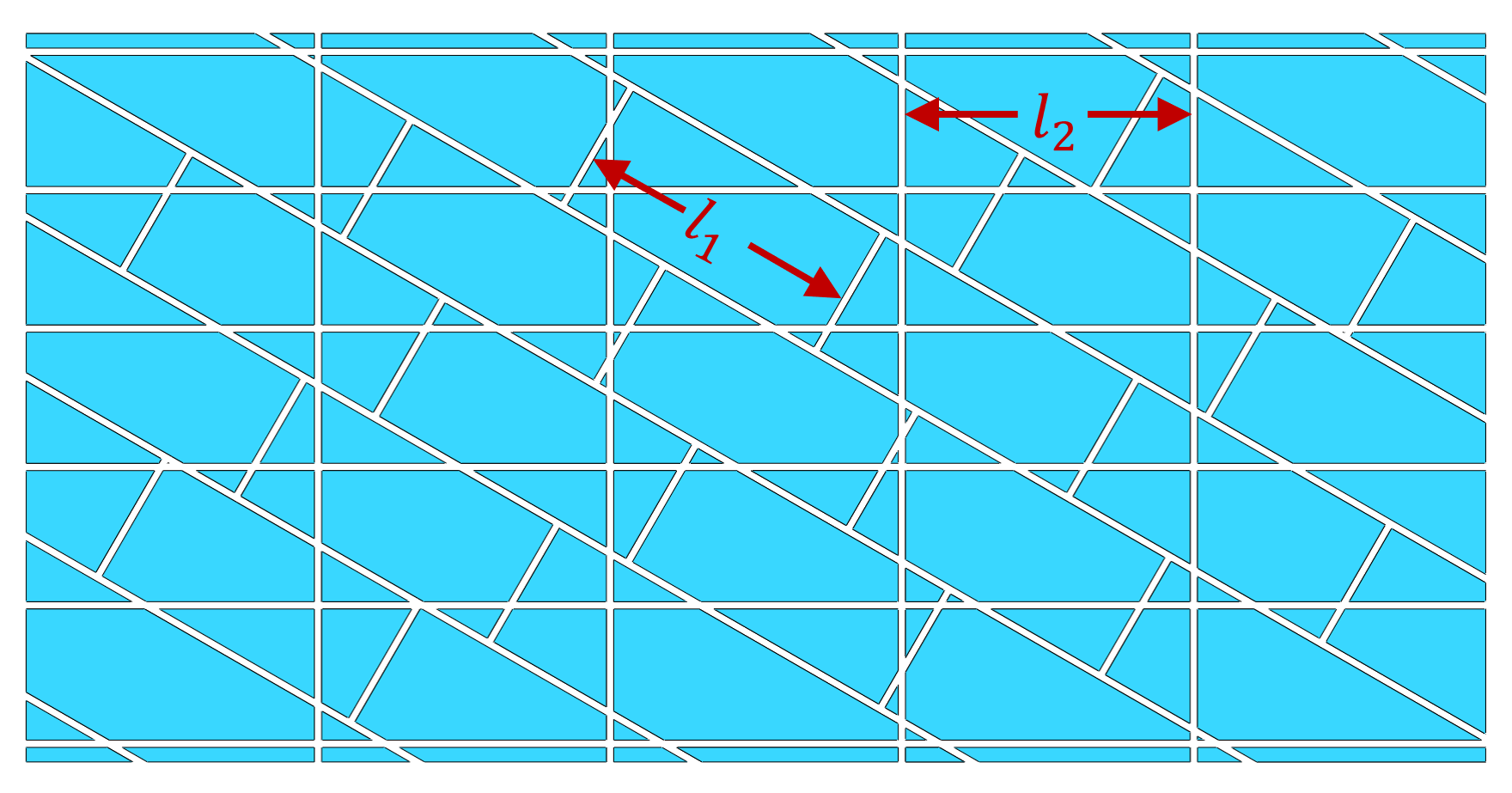}
		\caption{Second-step partition of the ply interface part in Fig. \ref{fig:delamination-part-1}.}
		\label{fig:delamination-part-2}
	\end{subfigure}
	\begin{subfigure}[h]{0.45\textwidth}
		\centering
		\includegraphics[width=\textwidth]{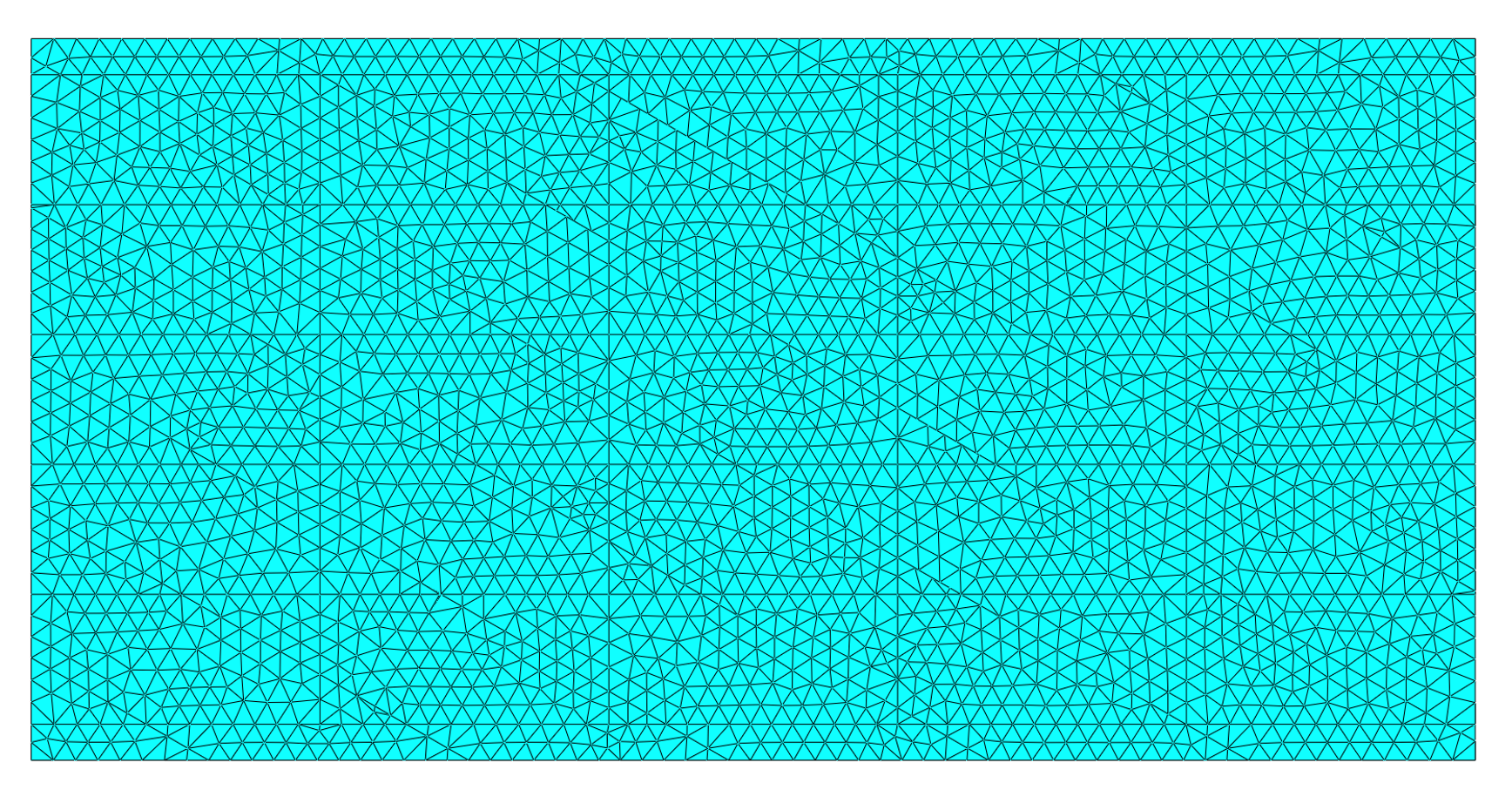}
		\caption{Mesh example of the partitioned ply interface part into delamination cracklets in Fig. \ref{fig:delamination-part-2} with small interface thickness ($t_m = t_f = 5\mu\text{m}$ )}
		\label{fig:delamination-mesh-1}
	\end{subfigure}
	\begin{subfigure}[h]{0.45\textwidth}
	\centering
	\includegraphics[width=\textwidth]{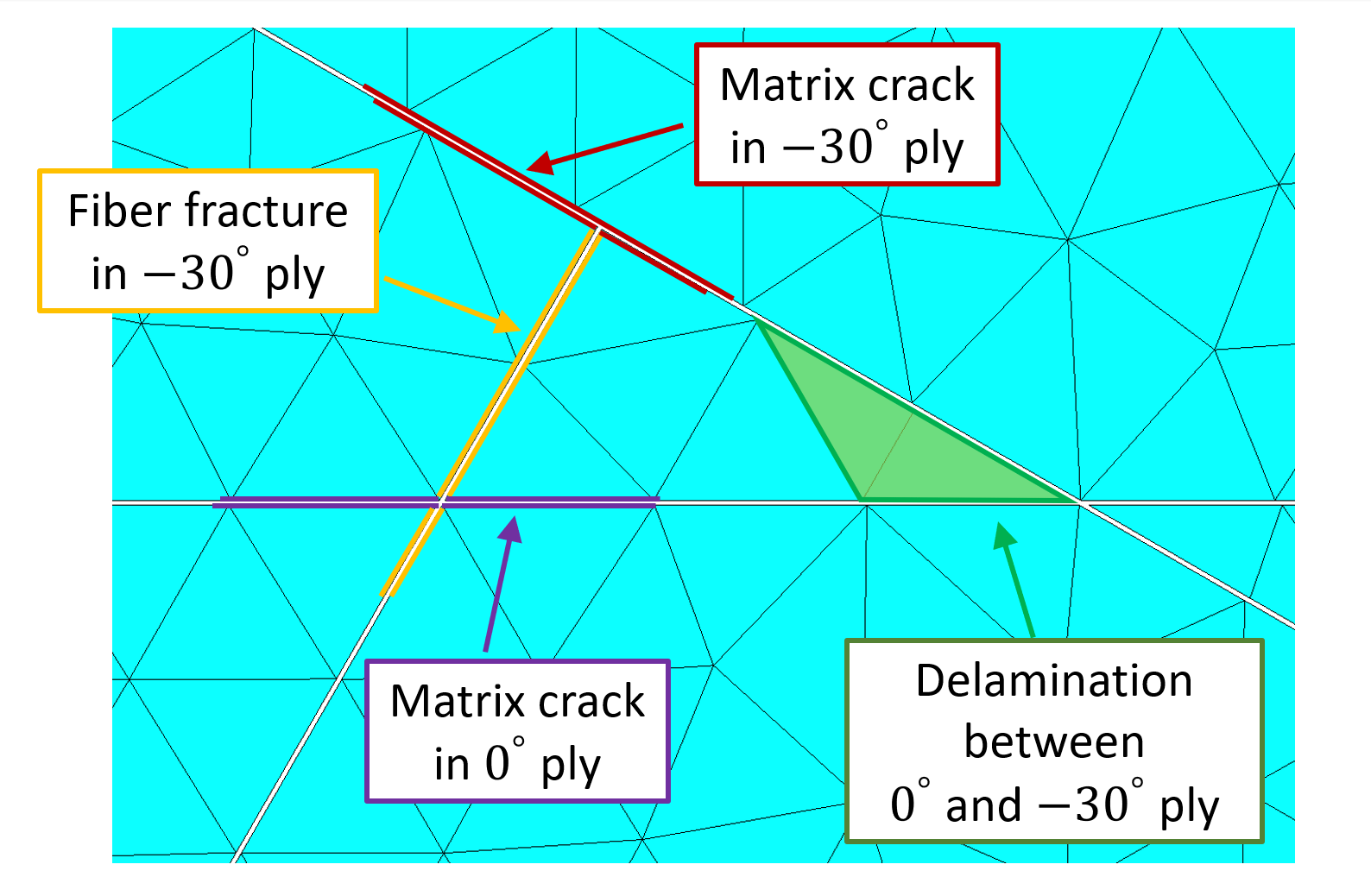}
	\caption{Magnified view of mesh in Fig. \ref{fig:delamination-mesh-1} to show the partitioning pattern of inter-ply delamination cracklets around matrix and yarn cracklets.}
	\label{fig:delamination-mesh-2}
	\end{subfigure}
	\caption{Partition and mesh example of a $ 20 \times 10 \, \text{mm} $ ply interface part between two adjacent plies with angles $ \theta_1 = -30 ^{\circ} $ and $ \theta_2 = 0 ^{\circ} $. Matrix crack densities are $ d_1=d_2= 1.8\text{mm}$, and yarn fracture spacings are $ l_1 = l_2 = 4.0\text{mm}$}
	\label{fig:delamination-part-figs}
\end{figure}

As an example, Fig. \ref{fig:delamination-part-1} shows the first-step partitioning of a inter-ply part between two adjacent plies to form matching nodes with potential matrix cracks, at $ 0^{\circ} $ and $ -30^{\circ} $, respectively, and where the intra-ply yarn-to-yarn interface thickness (matrix cracklet thickness) has been exaggerated ($ t_m = 0.1 \text{mm}$) for easier visualization of the partitioning pattern.
Fig. \ref{fig:delamination-part-2} shows the second-step partitioning of the ply interface in Fig. \ref{fig:delamination-part-1} to form matching nodes with potential yarn fractures (cracklets), again with exaggerated interface thickness ($t_b = 0.1 \text{mm}$);
Fig \ref{fig:delamination-mesh-1} shows a meshing example of the partitioned ply interface in Fig. \ref{fig:delamination-part-2} but with more realistic matrix and yarn cracklet thicknesses, $ t_f = t_m = 5 \mu \text{m}$, wherein the partitioning profile becomes very difficult to see. 
Fig. \ref{fig:delamination-mesh-2} shows a magnified view of the mesh in \ref{fig:delamination-mesh-1} at a potential crack bifurcation, and as expected, matching nodes are ensured for all potential crack bifurcations regardless of the type and pattern of cracks.

Once the partitioning process has been completed for all interface parts in the 2D projection plane (from a top-down view), the actual ply interface parts, partitioned into delamination cracklets, are generated in 3D modeling space by extruding their planar rendering with magnitude of the defined ply interface thicknesses, $t_d$. Then all cracklets for a specific interface part are translated and positioned to fill the 'thin gaps' between two specific plies in the global coordinate system of laminate model. After which corresponding surface-based tie constraints will be assigned.

To improve computational efficiency and accuracy, the code also checks the span and area of each resulting delamination cracklet and suppress any 'small block' with extremely small span or projected area. Such small segments are formed when at least three partitioning lines (representing yarn fractures or matrix crack surfaces from adjacent plies) 'almost intersect' at one point. According to conventional finite element theory, such small segments would typically result in extremely small elements, thus potentially introducing unnecessary residual error. Therefore, an examination process for such segments is incorporated into the code to ensure mesh quality. The user must input a threshold value such that a domain/segment will be removed from a partitioned ply interface if its area is below a certain threshold, for instance $10^{-9} \, \text{mm}^2$.

\subsection{Assembly of laminate from all yarn segments and cracklets (yarn, matrix, delamination) and assignment of surface-based tie constraints}\label{sec2.5:assmebly}

As explained in the above subsections, once a continuum or interface part (i.e., a yarn segment, or cracklet of three possible types) has been generated, it will be directly translated and positioned onto its specific location in the global coordinate system of the laminate. Thus, it remains to assign proper surface constraints between all adjacent parts to form an interconnected assembly constituting the entire laminate.

In order for the code to automatically assign all required surface-based tie formulations, each part and its surfaces must be explicitly labeled based on its global assembly location and geometry, such that the required number of contact pairs and two surfaces involved can be determined based on the ply order in the laminate and the ply-level discretization. For example, a matrix cracklet in AGDM may be labeled as 'layer-$i$-yarn-interface-$j$-segment-$k$-face-y-positive', which refers to the surface facing towards the positive $y$ direction in the $k$th partitioned segment of the $j$th yarn interface in the $i$th ply. Accordingly, it is understood to be the slave surface being tied to a master surface of a specific segment of a partitioned yarn (to be determined based on that yarn segment's label). A similar labeling mechanism in AGDM applies to all continuum yarn segments and all types (yarn, matrix and delamination) of interface cracklets. For brevity details are omitted here. In Abaqus, there are two available surfaced-based tie formulations, a surface-to-surface constraint and a node-to-surface constraint. The surface-to-surface tie constraint is the formulation that provides for optimized stress accuracy in both Abaqus/Standard and Abaqus/Explicit. \cite{DassaultSystemes2016}, and thus it is adopted in AGDM.  

\begin{figure}[H]
	\centering
	\includegraphics[width=0.8\textwidth]{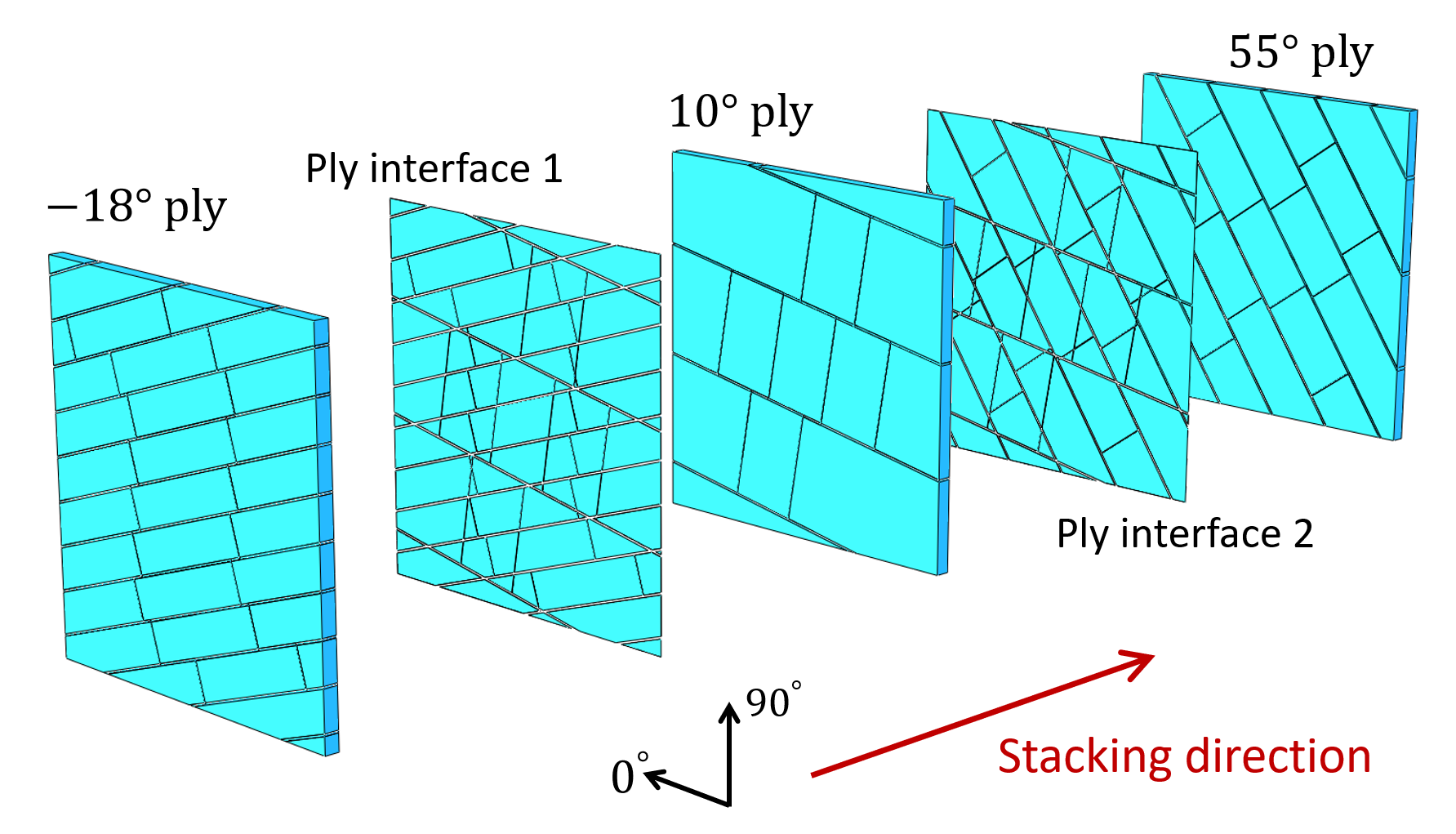}
	\caption{Expanded view of a three-ply laminate of an arbitrarily chosen layup $[\text{--}18^{\circ}/10^{\circ}/55^{\circ}]$, which has been reassembled from various laminate 'parts' that were subdivided into various yarn segments, yarn cracklets, matrix cracklets and delamination cracklets, and whose geometries were explicitly generated in AGDM. These resulting pieces were then reassembled to generate the model for laminate failure when subjected to loading in time, whereby the damage process evolves as a complex network of cracks leading to ultimate failure as the laminate breaks apart. Model parameters: $W=L=5 \text{mm}$; $d_{1,2,3}=0.5,1.5,0.8 \text{mm}$; $l_{1,2,3}=2,1,1.5 \text{mm}$; $h_{1,2,3}=0.2,0.15,0.2 \text{mm}$; $t_m=0.03 \text{mm}$; $t_f=0.01 \text{mm}$; $t_d=5e^{-3} \text{mm}$.}
	\label{fig:assembly-example-1}
\end{figure}

As an example, Fig. \ref{fig:assembly-example-1} shows the assembly of a square-shaped three-ply laminate with an arbitrarily chosen layup and dimensions. For illustration clarity, parts from each ply and associated ply interfaces are expanded in the laminate stacking direction.

\subsection{Materials and Mesh}

Since AGDM is geometry-based, all essential procedures (part generation and partitioning, laminate model assembly, and assignment of all surface-based tie constraints) do not require any material definition or predetermined mesh. Consequently, the material assignment and mesh technique are flexible once AGDM creates the laminate model. 
The user can either employ an built-in material library in Abaqus, or apply a user-defined material (such as UMAT, VUMAT), or user-defined elements (UEL) in conjunction with AGDM. 
The user can also determine the desired mesh technique and level of refinement appropriate to a specific application. The nodes and surfaces of each ply and the entire model are also grouped into node and surface sets, so that the user can directly apply specific boundary conditions or keep a track of outputs of a specific group. 

\subsection{Summary of fundamentals of AGDM and possible directions for further enhancements}

In this section, we have described the basic process and methodologies used in the proposed Auto-generated Geometry-based Discrete model (AGDM), and how it is realized in the preprocessing code. AGDM and several other notable modeling techniques are briefly summarized in Table \ref{tab:comparison-between-studies}. As shown, the AGDM approach realizes full capability in accurately modeling all three types of discrete cracks (with physical fracture surfaces) and their multi-mode couplings in composite laminates with an arbitrary stacking sequence, meanwhile avoids dependency of mesh sizes on crack spacings. All of these are key features required for effective and high-fidelity modeling of damage evolution in composite laminates. 

\begin{table}[H]
	\caption{Comparison of methods and capabilities of recent representative works on laminate progressive damage modeling} \label{tab:comparison-between-studies}
	\small
	\begin{tabularx}{\textwidth}{L{0.8} L{1.2} L{1.2} L{0.9} L{0.9}}
		\hline 
		Authors & Method and types of failures included & Simultaneous and accurate coupling of cracks & Arbitrary layup & Mesh dependency \\[1ex]
		\hline \hline
		Sun et al. \cite{Sun2013}; Hu et al. \cite{Hu2016}; Higuchi et al. \cite{Higuchi2017} & XFEM-based method for one type of intra-ply failure mode, cohesive elements with specially enriched nodes for delaminations & Interactions between cracks are approximated by modified enrichment functions, real intra-ply fracture surfaces still absent due XFEM-based approach & Capable & Crack spacing depends on mesh size \\[1ex] 
		\hline
		Higuchi et al. \cite{Higuchi2017a} & Inserted cohesive elements for matrix cracks and delaminations & Non-matching mesh between intra- and inter-ply cohesive elements, fiber failures are CDM-based without real fracture surfaces & Limited choices because the intra-ply cohesive elements shares nodes with continuum mesh & Crack spacing depends on mesh size \\[1ex]
		\hline
		Joseph et al. \cite{Joseph2018} & Modified Crack Band model for one type of intra-ply failure mode, cohesive elements for delaminations & Delamination is trigged by matrix cracks to account for their interaction, but intra-ply fracture surfaces are absent due to CDM-based crack band model & Capable & Reduced mesh dependency by artificially track and control crack propagation to maintain crack spacing \\[1ex]
		\hline
		Lu et al. \cite{Lu2018} & Element partitioning based on FNM for one type of intra-ply crack, SCE based on FNM for delaminations & SCE is only capable of matching at most two cracks (one from each adjacent ply) & Capable & Crack spacing depends on mesh size, could be reduced by artificial modification \\[1ex]
		\hline  
		Joosten et al. \cite{Joosten2018} & Predefined interface elements in priorly-generated structured mesh for matrix cracks, fiber failures, and delaminations & Capable & Limited ply-angle choices since highly skewed meshes are inevitable for arbitrary layups & Crack spacing depends on mesh size \\[1ex] 
		\hline
		Bouvet et al. \cite{Bouvet2009, Bouvet2013} & Predefined interface 'spring 'elements in priorly-generated mesh for matrix cracks and delaminations & Fiber failures are CDM-based without discontinuous fracture surfaces & Limited ply-angle choices because overlapping nodes between plies are required to define interface elements & Crack spacing depends on mesh size  \\[1ex]
		\hline
		Liu and Phoenix (present study) & Generation and assembly of yarn segments and interface cracklets for matrix cracks, yarn fractures and delaminations & Capable & Capable & Crack spacing can be different among plies, independent of mesh size \\
		\hline
	\end{tabularx}
\end{table}

Beyond the generation and assembly of the laminate model, the preprocessing code also integrates some essential steps required for a complete numerical analysis (though they are independent of AGDM and can be flexible), for example, selection or input of material properties and constitutive laws, specific mesh size and technique, desired loading and boundary conditions, and analysis job submissions. This significantly increases modeling efficiency, and is suitable for potential future applications involving virtual testing and parametric studies of laminates.

Regarding AGDM's functions, development of additional features are foreseen and essential, such as open hole/notched or round-shaped laminates, and compressive failure modes under torsion and bending. Another interesting aspect would be the statistical strength of materials and interfaces: for instance, the modeling of discrete yarn failures that have statistically varying material strengths from one to another, which is an essential feature of the composite damage processes at all scales, since tensile failure of carbon fibers and their epoxy-impregnated yarns is brittle and dominated by randomly occurring 'flaws' whose strengths typically follow Weibull, weakest-link statistics resulting in size effects, as explained in Phoenix et al. \cite{Phoenix1978, Phoenix1997, Mahesh2002}.

\section{Numerical Simulation Example and Comparison with Experimental Failure Process in Laminate Test Specimens}\label{sec3:numerical_example}

This section presents some numerical cases of failure process in laminates to compare the resulting behavior with that seen in comparable tests from an experimental study by Johnson and Chang \cite{Johnson2001}, in which strength data and X-radiograph of failure evolution of laminates with various layups were documented. Firstly, specimen geometries and boundary conditions, material properties and treatment of parameters for cohesive zone elements are discussed. Then, some numerical cases of an un-notched $[30/90/\text{--}30]_s$ laminate were studied to evaluate input parameters and mesh convergence behavior. Finally, with the selected input parameters based on the above cases, eight more numerical examples of laminates with other layups were investigated.  

\subsection{Geometry and boundary conditions}\label{geometry-and-bc}

\begin{figure}[H]
	\centering
	\begin{subfigure}[h]{0.8\textwidth}
		\includegraphics[width=\textwidth]{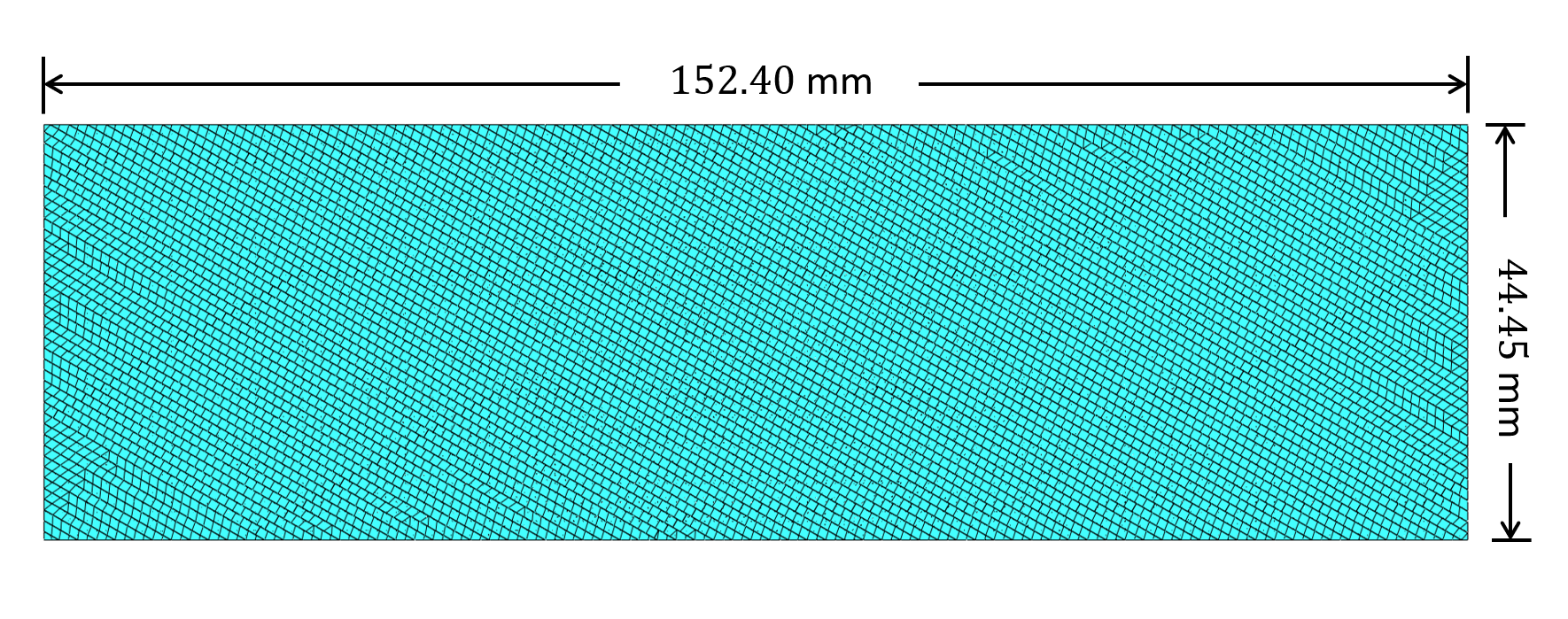}
		\caption{Top view of mesh and dimensions of the $[30/90/\text{--}30]_s$ laminate, where interfaces (cracklets) are not visible due to their small thickness (almost 2D dimensions).}
		\label{fig:model-mesh}
	\end{subfigure}
	\begin{subfigure}[h]{0.8\textwidth}
		\centering
		\includegraphics[width=\textwidth]{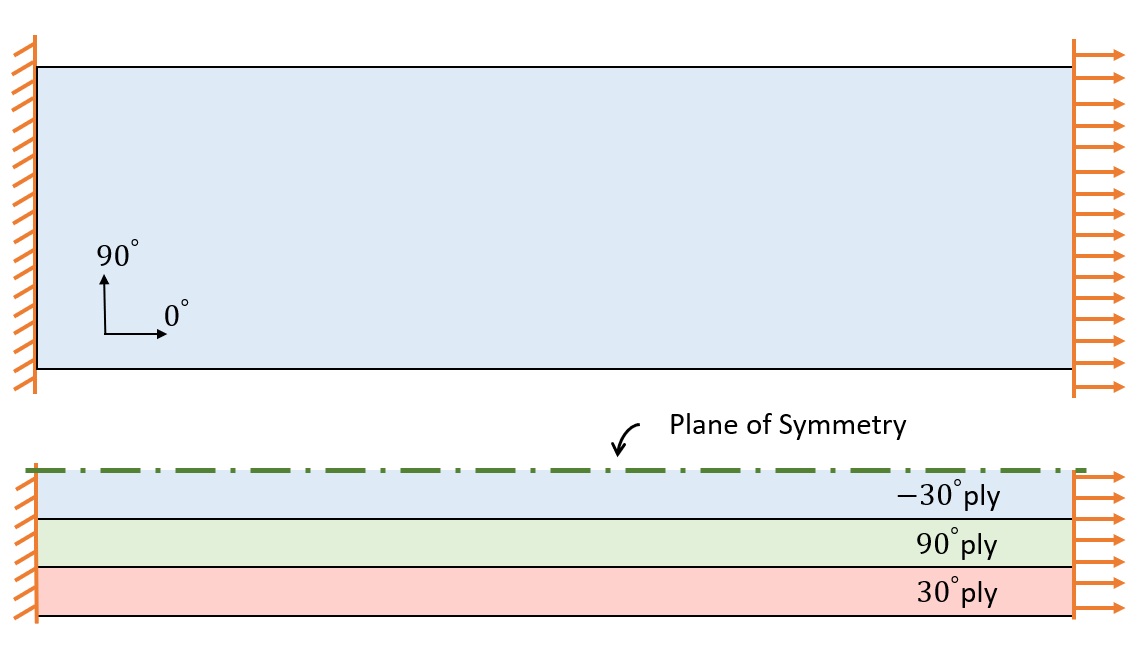}
		\caption{Boundary conditions applied in the numerical models ($[30/90/\text{--}30]_s$ laminate shown here), the ply thicknesses are artificially enlarged for schematic illustration. }
		\label{fig:model-bc}
	\end{subfigure}
	\caption{Dimensions and views of the boundary conditions and loading of the laminate in numerical simulations.}
\end{figure}

In the experimental study, all specimens have symmetric layup \cite{Johnson2001}, therefore by taking advantage of symmetry, only the bottom half of plies were modeled. For example, the layup of laminate in first series of numerical cases is $[30/90/\text{--}30]_s$, then the inputs for stacking sequence are: $N=3$, $\theta_1 = 30^{\circ}, \theta_2 = 90^{\circ}, \theta_3 = -30^{\circ}$. Correlated symmetric boundary conditions will be applied to the central plane. As an example, Fig. \ref{fig:model-mesh} shows dimensions, boundary conditions and stacking layup of the $[30/90/\text{--}30]_s$ specimen. 

All specimens are rectangular-shaped laminate with $W \times L = 44.45 \text{mm} \times 152.40 \text{mm}$, and the ply thicknesses are $h_{i}=0.129 \text{mm}$. Matrix crack spacing is set to $d_i = 1.0 \text{mm}$ for all plies, though they are not necessarily the same among plies. The effective spacing of yarn cracklets is $l_{i} = 25 \text{mm}$, and it could be deactivated in large angle (such as $90^{\circ}$) plies, realizing that yarn fractures are not expected in a $90^{\circ}$ ply under tensile load in the $0^{\circ}$ direction. Thickness of all interfaces (cracklets) is set to $t_f=t_m=t_d=10^{-3} \text{mm} = 1 \mu\text{m}$. For reference, with the above inputs and continuum yarn mesh size of $1 \text{mm}$ for all plies of $[30/90/\text{--}30]_s$, the resulting finite element model generated by AGDM contains $1325$ instances, $1886$ surface-based tie constraints, and $67414$ elements. For this specific case, it took the preprocessing commands about $2$ hours to generate the entire model and assign all local coordinates, material properties, and tie constraints. 

It should be noted that the mesh sizes can vary between different plies/yarns and are independent of matrix crack spacings $d_i$, as will be shown in the mesh and parametric study of some numerical examples. Besides, the chosen values of $d_i$ are based on commonly used inputs in comparable studies, and with the concern of computational cost. In real situations, much lower minimum crack spacings may be expected, as they are typically in the order of ply thickness. Given the large in-plane dimensions of the specimen, $d_i \approx 1\text{mm}$ should provide enough cracklets to capture effective high-density matrix crack networks. The continuum yarn mesh sizes, independent of matrix crack spacings, may be set to small values as well, such that high-quality stress profile and load transfer among cracks can be obtained. Of course, case studies on these features are worth further investigation, especially when local stress gradient and profile are of interest.

As shown in Fig. \ref{fig:model-bc}, the boundary conditions of specimens include: a fixed edge imposed along one end, a displacement loading in the $0^{\circ}$ direction applied on the other end, one free surface, and one plane of symmetry.

Considering the large number of required surface-based tie formulas and the type of material damage that evolves, Abaqus/Explicit solver is adopted for the numerical analysis.
Semi-automatic mass scaling throughout the loading steps are used with a target increment time of $10^{-6}$ seconds and a scaling frequency of every $10^4$ steps. For the $[30/90/\text{--}30]_s$ specimen as example, a total displacement of $1.2 \text{mm}$ is applied over $4$ seconds, while the loading rate is initially $0.0 \text{mm/s}$, then gradually increased to $0.01 \text{mm/s}$ at $0.1$ second, then raised to $0.3 \text{mm/s}$. The displacement loading increments have transient 'stops' every $0.1$ seconds to reduce any slight kinematic energy, such kind of 'stops' between loading were also applied in real experimental tests \cite{Johnson2001}. The modified loading curve can be defined as 'smooth steps' with specific controlling points in Abaqus \cite{DassaultSystemes2016}. The selected loading profile and associated mass scaling parameters were carefully determined based on various case studies, such that the kinetic energy would always be way less than $5\%$ of the strain energy until the point of complete failure. This way, the numerical simulation can be considered quasi-static.

\subsection{Materials and failure criterion}\label{subsec:materials}

The plies of laminate specimens in the experiment study are made by a commercially available product, T300/976 graphite/epoxy \cite{Torayca}. The material properties and fracture parameters associated with the three failure modes are described as follows.

\subsubsection{Continuum yarn segments and yarn cracklets}

The composite ply properties of orthotropic, linearly-elastic T300/976 are shown in Table \ref{tab:t300/976_yarn}, and are assigned to all yarn segments and yarn cracklets in their local coordinates based on the ply layup angle. First order, fully-integrated eight-node brick elements (C3D8), and 6-node triangular prism elements (C3D6) are used for continuum yarn segments and yarn cracklets. The properties of continuum yarns are available from built-in material library in Abaqus.

No stiffness degradation or failure occurs in the yarn segments. This way, the continuum yarn elements, which occupy almost the entire laminate volumetric space, remain orthotropic and linearly elastic throughout the numerical analysis. Yarn fractures occur as failed yarn cracklets, thus retaining the features of discrete cracks with real fracture surfaces. Specifically, for material points in a yarn cracklet, local damage initiation is triggered by excessive tensile stress in the fiber direction, based on the following criteria:
\begin{equation}\label{eq:fiber_failure}
\frac{\langle \sigma_{11} \rangle}{S_{11}} = 1
\end{equation}
where $\sigma_{11}$ is the current stress in the yarn/fiber direction and the symbol, $\langle \, \rangle$, represents the Macaulay bracket:
\begin{equation}
\langle x \rangle=\max(0, x)
\end{equation}
Upon damage initiation (once Eq. \ref{eq:fiber_failure} is satisfied for a material point), damage evolves such that the material stiffness is degraded based on a damage variable, $D^{f}$, defined as:
\begin{equation}\label{eq:damage_variable}
D^{f} = \frac{\epsilon_{11}^u(\epsilon_{11} - \epsilon_{11}^0)  }{\epsilon_{11}(\epsilon_{11}^{u}-\epsilon_{11}^0)}
\end{equation}
where $\epsilon_{11}, \epsilon_{11}^u$ and $\epsilon_{12}^0$ are current strain, ultimate failure strain and damage-initiating strain in the fiber direction, respectively. Damage variables range from $0$ to $1$ where '$0$' refers to an intact material and '$1$' represents its complete failure, reflected by the initial stiffness being multiplied by $(1-D^{f})$. 

It should be noted that the constitutive properties and failure behavior for yarn cracklet elements are realized by user-defined material (requires VUMAT file written in Fortran commands), whose stiffness values are equivalent to the continuum yarn elements, however, due to the small thickness of yarn cracklets, the strain of a yarn cracklet element in fiber direction ($\epsilon_{11}$) is regarded as the element's real separation, i.e., the corresponding 'strain' is computed as if the element has a thickness of one unit ($1 \text{mm}$ as in the unit system of this study). This is a commonly seen treatment in numerical models containing zero or small-thickness interface elements (for zero thickness elements, strain is not even practical due to 'division by zero'). Indeed, built-in cohesive elements (COH3D8, COH3D6) in Abaqus, by default, will have a nominal thickness of one, as is the case for matrix and delamination cracklet elements to occur in the next subsection. Here, the yarn cracklets are modeled with C3D8 elements, but the computation of $\epsilon_{11}$ in VUMAT is adopted with common treatment for small-thickness elements.

\begin{table}[H]
	\centering
	\caption{Mechanical properties of T300/976 plies and their epoxy impregnated yarns \cite{Johnson2001a, Lu2018, Wang2009, Hosseini-Toudeshky2013, Torayca}. Note that yarn cracklet elements are modeled with user-defined materials (VUMAT) wherein a 'nominal thickness' of one unit is adopted, namely, 'strains' are real interfacial separations in fiber direction.} \label{tab:t300/976_yarn}
	\begin{tabular}{l l}
		\hline 
		Property & Value \\
		\hline \hline
		Density: $\rho$ ($\text{g/cm}^3$) & 1.76 \\
		Longitudinal Young's Modulus: $E_{11}$ (GPa) & 139.2 \\
		Transverse Young's Moduli: $E_{22}$, $E_{33}$ (GPa) & 9.72 \\
		Shear Moduli: $G_{12}$, $G_{13}$ (GPa) & 5.58 \\
		Shear Modulus: $G_{23}$ (GPa) & 3.45 \\
		Poisson's ratios: $\nu_{12}$, $\nu_{13}$ & 0.29 \\
		Poisson's ratio: $\nu_{23}$ & 0.4 \\
		& \\
		Tensile strength in fiber direction: $S^{11}$ (MPa) & 1515 \\
		Damage-initiating tensile strain/separation in fiber direction: $\epsilon_{11}^{0}=S_{11}/E_{11}$ & 0.0109 \\
		Ultimate failure strain/separation in fiber direction: $\epsilon_{11}^u$ & 0.013 \\
		\hline
	\end{tabular}
\end{table}

As alluded to earlier, compressive failure may occur in certain situations (mostly in compression tests), which is characterized by the formation of inclined 'kink bands'. However, this phenomenon (and associated failure mode) is not treated in the present version of AGDM, but it is also not a feature of the failure process in the tensile-loaded specimens \cite{Johnson2001} to which we compare here. Thus Table \ref{tab:t300/976_yarn} does not list certain quantities that would be necessary for compressive load modeling, however, such features can be easily included in AGDM by enabling corresponding stiffness degradation and damage development in continuum yarn elements by user defined materials (VUMAT, UMAT) or elements (UEL).

\subsubsection{Matrix cracklets and delamination cracklets}
Both matrix cracklets and delamination cracklets are modeled with either eight-node or six-node cohesive elements (COH3D8 and COH3D6), and incorporate cohesive zone models have uncoupled traction-separation response:
\begin{equation}
t=\begin{bmatrix}
\tau_n \\
\tau_{s} \\
\tau_{t}
\end{bmatrix} = \begin{bmatrix}
K_n & 0 & 0 \\
0 & K_s & 0\\
0 & 0 & K_t
\end{bmatrix} \begin{bmatrix}
\delta_n \\
\delta_s \\
\delta_t
\end{bmatrix}
\end{equation}
where $\tau_n, \tau_s, \tau_t$ are the traction in normal and two orthogonal shear directions, respectively; $K_n, K_s, K_t$ are penalty stiffness in the corresponding directions, and $\delta_n, \delta_{s}, \delta_{t}$ are the relative displacement (opening) of the interface cohesive elements in the normal and shear directions, respectively. 

Damage initiation is triggered following a quadratic criterion, where the compressive traction (negative in value) along the normal direction is naturally suppressed through the notation $\langle \, \rangle$:
\begin{equation}
\left(\frac{\langle \tau_n \rangle}{T_n} \right)^2 + \left( \frac{\tau_s}{T_s} \right)^2 + \left( \frac{\tau_t}{T_t} \right)^2 = 1  
\end{equation}
where $T_n$, $T_s$ and $T_t$ are the cohesive strengths in the normal and two shear directions, respectively. Further interfacial damage evolution then follows the B-K (Benzeggagh and Kenane) law  \cite{Benzeggagh1996} based on energy dissipation, and is applied for matrix cracklets and delamination cracklets:
\begin{equation}\label{b-k law}
G_C = G_{IC} + (G_{IIC}-G_{IC}) B^{\eta}
\end{equation}
where $G_C$ is the mixed-mode fracture energy, $G_{IC}$ and $G_{IIC}$ are the Mode-I and Mode-II critical strain energy release rates, $B$ is a measure of the mode mixity and $\eta$ is a material parameter. This law is useful when the critical fracture energies are the same in the two shear directions, as is the case for the material used here. An effective variable to assess the overall damage extent, $D^e$, is computed based on the damage evolution law and mode mixity, then the stiffness of the interface is reduced accordingly in a similar manner to the damage variable used for fiber fracture, $D^f$ in Eq. \ref{eq:damage_variable}. Detailed descriptions of the aforementioned constitutive equations can be found in the Abaqus technical documentation \cite{DassaultSystemes2016}, and a study by Camanho and Davila \cite{Camanho2002}.

The effective damage variable, $D^e$, usually dominated by the leading fracture mode, is applied to reduce all penalty stiffnesses. For example, if an interface element is in pure Mode-II shear ($|\delta_s| >0, \delta_n=\delta_t=0$), then $D^e$ is only contributed by the current shear separation, strength and toughness in Mode-II fracture, but the interfacial penalty stiffnesses are all degraded by a factor of $(1-D^e)$. This will guarantee that all stress components vanish simultaneously as a damaged material point fails. It is also reasonable to assume that $D^e \in [0, 1]$ is a monotonically increasing variable, meaning that all dissipated energy or degraded stiffness can not be recovered. 

The above power law for failure initiation and B-K law for mix-mode damage evolution are mainly based on effective traction-separation laws and energy release rates (i.e., non-potential based), usually requiring a material parameter determined by experimental studies, such as $\eta$ in Eq. \ref{b-k law}. In another study, Nguyen and Waas \cite{Nguyen2016} discussed commonly used cohesive rules, then proposed a generalized potential-based mixed-mode formulation, wherein normal and tangential traction-separation laws under mixed-mode are scaled with two effective piecewise functions constructed from pure mode traction-separation laws. The proposed method ensures that stress components vanish simultaneously as crack propagates, while only cohesive strength and fracture toughness values are required. This kind of potential-based laws could be especially useful when material or non-physical parameters in non-potential based laws are hardly accessible or difficult to obtain by experimental studies. 

Attention should be paid to the chosen interface element size, here the assigned element size for interface cracklets is approximately $1 \text{mm}$, whereas the cohesive zone length of the T300/976 interface is about $0.23 \text{mm}$ based on formulas given in Turon et al. \cite{Turon2007}. Therefore, a $1 \text{mm}$ element size is considered to yield a 'coarse' mesh for cohesive zone models. 
A suggested solution \cite{Turon2007} is to reduce the interface strength as:
\begin{equation}\label{eq:reduced_interface_strength}
\bar{\tau}^0 = \sqrt{\frac{9\pi E G_c}{32 N_e l_e}} \approx \sqrt{\frac{ E G_c}{ N_e l_e}}
\end{equation}
where $N_e$ is the 'specified (desired) number of elements within the cohesive zone', $E$ is the Young's modulus (which should be interpreted as $E_{22}$ for orthotropic materials), $l_e$ is the element size, and $\bar{\tau}^0$ is the reduced interface strength.
Harper and Hallett \cite{Harper2008} investigated methods for determining appropriate values of interface strengths for delamination and discussed the limitations of existing formulas for interface strength reduction. 

Such modifications are commonly seen in literature, though a universal rule (probably impractical) has not been explicitly established: specifically for T300/976, a reduction factor of $2$ is used in \cite{Lu2018}, corresponding to a value of $N_e$ around $1.3$; Hossein et al. \cite{Hosseini-Toudeshky2013} used $N_e=5$ for delamination growth modeling of laminates according to an earlier suggestion by Mi et al. \cite{Mi1998} and Falk et al. \cite{Falk2001}. In case studies of delamination modeling by Tu et al. \cite{Turon2007}, $N_e=5$ was found to yield relatively accurate results comparing to experiments, while $N_e=10$ gave slightly more accurate results, the additional improvement was insufficient to justify the increased computational cost. Based on the above investigations, $N_e=5$ is adopted here in modeling the T300/976 interfaces, but $N_e=1$ is also applied in some simulations for comparison. The chosen parameter values for the cohesive elements are shown in Table \ref{tab:t300/976_interface}, in which the reduced strengths are computed based on Eq. \ref{eq:reduced_interface_strength} with $N_e=5$ and $l_e=1 \text{mm}$, and where no adjustment is needed for the fracture toughness values. The above CZM-based properties are directly available from a built-in material library in Abaqus, and thus no user-defined materials are needed for matrix cracklets and delamination cracklets.

\begin{table}[H]
	\centering
	\caption{Interfacial properties of T300/976 for cohesive elements in matrix cracklets and delamination cracklets \cite{Wang2009, Hosseini-Toudeshky2013, Lu2018, Torayca, Johnson2001a}. } \label{tab:t300/976_interface}
	\begin{tabular}{l l}
		\hline 
		Property & Value \\[1ex]
		\hline \hline
		Interfacial normal strength: $T_n^{*}$ & 44.5 MPa \\ 
		Interfacial shear strength: $T_s^{*}$ & 106.9 MPa \\
		Interfacial Mode-I fracture toughness: $G_{IC}, G_{IC}^{m}, G_{IC}^{d}$ & $0.0876$ $\text{kJ/m}^2$ \\
		Interfacial Mode-II fracture toughness: $G_{IIC},  G_{IIC}^{m}, G_{IIC}^{d}$ & $0.315$ $\text{kJ/m}^2$ \\[1ex]
		Adjusted normal strength of matrix and delamination cracklets ($Ne=5$): $T_n^{m}, T_n^{d}$ & 13.0 MPa \\
		Adjusted shear strength of matrix and delamination cracklets ($Ne=5$): $T_s^{m}, T_s^{d}$ & 24.7 MPa \\[1ex]
		Adjusted normal strength of matrix and delamination cracklets ($Ne=1$): $T_n^{m}, T_n^{d}$ & 29.2 MPa \\
		Adjusted shear strength of matrix and delamination cracklets ($Ne=1$): $T_s^{m}, T_s^{d}$ & 55.3 MPa \\[1ex]
		Interfacial penalty stiffness: $K_n, K_s, K_t$ & $10^6$ $\text{N/mm}^3$ \\
		Material parameter for mixed-mode B-K law: $\eta$ & 2.68 \\
		\hline
	\end{tabular}
\end{table}

\subsection{Results and discussion}

The experimental observations from the failure process in the tested $[30/90/\text{--}30]_s$ specimen \cite{Johnson2001} being for comparison, reveal that initial damage growth consists of small edge delaminations bounded by short matrix cracks, and soon followed by abrupt catastrophic failure involving massive delamination that spreads from the initial edge delaminations. Several yarn breaks are also observed at final failure stage. The small edge delaminations have virtually no effect on the specimen stiffness, i.e., little stiffness loss was observed before sudden failure of the entire specimen.

\subsubsection{Damage profile and failure evolution}

\begin{figure}[H]
	\centering
	\includegraphics[width=0.8\textwidth]{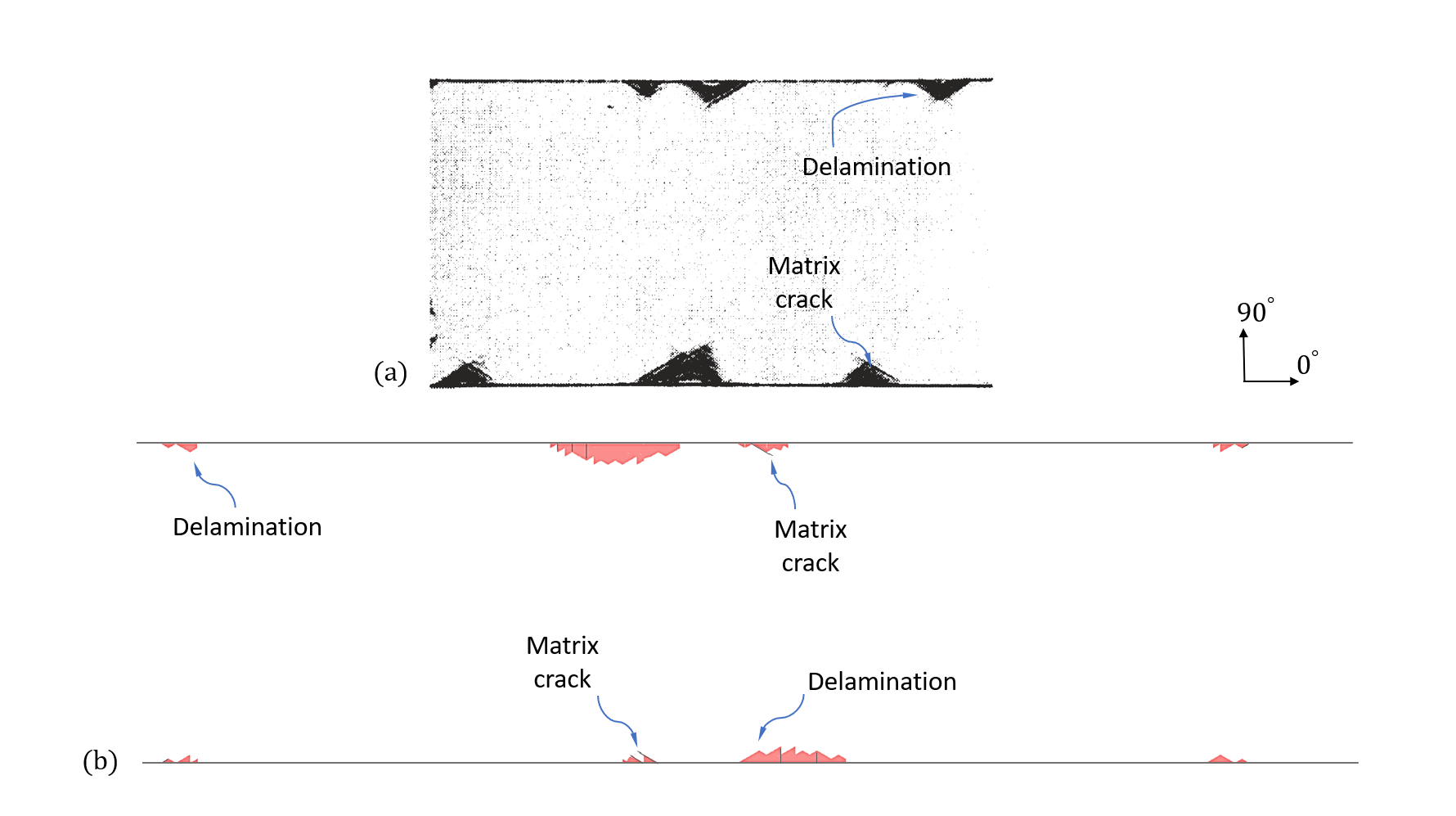}
	\caption{Damage profiles of simulated and experimental specimens ($[30/90/\text{--}30]_s$) at peak load right before catastrophic failure. (a) X-radiograph of tested specimen in experiments, just prior to the abrupt failure, by Johnson and Chang \cite{Johnson2001}; (b) Numerical simulation of damage profile by AGDM, virtually at its peak load, all failure modes are shown in reference (undeformed) laminate shape.}
	\label{fig:commparison-peak}
\end{figure}

\begin{figure}[H]
	\centering
	\includegraphics[width=0.8\textwidth]{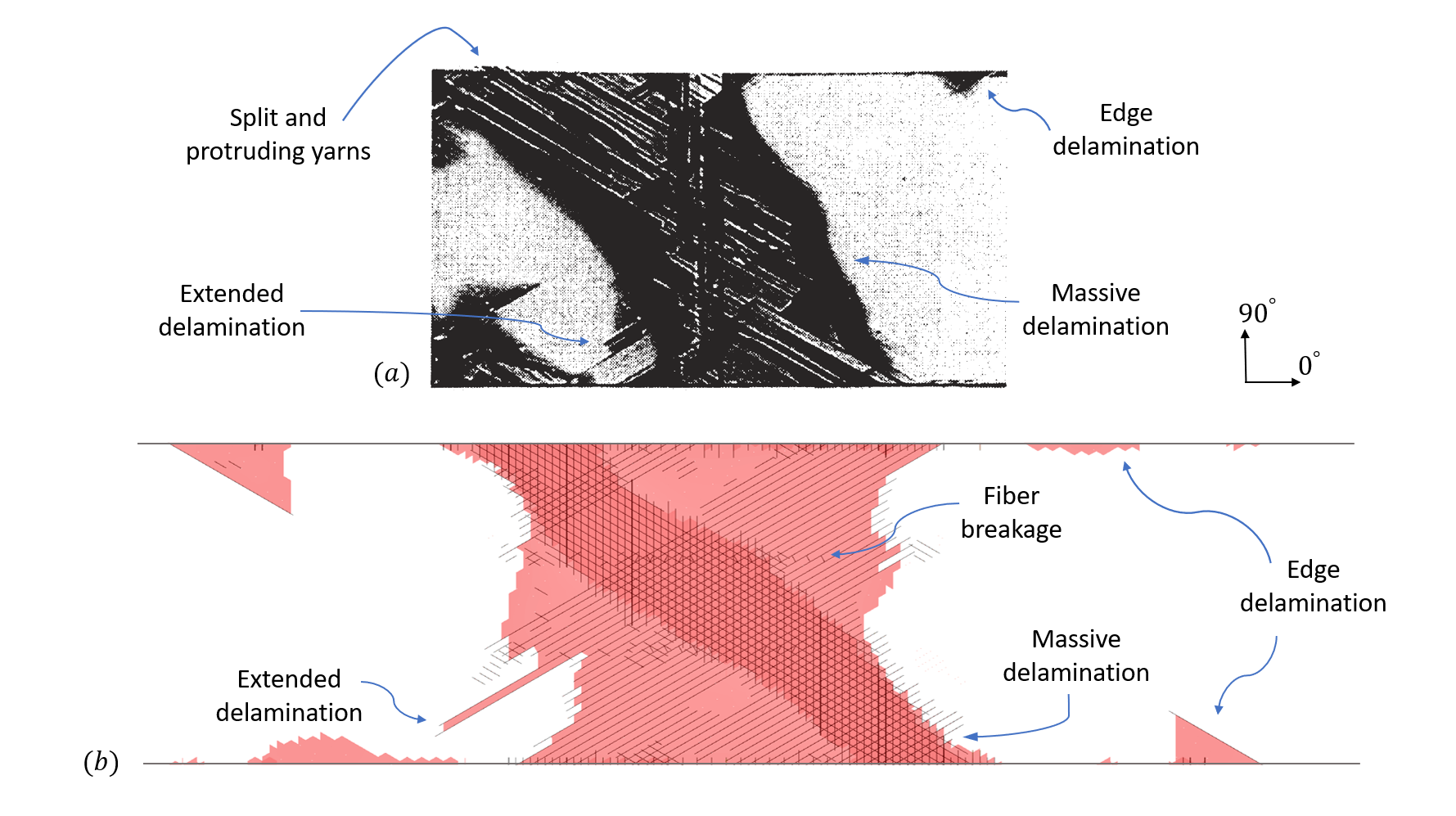}
	\caption{Damage profiles of simulated and experimental specimens ($[30/90/\text{--}30]_s$) as they fail catastrophically. (a) X-radiograph of failed experimental specimen as it passes the peak load, by Johnson and Chang \cite{Johnson2001}; (b) Damage profile seen in the AGDM numerical simulation of the specimen at peak load just as catastrophic failure occurs, and where all damages are shown in the reference (undeformed) laminate shape. }
	\label{fig:commparison-failed}
\end{figure}

For the damage profile and progressive failure evolution, the numerical prediction by AGDM is in close agreement with the experimental observations. Fig. \ref{fig:commparison-peak} shows a comparison between the X-radiograph and numerical prediction of the $[30/90/\text{--}30]_s$ specimen just before final abrupt failure. Fig. \ref{fig:commparison-failed} shows comparison between X-radiograph and numerical prediction of the abruptly-failed specimen right at peak load.

In the numerical analysis, several isolated small edge delaminations are observed midway along the specimen, which occur just prior to the abrupt failure. Tiny edge delaminations are also initiated close to the ends of the specimen but are not a significant factor in the final failure. Several matrix cracks associated with edge delaminations are also observed. Almost instantly, upon having just reached the peak load, massive delamination occurred in all plies and formed a wide delamination band that initiated from previously small edge delaminations. Intensive matrix cracks occurred inside and around the border of the massive delamination region. Around the edge of massive delamination region, several split yarns are formed by bounded matrix cracks and local delaminations, and have noticeably protruded out of the laminate edge. Several extended delamination blocks are also observed, which are bounded by matrix cracks that propagated further than other neighboring cracks. Smaller isolated edge delaminations away from the massive delamination band occurred as well, but do not dominate the abrupt failure of the laminate. 
Finally, several yarn fractures (fiber breaks) have taken place near the center of the delamination band in both the $30^{\circ}$ and $-30^{\circ}$ plies.

\begin{figure}[H]
\centering
	\begin{subfigure}[t]{0.8\textwidth}
		\includegraphics[width=\textwidth]{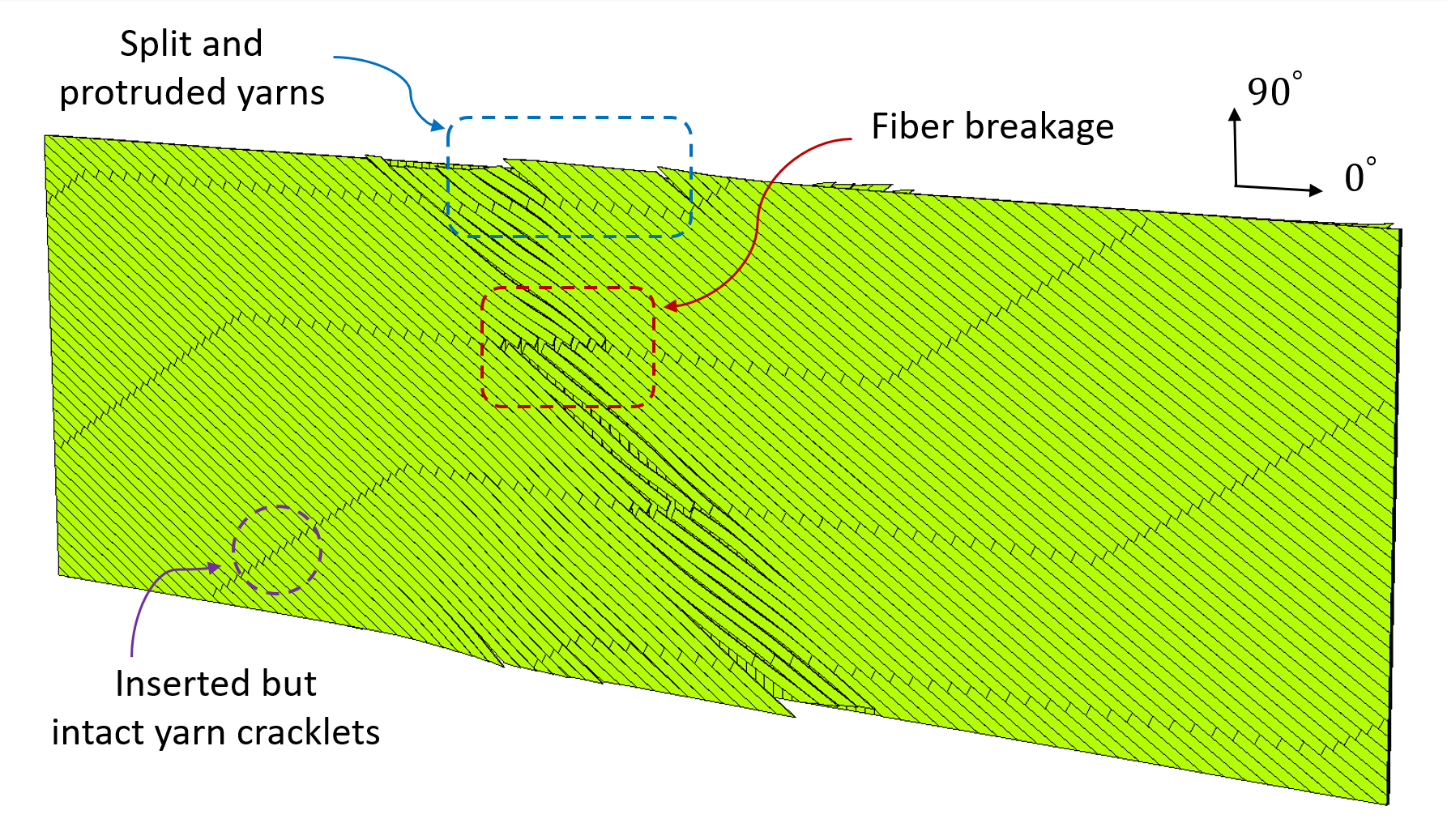}
		\caption{}
		\label{fig:yarns-failed}
	\end{subfigure}
	\begin{subfigure}[b]{0.5\textwidth}
		\includegraphics[width=\textwidth]{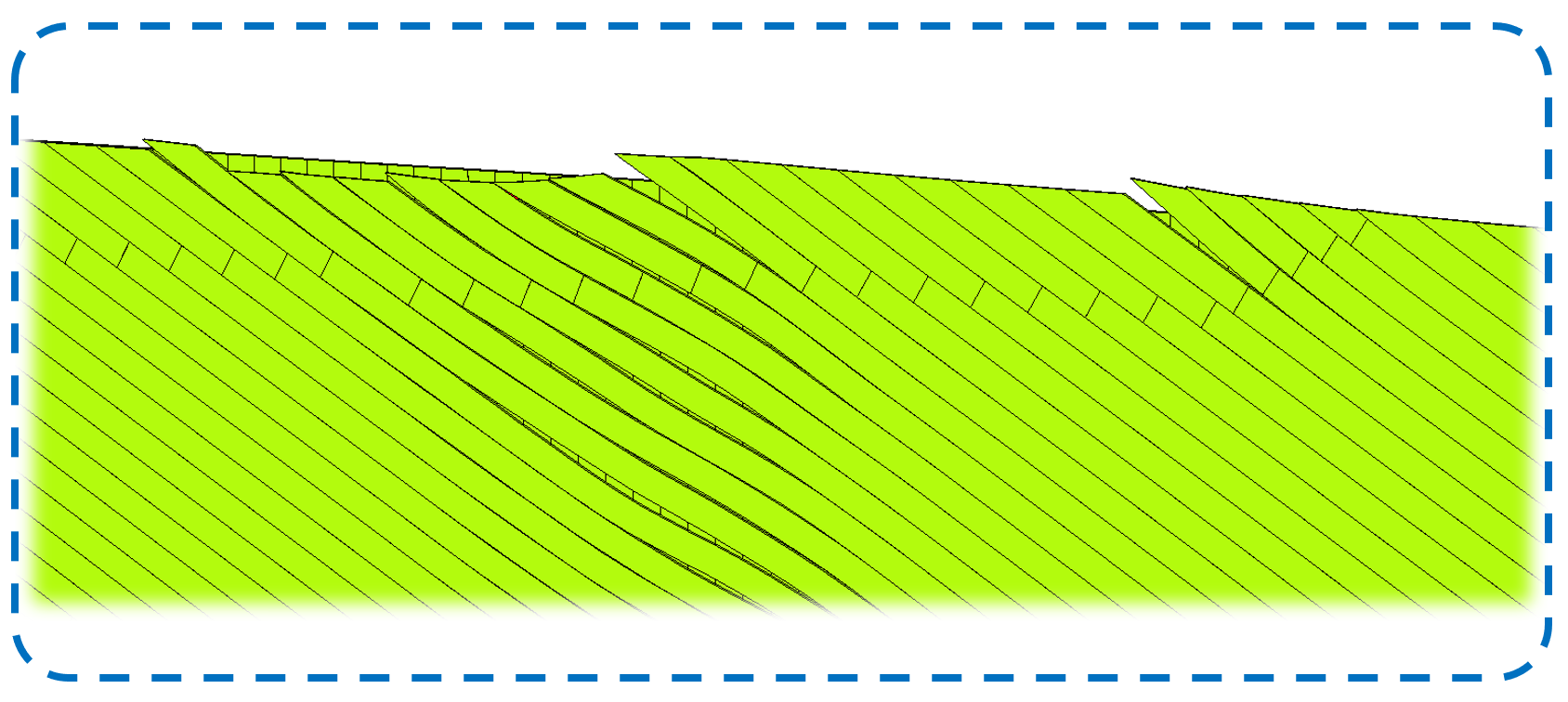}
		\caption{}
		\label{fig:protruded-yarns}
	\end{subfigure}
	\begin{subfigure}[b]{0.4\textwidth}
		\includegraphics[width=\textwidth]{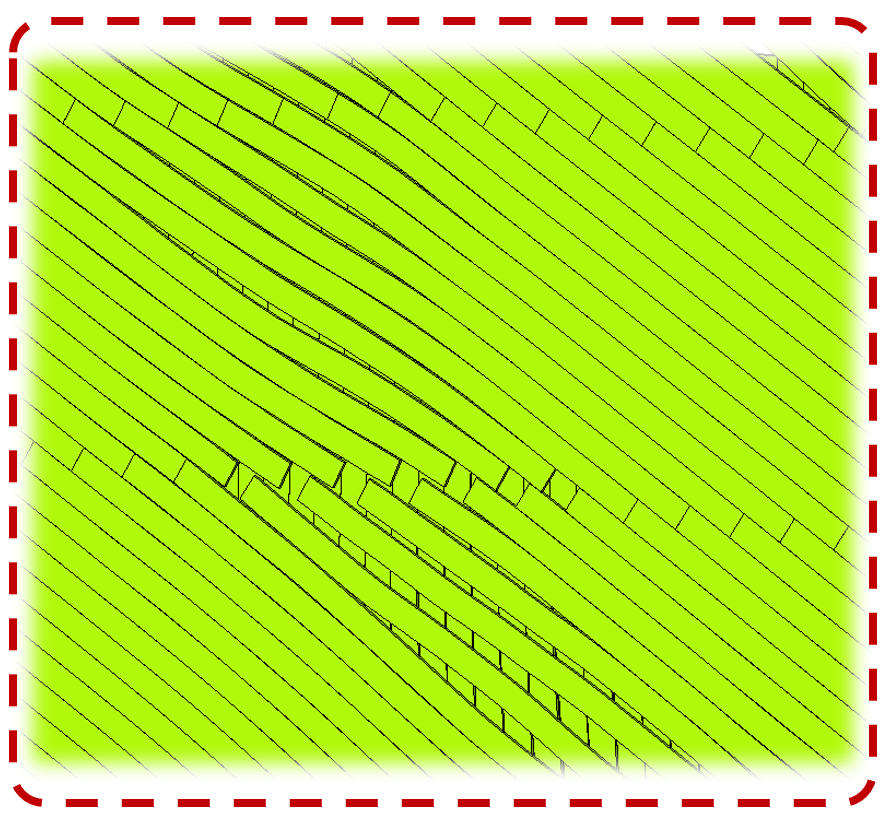}
		\caption{}
		\label{fig:yarn-breaks}
	\end{subfigure}
\caption{Damage pattern of $[30/90/\text{--}30]_s$ specimen in true deformed shape. (a) Failed laminate at peak load, showing yarns separated by matrix cracks and protruding from the sides of massive delamination band. Yarn breaks have also occurred, which agrees with observations in experiments \cite{Johnson2001}. (b) Close-up view of split and protruded yarns at one side of a delamination band in Fig. \ref{fig:yarns-failed}. (c) Close-up view of yarn breaks near the middle of massive delamination band in Fig. \ref{fig:yarns-failed}. }
\label{fig:yarn-damages-global-and-close}
\end{figure}

Fig. \ref{fig:yarn-damages-global-and-close} shows the failed specimen (corresponding to the same output frame shown Fig. \ref{fig:commparison-failed}), but in its true deformed shape, where failed yarns, and patterns of matrix cracks and delaminations at interfaces can be clearly seen. This figure also includes close-up views of the split and protruded yarns at the edges of a massive delamination band, as well as yarn fractures in the middle of the delamination band. 

It should be emphasized that the visualized damage patterns are based on user-defined graphical outputs, and at times that are specified in advance. For each numerical analysis of an $[30/90/\text{--}30]_s$ specimen, a total of $400$ outputs are saved at evenly spaced time marks over the total test time ($4$ seconds). Thus with a targeted time increment of $10^{-6}$ seconds, around $4$ million increments are required, it is highly impractical to output all of them. The modeled specimen completely failed in two consecutive outputs, and therefore the damage pattern in Fig. \ref{fig:commparison-peak}(b) is based on the output at peak load, while as shown in the next available output in Fig. \ref{fig:commparison-failed}(b), massive delamination and abrupt failure had already occurred. 
If more intervals were saved as outputs near the point of peak load, a damage pattern somewhere between Fig. \ref{fig:commparison-peak}(b) and Fig. \ref{fig:commparison-failed}(b) would have been available to better match the X-radiograph of damage in the experimental specimen. Nevertheless, close examination and comparison of the damage patterns in two consecutive outputs indicates that all important features are captured well in the simulation, despite the limitations in available output frames for comparison.     

In brief, the features of the above damage profiles observed in the numerical simulation are in excellent qualitative agreement with the experimental observations. As evidenced in both Fig. \ref{fig:commparison-failed} and Fig. \ref{fig:yarn-damages-global-and-close}, AGDM not only successfully predicted the shape and pattern of the massive delamination bands, but also captured other damage modes such as split and protruded yarns at the specimen edges, extended delaminations, isolated edge delaminations and yarn breakage. Qualitatively, almost all the essential damage modes and their corresponding shapes, patterns and locations are well predicted by AGDM. Also note that predicted failure pattern is generally centro-symmetric in the numerical simulation, while in reality, the progressive damage pattern is expected to be more complicated due to multiple random factors, such as material imperfections and the statistical strength distribution of carbon fibers and yarns in composites, and the loading profile (loading rate, and whether a fatigue or static test). These aspects are worthy of study but are beyond the scope of present work.

\subsubsection{Reaction load vs applied displacement}\label{subsec:load-vs-disp}

\begin{figure}[H]
	\centering
	\includegraphics[width=0.8\textwidth]{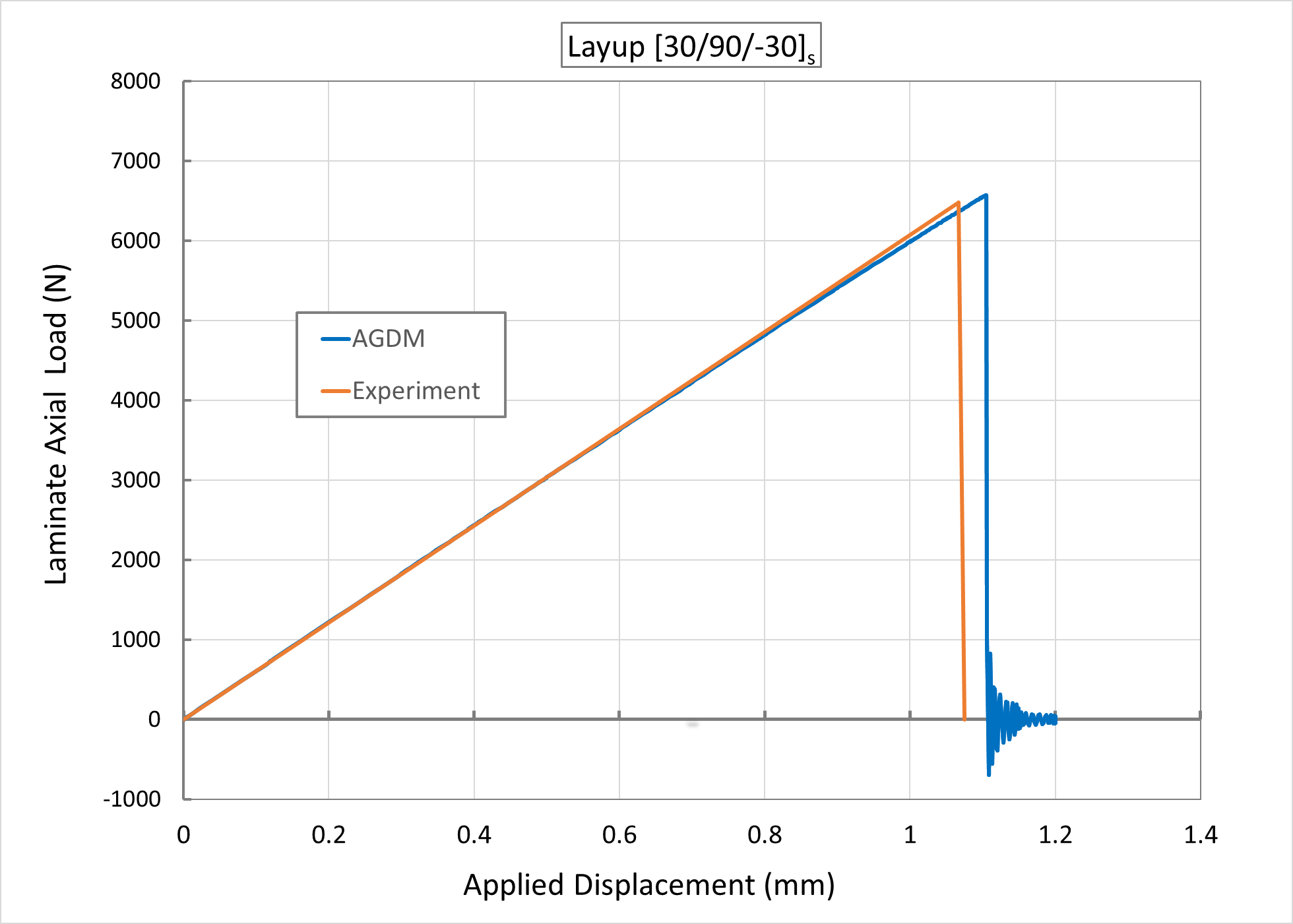}
	\caption{Reaction load versus applied displacement: a comparison between numerical simulation and data from tested $[30/90/\text{--}30]_s$ specimen.}
	\label{fig:load-vs-disp-1}
\end{figure}

Fig. \ref{fig:load-vs-disp-1} shows and compares the reaction force versus applied displacement in numerical simulation to the experimental data of the $[30/90/\text{--}30]_s$ specimen. The simulated laminate also suffers catastrophic failure, where the reaction load drops rapidly from the peak load to almost zero in less than two numerical outputs, where one output approximately refers to an end-to-end displacement increment of $0.003 \text{mm}$. The peak reaction load and applied displacement at peak load in the experimental specimen are $6483 \text{N}$ and $1.067 \text{mm}$, respectively, based on the data provided in \cite{Johnson2001}. In numerical simulation, the specimen failed at $6574\text{N}$ with corresponding displacement of $1.104 \text{mm}$, and prior to that, the stiffness was almost constant until sudden failure, which agrees well with the experimental observations.

The average stiffness based on the final peak load and displacement is only $2.0\%$ less than the stiffness computed from experimental data, this value agrees with the analysis by Turon et al. \cite{Turon2007}: Simply speaking, such stiffness loss is inevitable when cohesive surfaces or small/zero thickness interface elements are applied between plies. However, as long as the penalty stiffness of a cohesive interface is properly chosen (large enough), stiffness loss of the laminates due to the presence of cohesive interfaces will be less than $2.0$\%, which is sufficiently accurate for most applications. Note that an extremely large penalty stiffness value will cause less reduction of the simulated structures in component scale, but could introduce numerical problems such as spurious oscillation of interfacial tractions \cite{Turon2007}. 

Affected by the loading rate and artificial mass scaling applied in the numerical simulation, a moderate and decaying dynamic oscillatory response occurred following the sudden load-drop after achieving the peak load. The introduced kinematic energy (dynamic oscillation) were then gradually decayed with a rate controlled by the default small damping factors \cite{DassaultSystemes2016}. Given the fact that the specimen failed abruptly, this kind of dynamic effect is virtually inevitable in a Dynamic/Explicit simulation, however, in this specific study, such a post failure reaction is not a feature with physical significance. As previously discussed in Section \ref{geometry-and-bc}, the mass scaling and loading profiles are carefully chosen to ensure that the quasi-static features are well persevered in the time and load region of interest.

\subsubsection{Mesh convergence study}

\begin{figure}[H]
	\centering
	\includegraphics[width=0.8\textwidth]{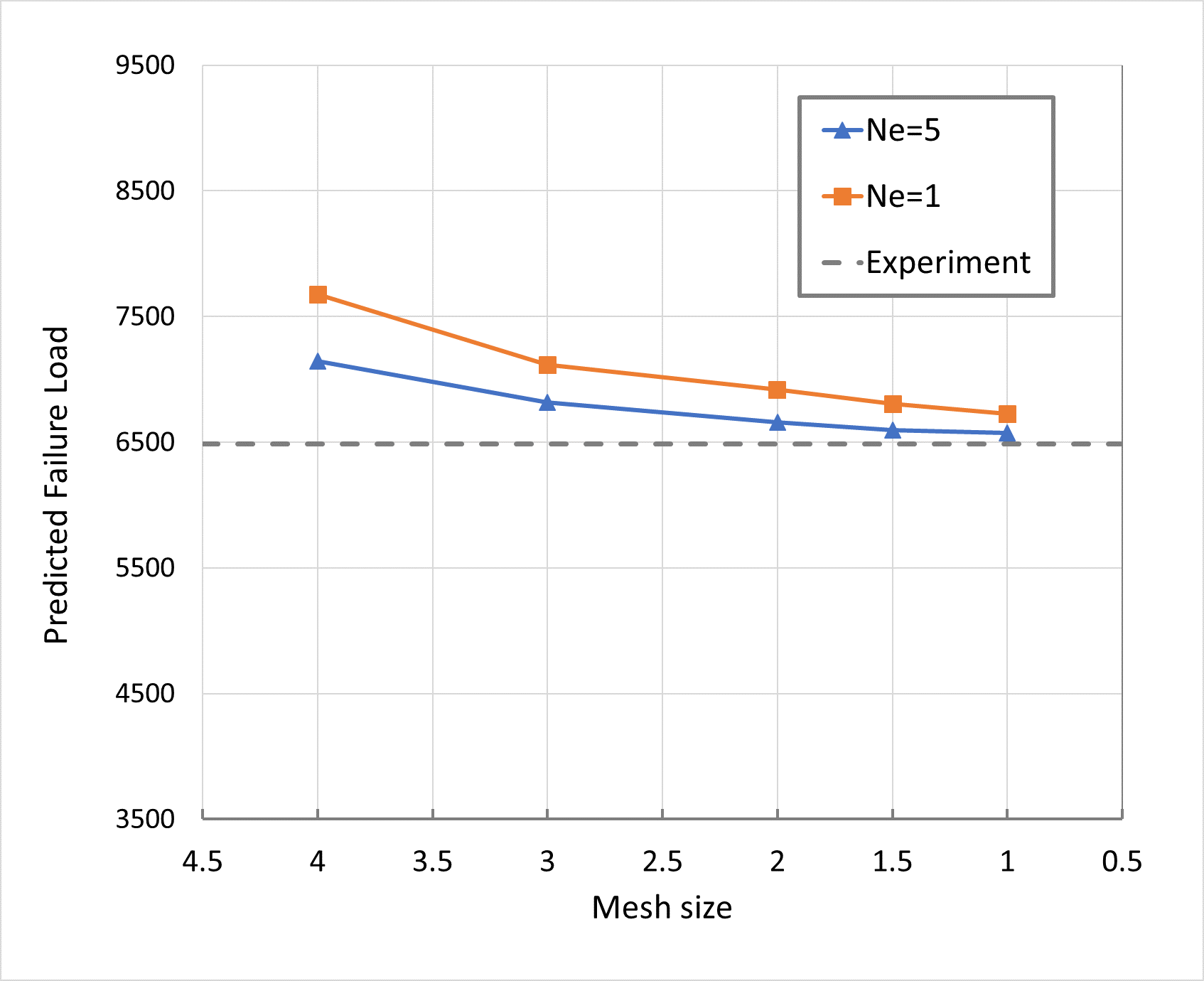}
	\caption{Mesh convergence study on $[30/90/\text{--}30]_s$ specimens with various mesh sizes and two choices of $N_e$.}
	\label{fig:mesh-study}
\end{figure}

A brief mesh convergence study on the $[30/90/\text{--}30]_s$ laminate specimen, involving different choices of mesh sizes and two choices of the 'specified number of elements in the cohesive zone', $N_e$ in Eq. \ref{eq:reduced_interface_strength}, is shown in Fig. \ref{fig:mesh-study}. The predicted failure load stabilizes and shows sign of convergence as element size becomes smaller than $2 \text{mm}$. As expected, for the same mesh size, the predicted ultimate loads are lower when $N_e=5$ as compared to the case when $N_e=1$ due to the further reduced interfacial strengths. The differences in predicted failure load due to the choice of $N_e$ are less significant as element size becomes smaller.

In the cases shown in Fig. \ref{fig:mesh-study}, mesh element sizes are the same between continuum yarns and cohesive interfaces, which is not necessarily the case in AGDM. For example, if extra fine mesh and output results involving interfaces are required, then the yarn width (matrix crack spacing) may remain as $1 \text{mm}$, while the interface (cracklets) elements can be reduced to sizes smaller than the cohesive zone length of the material, such as $0.125 \text{mm}$. 

Several simulations using different mesh sizes between yarns and interfaces (cracklets) are also investigated, the results together with some results in Fig. \ref{fig:mesh-study} are listed in Table \ref{tab:result_table_1}. It is observed that the mesh size of interface cohesive elements, together with the choice of $N_e$, have more noticeable effect on the predicted value of ultimate strength than other parameters, this phenomena is expected since all discrete failures are accounted for by the interface elements (inserted cracklets). The predicted progressive damage pattern is less sensitive to mesh sizes and becomes consistent if all mesh sizes are smaller than $2 \text{mm}$. The predicted values of ultimate strength tend to stabilize as element sizes are approaching $1 \text{mm}$. 

\begin{table}[H]
	\centering
	\caption{Mesh convergence studies of $[30/90/\text{--}30]_s$ laminate to compare with experiments by Johnson and Chang \cite{Johnson2001}.} \label{tab:result_table_1}
	\begin{tabularx}{\textwidth}{C{0.5} C{1.25} C{1.0} C{1.25} C{1.25} C{1.25} C{1.0} C{0.5}}
		\hline 
		Case & Laminate layup & Continuum yarn mesh size (mm) & Cohesive interface mesh size (mm) & Specified number of elements in cohesive zone, $N_e$ & Average strength of tested specimens (MPa) & Strength predicted by AGDM (MPa) & Error ($\%$) \\[0.5ex]
		\hline \hline
		1.1 & $[30/90/\text{--}30]_s$ & $2.0$  & $2.0$  & 1 & $377.6$ & $402.9$ & $+6.7$ \\
		1.2 & $[30/90/\text{--}30]_s$ & $2.0$  & $2.0$  & 5 & $377.6$ & $387.8$ & $+2.7$ \\
		1.3 & $[30/90/\text{--}30]_s$ & $1.5$  & $1.5$  & 1 & $377.6$ & $396.2$ & $+4.9$ \\
		1.4 & $[30/90/\text{--}30]_s$ & $1.5$  & $1.5$  & 5 & $377.6$ & $385.3$ & $+2.0$ \\
		1.5 & $[30/90/\text{--}30]_s$ & $1.0$  & $1.0$  & 1 & $377.6$ & $391.8$ & $+3.8$ \\
		1.6 & $[30/90/\text{--}30]_s$ & $1.0$  & $1.0$  & 5 & $377.6$ & $382.9$ & $+1.4$ \\
		1.7 & $[30/90/\text{--}30]_s$ & $1.5$  & $1.0$  & 1 & $377.6$ & $392.8$ & $+4.0$ \\
		1.8 & $[30/90/\text{--}30]_s$ & $1.5$  & $1.0$  & 5 & $377.6$ & $384.1$ & $+1.7$ \\
		1.9 & $[30/90/\text{--}30]_s$ & $2.0$  & $1.0$  & 1 & $377.6$ & $393.1$ & $+4.1$ \\
		1.10 & $[30/90/\text{--}30]_s$ & $2.0$  & $1.0$  & 5 & $377.6$ & $384.4$ & $+1.8$ \\
		\hline
	\end{tabularx}
\end{table}

\subsubsection{Necessity of having matching meshes and multiple types of failure modes}

The same $[30/90/\text{--}30]_s$ specimen was also modeled by Lu et al. \cite{Lu2018}, in which the numerical model with standard cohesive elements (non-matching mesh) did not the capture correct damage evolution and significantly over-predicted the failure load with an error of $+51\%$, while the model with proposed separable cohesive elements (SCE) by them \cite{Lu2018} realized the matching mesh between inter-ply delaminations and intra-ply matrix cracks, and improved the prediction to $+16\%$ with a mesh size of $1 \text{mm}$ and $N_e$ of about $1.3$. This comparison indicates the necessity and importance of having matching meshes between cracks.

However, since only one type of intra-ply failure is supported in SCE, potential fiber fractures in the $[30/90/\text{--}30]_s$ specimen are absent. Damage in laminates are only represented by matrix cracks and delaminations, and therefore cannot degrade as rapidly, which resulted in a more gradual damage propagation and over-prediction of the specimen failure load, as discussed in Section \ref{multi-discrete}. For $0^{\circ}$ plies, the one type of intra-ply failure is commonly assigned to fiber fractures, as is also the case in \cite{Lu2018}, but this in turn leads to the absence of potential matrix cracks in those plies. 

In this paper, all types of discrete failures are supported in the $[30/90/\text{--}30]_s$ model generated by AGDM. 
The mesh sizes and material properties in Case 1.5 is comparable to the inputs used by Lu et al. \cite{Lu2018} for the same laminate, and the prediction error here is reduced from $+16\%$ to $+3.8\%$. It is believed that a major reason for such improvement is the added modeling capability of yarn fractures (triggered by clusters of tensile fiber breaks) in AGDM: since several fiber/yarn breakage in the $30^{\circ}$ and $\text{--}30^{\circ}$ plies did occur during abrupt failure of the specimen (as described by Johnson and Chang in the experimental study \cite{Johnson2001}), and were successfully captured by failure of several yarn cracklets near the center of the massive delamination band (Fig. \ref{fig:commparison-failed} and \ref{fig:yarn-damages-global-and-close}). In Case 1.6, by using $N_e=5$ as suggested in \cite{Hosseini-Toudeshky2013,Mi1998,Falk2001}, the error of numerical prediction is further reduced from $+3.8\%$ to $+1.4\%$. 

\subsubsection{Numerical case studies on other specimens with different layups}

Eight other layups in the same series of experimental tests by Johnson and Chang \cite{Johnson2001} have been simulated by AGDM and the results are summarized in Table \ref{tab:result_table_2}. The mesh sizes and material parameters are determined based on the simulation and mesh convergence studies of the $[30/90/\text{--}30]_s$ laminate as shown in Fig. \ref{fig:mesh-study} and Table \ref{tab:result_table_1}. More specifically, since the inputs used in Case 1.6 resulted in the closet prediction comparing to reported strengths from experiments, they are adopted here for numerical studies of other laminates, while the mechanical properties for all materials remain the same as listed in Table \ref{tab:t300/976_yarn} and \ref{tab:t300/976_interface}. Each model (with $1 \text{mm}$ mesh size) requires around $20$ minutes to $3$ hours for the preprocessing commands to generate, and then takes about 48 to 70 hours for the Abaqus/Explicit solver to run on a 32-core desktop workstation. Among which, the $[50_2/75_2/ \text{--}75_2/\text{--}50_2]_s$ laminate (Case 9) took the longest time to complete: about 3 hours to create the entire finite element model, then 70 hours and 11 minutes to compute the numerical simulation executed in parallel mode with 24 CPU cores (48 threads); the $[60_4/\text{--}60_4]_s$ model (Case 5) took the least amount of time: about $25$ minutes to create, then 48 hours to complete in parallel mode with 16 CPU cores (32 threads).  

\begin{table}[H]
	\centering
	\caption{Case studies of selected specimens with various layups to compare with experiments \cite{Johnson2001}} \label{tab:result_table_2}
	\begin{tabularx}{\textwidth}{C{0.5} C{2.0} C{1.0} C{1.0} C{1.0} C{1.0} C{1.0} C{0.5}}
	\hline 
	Case & Laminate layup & Continuum yarn mesh size (mm) & Cohesive interface mesh size (mm) & Specified number of elements in cohesive zone, $N_e$ & Average strength of tested specimens (MPa) & Strength predicted by AGDM (MPa) & Error ($\%$) \\[0.5ex]
	\hline \hline
	2 & $[45_2/ \text{--}45_2]_s$ & $1.0$  & $1.0$  & 5 & $141.9$ & $147.0$ & $+3.6$ \\
	3 & $[45_3/ \text{--}45_3]_s$ & $1.0$  & $1.0$  & 5 & $116.4$ & $120.2$ & $+3.3$ \\
	4 & $[60_2/ \text{--}60_2]_s$ & $1.0$  & $1.0$  & 5 & $69.6$ & $65.5 $ & $-5.9$ \\
	5 & $[60_4/ \text{--}60_4]_s$ & $1.0$  & $1.0$  & 5 & $48.2$ & $45.1$ & $-6.4$ \\
	6 & $[0/45/ \text{--}45/90]_s$ & $1.0$  & $1.0$  & 5 & $487.1$ & $499.4$ & $+2.5$ \\
	7 & $[90_2/ 45_2/ \text{--}45_2/ 0_2]_s$ & $1.0$  & $1.0$  & 5 & $510.5$ & $534.6$ & $+4.7$ \\
	8 & $[55_2/ 75_2/ \text{--}75_2/ \text{--}55_2]_s$ & $1.0$  & $1.0$  & 5 & $75.1$ & $68.3$ & $-9.1$ \\
	9 & $[50_2/ 75_2/ \text{--}75_2/ \text{--}50_2]_s$ & $1.0$  & $1.0$  & 5 & $101.3$ & $92.9$ & $-8.2$ \\		
	\hline
	\end{tabularx}
\end{table}

The predicted values of specimen strength are generally in excellent agreement with the experiments, with the largest absolute value of error being $9.1\%$ (Case 8). It should be emphasized here that this kind of simulations usually involve a large number of elements, interactions and inputs, thus the numerical outputs are potentially affected by various factors like mesh sizes, selected value of interfacial cohesive strength and fracture toughness, damage evolution type and mode mixity parameters, etc. For example, if $N_e$ is determined as $1$ instead of $5$ for the cases in Table \ref{tab:result_table_2}, the cohesive zone strength will increase (Eq. \ref{eq:reduced_interface_strength}), then the prediction accuracy of laminates with larger layup angles (Case 4,5,8 and 9) should be improved, while in turn the error of several laminates (Case 2,3,6 and 7) may become larger. In certain cases, inputs/factors that would over-predict the results might cancel out with ones that would under-predict the strength. Therefore, the exact value of prediction error may not be a major concern here. 

\begin{figure}[H]
	\centering
	\begin{subfigure}[t]{0.6\textwidth}
		\includegraphics[width=\textwidth]{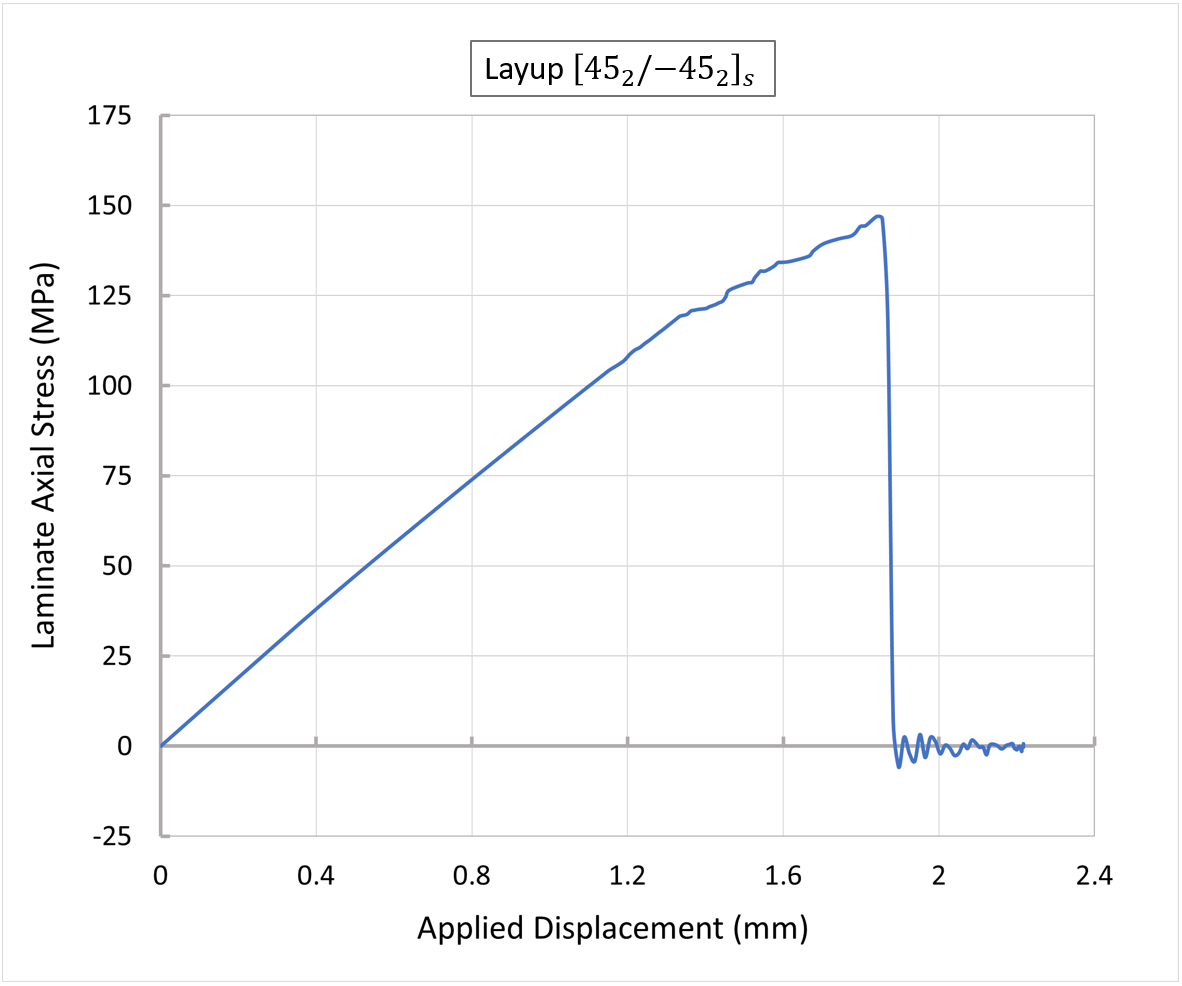}
		\caption{Reaction stress versus applied displacement load in numerical simulation (uniaxial tensile test).}
		\label{fig:load-vs-disp-2}
	\end{subfigure}
	\begin{subfigure}[b]{0.6\textwidth}
		\includegraphics[width=\textwidth]{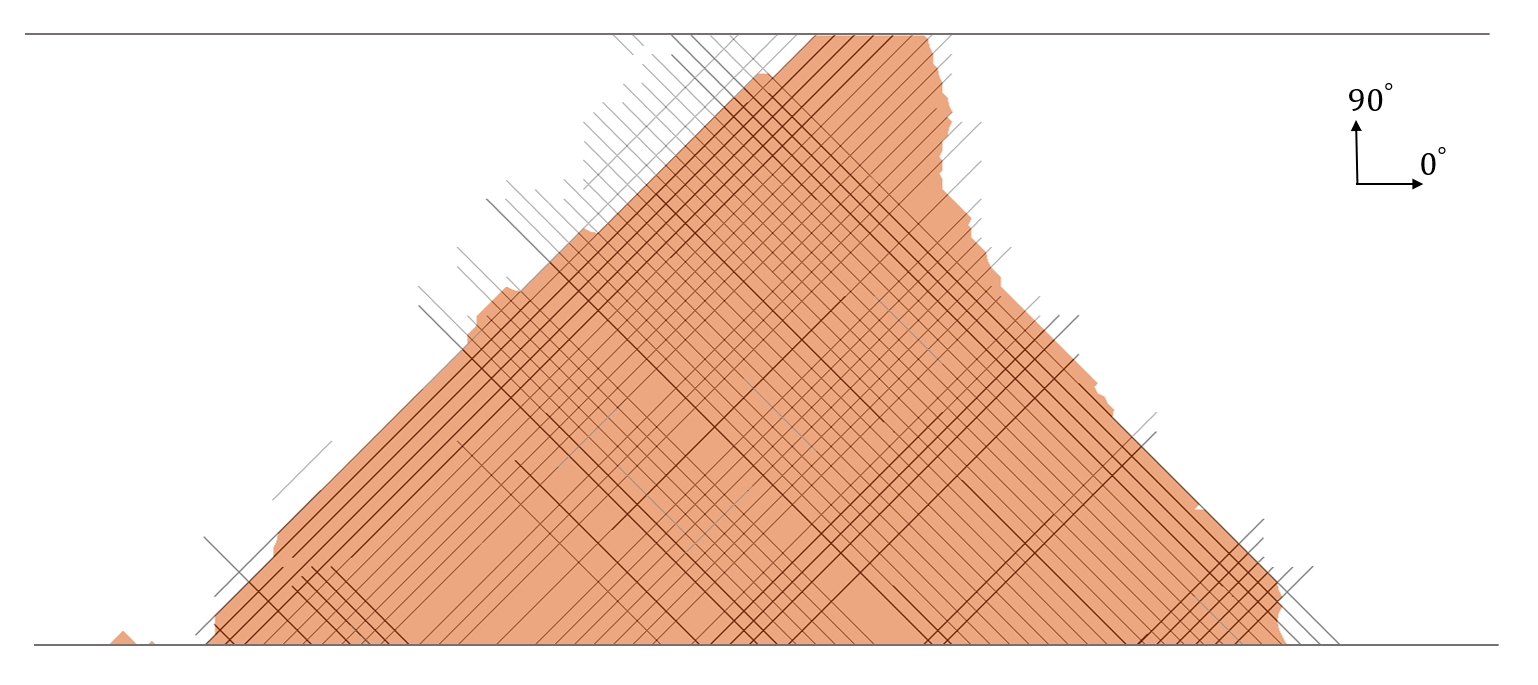}
		\caption{Predicted damage pattern of the specimen after ultimate failure.}
		\label{fig:damage-pattern-45_-45}
	\end{subfigure}
	\caption{Numerical prediction of laminate axial stress and damage pattern of the $[45_{2}/\text{--}45_{2}]_s$ specimen (Case 2).}
\end{figure}

It is also noted that $[45_{i}/\text{--}45_{i}]_s$ specimens in the experimental studies showed non-linear force-vs-load response, which was believed to be attributed by non-linear shear response of matrix rather than damage accumulations \cite{Johnson2001}. 
In this study, non-linear material response is not included into the corresponding constitutive laws of matrix and delamination cracklets, therefore the predicted response of such laminates (Case 2 and 3) does not show any apparent sign of non-linear/yielding behavior. For example, the stress-vs-displacement curve of $[45_{2}/\text{--}45_{2}]_s$ laminate is shown in Fig. \ref{fig:load-vs-disp-2}, only very slight stiffness degradations are observed in the later stage of the tensile loading; the predicted damage pattern of the failed specimen is as shown in Fig. \ref{fig:damage-pattern-45_-45}. 

In general, non-linear response of composite laminates in tensile tests could be attributed by various factors, such as stacking sequence, shape and initial notches, constitutive laws for interfaces, accumulation of damages, loading rate and type (monotonic or cyclic). For example, in certain specimens, non-linear response may be directly resulted from the softening shear behavior of the epoxy matrix; while in other specimens, non-linear response may be mainly contributed by gradual stiffness reduction due to accumulation of matrix cracks and delaminations, or if the specimen contains damages from previous loads (in a cyclic loading test). Potential studies involving these phenomena are complicated and can be adopted in future applications of AGDM, for example, by adding non-linear shear response into matrix and delamination cracklets through user-defined material subroutines (UMAT/VUMAT).  

A comprehensive literature search on the reported material properties of T300/976 and suggested interfacial cohesive zone parameters are presented in Section \ref{subsec:materials}, which indicates that the inputs adopted in the simulations are reasonable and accurate. In fact, all essential parameters used here are generally the same as values used in other numerical studies of the same type of material (T300/976), as can be found in \cite{Wang2009, Hosseini-Toudeshky2013, Johnson2001a, Lu2018}. As discussed above, a direct comparison between AGDM and other methods on the same specimens (in Case 1) with comparable inputs indicates that the prediction is significantly improved by modeling all types of failures meanwhile ensuring explicit coupling between them.

\begin{figure}[H]
	\centering
	\includegraphics[width=0.6\textwidth]{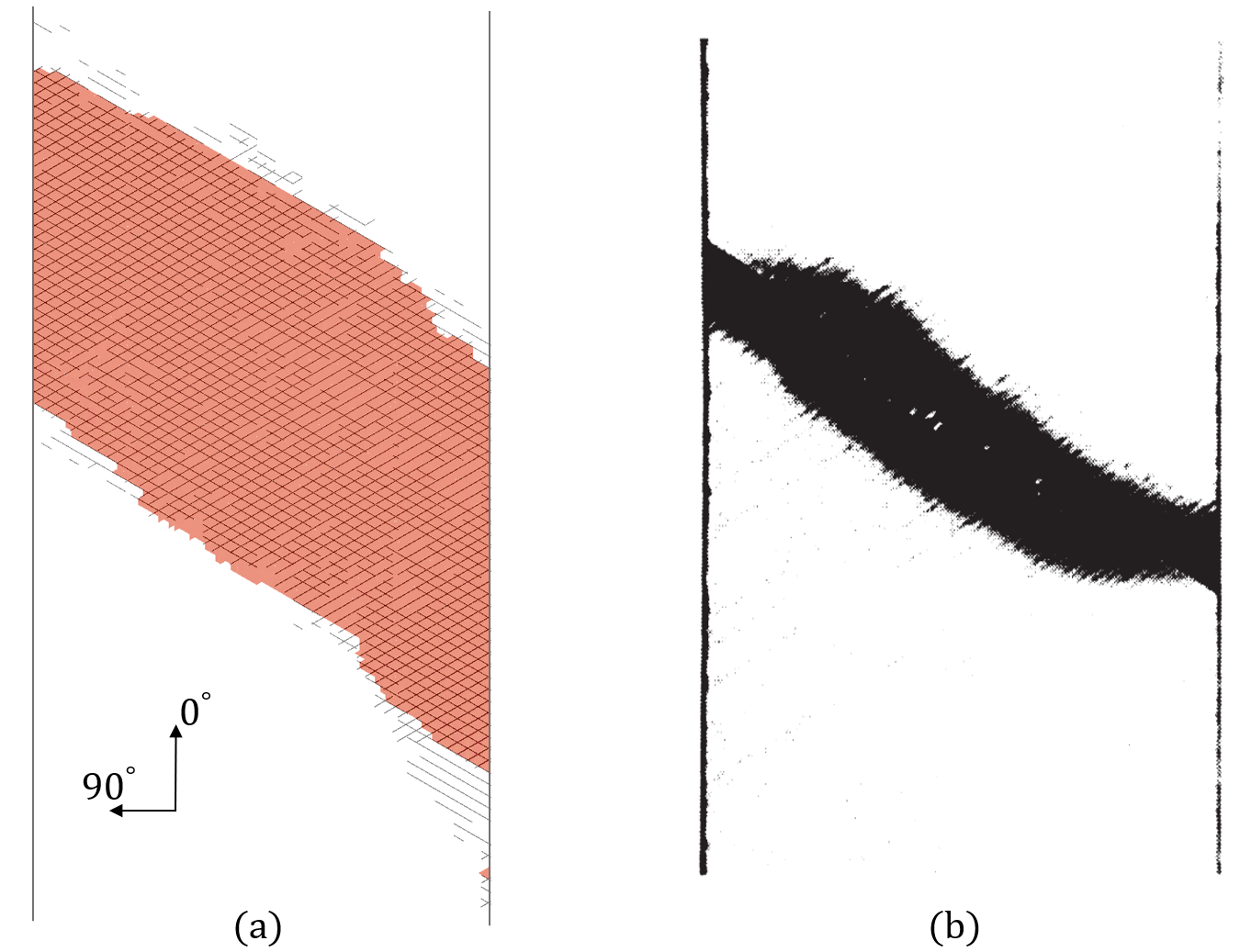}
	\caption{Damage profile of $[60_2/\text{--}60_2]_s$ laminate specimens. (a) X-radiograph of a damaged specimen, by Johnson and Chang \cite{Johnson2001}; (b) Predicted damage pattern of the specimen (Case 4) after ultimate failure.}
	\label{fig:damage_pattern_60-60}
\end{figure}

Fig. \ref{fig:damage_pattern_60-60} shows the X-radiograph of a damaged  $[60_2/ \text{--}60_2]_s$ laminate from experimental study and the predicted ultimate failure pattern of the same laminate by AGDM (Case 4). The $[60_{i}/ \text{--}60_{i}]_s$ specimens, despite having similar stacking formation and close layup angles to the $[45_{i}/ \text{--}45_{i}]_s$ laminates, do not show non-linear response in the experimental study. It also seems that the predicted damage patterns for these simple layup laminates (Case 2 to 5) have more 'widely spread' delaminations induced with matrix cracks, this is believed to be caused by the increased loading rate and 'uniform' material strength in the numerical simulations: as explained above, more objective and explicit studies involving the progressive damages should be performed with statistical material and interfacial strengths, which is part of the future work.

\begin{figure}[H]
	\centering
	\includegraphics[width=0.8\textwidth]{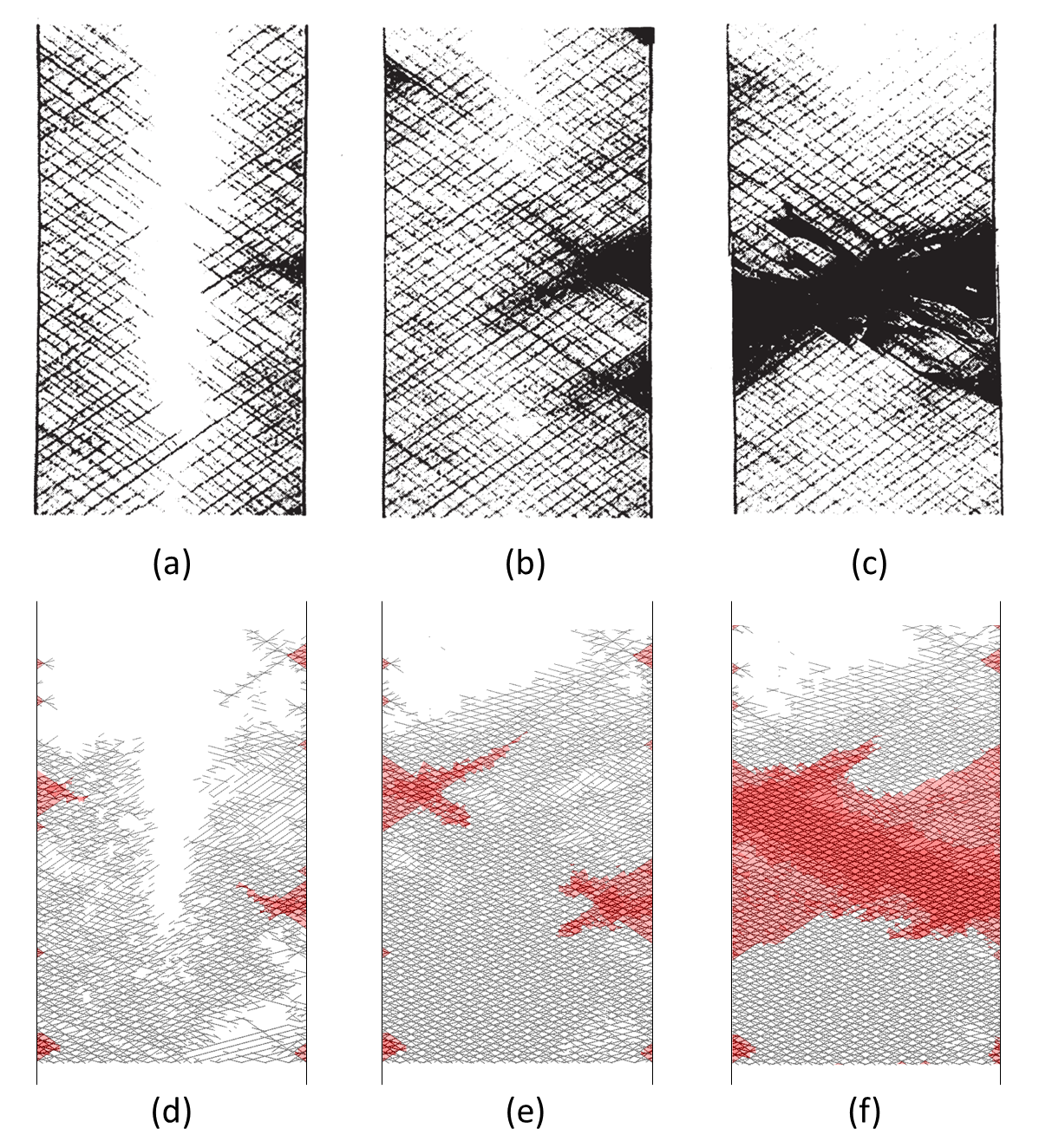}
	\caption{Progressive damage profile in experiment and its comparison with prediction by numerical simulations on $[55_2/75_2/ \text{--}75_2/\text{--}55_2]_s$ laminate. (a), (b) and (c) are X-radiograph of specimen in experimental study by Johnson and Chang \cite{Johnson2001}; (d), (e) and (f) are predicted damage patterns in numerical simulation by AGDM (Case 8) at $97.4\%$, $99.1\%$ and $100 \%$ of the ultimate axial stress (68.3 MPa), respectively.}
	\label{fig:damage_4layer}
\end{figure}

As another example, Fig. \ref{fig:damage_4layer} shows the comparison between X-radiograph and predicted damages of a $[55_2/ 75_2/ \text{--}75_2/ \text{--}55_2]_s$ laminate specimen (Case 8) by AGDM. The predicted failure profiles and progressive damage patterns agree well with the experimental observations: Matrix cracks and delaminations were largely initiated and propagated in the $55$ and $\text{--}55$ plies, then quickly developed throughout the laminate, significant delaminations were observed between all plies, fiber breaks (yarn fractures) did not occur (i.e., all inserted yarn cracklets remained unbroken).
Prediction on other cases are all in good qualitative agreement with available X-radiograph of tested specimens, detailed images and descriptions are left out here.

\section{Conclusions}\label{sec4:conclusion}

In this paper, a brief review of recent studies involving single or multiple types of damage evolution in laminates is given, which indicates that it is not only essential to model all three types of discrete failure modes with explicit fractures surfaces, but also critical to ensure accurate coupling and interaction between different types of cracks. To overcome various limitations of current methods in literature, a newly developed Auto-generated Geometry-based Discrete Model (AGDM) is presented.

The model has following features to realize explicit modeling of progressive damage in composite laminates: (i) supports all three types of failure modes, i.e., fiber/yarn fractures, transverse matrix cracks, and delaminations; (ii) contains high-fidelity representations of all cracks with real fracture surfaces in physical space; (iii) creates matching meshes at all potential crack bifurcation and intersections, which ensures correct displacement jump and accurate load transfer among coupling cracks. 
Moreover, the model also provides users with flexible options such as: deactivation of yarn fracture cracklets in a specific ply or an entire laminate; plies with different matrix crack spacing; random number of plies with arbitrary stacking sequence; controllable mesh size and profile regardless of matrix crack spacing. These features and prospective future functions will provide AGDM with significant efficiency and potential for various applications.

Unlike the existing methods that are mostly relied on user-defined modification and treatment of mesh (requires priorly-generated finite elements), the presented AGDM applies a novel method to create distinct, specially partitioned 'parts' with explicitly defined tie constraints, and thus no prior mesh or elements are needed.
The entire AGDM and all essential steps required for numerical analysis are executed in an intensively developed preprocessing code, which contains more than two-thousand lines of Python commands, and therefore may not be suitable to include here. In the future, AGDM is potentially available as an Abaqus plug-in.

While the focus of this paper is to review methodologies of current numerical models in literature and explain the necessity and advantages of the newly developed AGDM, many numerical examples are performed as well for evaluation. The first series of numerical cases consists of parametric and mesh convergence studies of an un-notched $[30/90/\text{--}30]_s$ laminate specimen made of T300/976, in which AGDM successfully captured the key phenomena for all three types of damage modes in the specimen. The predicted failure pattern of initial and final damage stages are both in excellent agreement with X-radiograph from experiments. Then eight numerical cases of other specimens with various layups are investigated, in which the AGDM predictions of progressive damage patterns and strengths of specimens are all in considerable agreement with X-radiograph and data from experiments. Nevertheless, several parameters/inputs potentially affecting the numerical results should be carefully investigated; meanwhile, functions including statistical material strengths and refined interfacial constitutive laws should be introduced in future versions of AGDM, since these are also important features for robust and accurate prediction on failure pattern and ultimate strength.   

In brief, this study demonstrates the excellent capability and potential of AGDM in high-fidelity modeling of progressive failure and damage evolution in composite laminates. 
Current and future work includes two main directions: first is on the development of AGDM utilities to support more features, as have been discussed above; second is about the specific applications that could be realized by AGDM, such as virtual testing and parametric studies in meso-scopic modeling of composite overwrapped pressure vessels (COPVs), where local winding angles/layup and prepreg yarn geometries of overwrap are highly complicated.

\section*{Acknowledgments}
This study was partially supported by research grant from National Institute of Standards and Technology (NIST) Agreement ID 70NANB14H323. The first author also acknowledges support from Cornell Graduate School.

{\footnotesize
\bibliography{mybibfile}
}

\listoffigures

\listoftables

\end{document}